\documentclass[12pt]{article}
\usepackage{amssymb}
\usepackage{amsmath}
\usepackage{amsthm}

\input pictex.tex
\newcommand{\sg}{\mathrm{sign}}

\newcommand{\precede}{\preceq}

\newcommand{\con}{\Gamma_{\preceq}^{\esp \esp c}}

\newcommand{\esp}{\hspace{0.02cm}}
\theoremstyle{definition}
 
\newtheorem{thm}{Theorem}[section]

\newtheorem{prop}[thm]{Proposition}

\newtheorem{lem}[thm]{Lemma}

\newtheorem{rem}[thm]{Remark}

\newtheorem{defn}[thm]{Definition}

\newtheorem{ex}[thm]{Example}

\newtheorem{question}[thm]{Question}

\input epsf

\date{}
\author{Andr\'es Navas}

\begin{document}

\title{On the dynamics of (left) orderable groups\\
Sur la dynamique des groupes ordonnables}
\maketitle

\vspace{-0.21cm}

\noindent{\bf Abstract.} We develop dynamical methods for studying left-orderable 
groups as well as the spaces of orderings associated to them. We give new and 
elementary proofs of theorems by Linnell (if a left-orderable group has infinitely many 
orderings, then it has uncountably many) and McCleary (the space of orderings of the 
free group is a Cantor set). We show that this last result also holds for countable 
torsion-free nilpotent groups which are not rank-one Abelian. Finally, we apply our 
methods to the case of braid groups. In particular, we show that the positive cone 
of the Dehornoy ordering is not finitely generated as a semigroup. To do this, we 
define the Conradian soul of an ordering as the maximal convex subgroup restricted 
to which the ordering is Conradian, and we elaborate on this notion.  

\vspace{0.281cm}

\noindent{\bf R\'esum\'e.} Nous d\'eveloppons des m\'ethodes dynamiques pour \'etudier 
les groupes ordonnables ainsi que leurs espaces d'ordres associ\'es. Nous donnons des 
preuves nouvelles et \'el\'ementaires de th\'eor\`emes d\^us \`a Linnell (si un groupe 
ordonnable poss\`ede une infinit\'e d'ordres, alors il poss\`ede une infinit\'e non 
d\'enombrable) et McCleary (l'espace des ordres du groupe libre est un Cantor). Nous 
montrons que ce dernier r\'esultat est valable aussi pour les groupes nilpotents 
d\'enombrables et sans torsion qui ne sont pas ab\'eliens de rang un. Finalement, 
nous appliquons nos m\'ethodes au cas des groupes de tresses. En particulier, nous 
d\'emontrons que le cone positif de l'ordre de Dehornoy n'est pas de type fini 
en tant que semi-groupe. Pour ce faire, nous d\'efinissons le noyau conradien d'un 
ordre comme \'etant le plus grand sous-groupe convexe sur lequel la relation est 
conradienne, et nous travaillons avec cette notion.

\vspace{0.381cm}

\noindent{\bf Keywords:} orderable groups, Conradian ordering, actions on the line.

\vspace{0.281cm}

\noindent{\bf Mots cl\'es:} groupes ordonnables, ordre conradien, actions sur la droite.

\vspace{0.381cm}

\noindent{\bf Subject classification AMS (2010)}: 06F15, 20F36, 20F60, 22F50,

\newpage

%%%%%%%%%%%%%%%%%%%%%%%%%%%%%%%%%%%%%%%%%%%%%%%%%%%%%%%%%%%%%%%%%%%%%%%%%%%%%%%%%%%%%%%

\tableofcontents

\newpage

\section*{Introduction}

\addcontentsline{toc}{section}{Introduction}

\vspace{0.25cm}

\hspace{0.35cm} The theory of orderable groups (that is, groups admitting a 
left-invariant total order relation) is a well developed subject in group theory 
whose starting points correspond to seminal works by Dedekind and H\"older at 
the end of the nineteenth century and the beginning of the twentieth century, 
respectively. Starting from the fifties, this theory was strongly pursued by 
several mathematical schools. Widely known modern references for all of this 
are the books \cite{botto} and \cite{koko}. (We should point out that, 
in general, this theory is presented as a particular subject of the much 
bigger one of lattice-orderable groups \cite{darnel,glass,kopi}.) In the recent 
years, the possibility of ordering many interesting groups (Thompson's 
group F \cite{nos}, braid groups \cite{pan-synt}, mapping class groups of 
punctured surfaces with boundary \cite{paper-in-l'enseignement}, fundamental 
groups of some hyperbolic 3-dimensional manifolds \cite{fourier,cal,pr2,RW}, etc), 
and the question of knowing whether some particular classes of groups can be ordered 
(higher rank lattices \cite{lif-morris,lifschits-witte,witte2}, groups with 
Kazhdan's property (T) \cite{valette,yo-ENS}, etc), have attracted the interest 
to this area of people coming from different fields in mathematics as low dimensional 
geometry and topology, combinatorial and geometric group theory, rigidity theory, 
mathematical logic, and model theory. 

Orderable groups have mostly been studied using pure algebraic methods. Nevertheless, 
the whole theory should have a natural dynamical counterpart. Indeed, an easy and 
well-known argument shows that every countable orderable group admits a faithful 
action by orientation-preserving homeomorphisms of the real line; moreover, the 
converse is true even without the countability hypothesis (see Proposition 
\ref{thorden}). Quite surprisingly, this very simple remark has not been 
exploited as it should have been, as the following examples show:

\vspace{0.1cm}

\noindent -- The first example of an orderable group which is non locally indicable is 
generally attributed to Bergman \cite{bergman} (see also \cite{tararin2}). This group is 
contained in $\widetilde{\mathrm{PSL}}(2,\mathbb{R})$, and it corresponds to the universal 
cover of the $(2,3,7)$-triangle group. Nevertheless, the fact that this group acts on 
the line and its first cohomology is trivial had been already remarked (almost twenty years 
before) by Thurston in relation to his famous stability theorem for codimension-one  
foliations \cite{thurston}.

\vspace{0.1cm}

\noindent -- A celebrated result by Dehornoy establishes that braid groups $B_n$ are 
orderable (see for instance \cite{dehornoy-libro}). However, readily soon after Dehornoy's 
work, Thurston pointed out to the mathematical community that the fact that these groups act 
faithfully on the line had been already noted by Nielsen in 1927 (see for instance the 
remark at the end of \cite{kassel}). Indeed, the geometric techniques by Nielsen allow to 
produce many (left-invariant and total) orders on $B_n$, and it turns out that one of 
them coincides with Dehornoy's ordering \cite{paper-in-l'enseignement}. We refer 
the reader to \cite{pan-synt} for a nice exposition of all of these ideas.

\vspace{0.1cm}

\noindent -- In the opposite direction, many results about the existence of invariant 
Radon measures for actions on the line are closely related to the prior algebraic 
theory of Conradian orders: see \S \ref{todo-conrad} for more explanation on this.

\vspace{0.1cm}

This work represents a systematic study of some of the aspects of the theory of orderable 
groups. This study is done preferably, though not only, from a dynamical viewpoint. In 
\S 1, we begin by revisiting some classical orderability criteria, as for instance the 
decomposition into positive and negative cones. We also recall the construction of the 
{\em space of orderings} associated to an orderable group, which corresponds to a 
(Hausdorff) topological space on which the underlying group acts naturally by 
conjugacy (or equivalently, by right multiplication). Roughly, two orderings are 
close if they coincide over large finite subsets. Although the author learned 
this idea from Ghys almost ten years ago, the first reference on this  
is Sikora's seminal work \cite{sikora} (see also \cite{pr}). The main 
issue here is to establish a relationship with a classical criterion 
of orderability due to Conrad, Fuchs, Lo\'s, and Ohnishi. This approach 
allows us, in particular, to give a short and simple proof of the known 
fact that every locally indicable group admits a left-invariant total order 
satisfying the so called {\em Conrad property} ({\em c.f.} Proposition \ref{notar}). 

In \S 2, we recall the classical dynamical criterion for orderability 
of countable groups. After elaborating a deep further on this, we use 
elementary perturbation type arguments for giving a new proof of the 
following result first established (in a different context) by 
McCleary \cite{MC}.\footnote{Added in Proof: Notice that 
Theorem A was presented as a conjecture in \cite{sikora}. 
Although it was already known, we have decided to include 
our proof here in order to illustrate our methods. Let us point 
out that Clay has recently shown that the space of orderings of 
$F_n$ contains points which are recurrent for the dynamics of the 
conjugacy action and whose orbits are dense, thus straightening 
Theorem A (see \cite{clay-lindo}). A dynamical proof of this result 
(inspired on our dynamical ideas) appears in \cite{rivas-dense}.}

\vspace{0.4cm}

\noindent{\bf Theorem A.} {\em For every integer $n \geq 2$, the space 
of orderings of the free group $F_n$ is homeomorphic to the Cantor set.}

\vspace{0.4cm}
 
Using a short argument due to Linnell \cite{linnell}, this allows 
us to answer by the affirmative a question from \cite{smith}.
 
\vspace{0.4cm} 

\noindent{\bf Corollary.} {\em If $\preceq$ is a left-invariant total order 
relation on $F_n$ (where $n \geq 2$), then the semigroup formed by the elements 
$g \in F_n$ satisfying \esp $g \succ id$ \esp is not finitely generated.}

\vspace{0.4cm}

In the general case, if the space of orderings of an orderable group is infinite, then it 
may have a very complicated structure. A quite interesting example illustrating this fact 
is given by braid groups which, according to a nice construction by Dubrovina and 
Dubrovin \cite{dub}, do admit orders that are isolated (in the corresponding space of 
orders). The rest of this work is a tentative approach for studying this type of 
phenomenon. For this, in \S 3 we revisit some classical 
properties for orders on groups. We begin by recalling H\"older's theorem 
concerning Archimedean orders ({\em c.f.} Proposition \ref{completa}) and free actions on the line 
({\em c.f.} Proposition \ref{existeordenarq}). In the same spirit, Proposition \ref{casi-obvio} 
shows (for countable groups) the equivalence of being bi-orderable and admitting 
{\em almost free} actions on the line. Very important 
for our approach is the dynamical counterpart of the Conrad property for left-invariant 
orders, namely the nonexistence of {\em crossed elements} (or {\em resilient orbits}) 
for the corresponding actions ({\em c.f.} Propositions \ref{conrad-yo} and \ref{conrad-yo2}). 
We then define the notion of {\em Conradian soul} of an order as the maximal convex 
subgroup such that the restriction of the original order to it satisfies the Conrad property. 
The pertinence of this concept is showed by providing an equivalent dynamical definition 
for countable orderable groups ({\em c.f.} Proposition \ref{reca}). Section 3 finishes 
with a little discussion on the notion of right-recurrence for orders, which has been 
introduced by Morris-Witte in his beautiful work on amenable orderable groups \cite{witte}.

In \S 4, we study of the structure of spaces of orderings for general 
orderable groups. In \S \ref{caso-conrad}, we begin by using pure algebraic 
arguments to show that, if $\preceq$ is a Conradian ordering on a group $\Gamma$, 
then $\preceq$ cannot be isolated when $\Gamma$ has infinitely many orders 
({\em c.f.} Proposition \ref{dorilita}). As a consequence we obtain the following 
result, which extends \cite[Proposition 1.7]{sikora}. For the statement, recall 
that the {\em rank} of a torsion-free Abelian group is the minimal dimension of 
a vector space over $\mathbb{Q}$ in which the group embeds.

\vspace{0.4cm}

\noindent{\bf Theorem B.} {\em The space of orderings of every (non-trivial) 
countable torsion-free nilpotent group which is not rank-one Abelian is 
homeomorphic to the Cantor set. Consequently, for each left-invariant total 
order $\preceq$ on such a group $\Gamma$, the semigroup formed by the elements 
$g \!\in\! \Gamma$ satisfying $g \succ\! id$ is not finitely generated.}

\vspace{0.4cm}

Continuing in this direction, in \S \ref{caso-trivial} we use the results of \S 3.3 to give a 
very short proof of the fact that, if a left-invariant total order $\preceq$ on a countable group 
$\Gamma$ has trivial Conradian soul, then $\preceq$ is not isolated in the space of orderings 
of $\Gamma$ ({\em c.f.} Proposition \ref{unilita}). Finally, by elaborating on the arguments 
of \S \ref{caso-conrad} and \S \ref{caso-trivial}, in \S \ref{caso-general} we give a 
slightly different (though equivalent) version of a recent result of 
Linnell.\footnote{Added in Proof: This corresponds essentially 
to \cite[Proposition 1.7]{linnell}, and is included in 
\cite{linnell-2}. Let us point out that a different proof 
covering the case of uncountable groups was subsequently 
given in \cite{crossings}.}

\vspace{0.4cm}

\noindent{\bf Theorem C.} {\em The space of orderings of a countable (orderable) 
group is either finite or contains a homeomorphic copy of the Cantor set.}

\vspace{0.4cm}

Perhaps more interesting than the statement above are the techniques involved in the 
proof, which are completely different from those of Linnell. These techniques allow us 
to identify (and partially understand) a very precise bifurcation phenomenon in some 
spaces of orderings. Indeed, if an ordering is isolated inside an infinite space of 
orderings, then its Conradian soul is non-trivial but admits only finitely many 
orderings. Thus, one can consider the finitely many associated 
orderings on the group obtained by changing the original one on the Conradian soul and 
keeping it outside (this procedure of {\em convex extension} is classical: see \S 
\ref{extension}). It appears that at least one of these new orderings is an 
accumulation point of its orbit under the action of the group ({\em c.f.} Proposition 
\ref{sandra}). For instance, for the case of Dubrovina-Dubrovin's ordering on $B_3$, 
the Conradian soul is isomorphic to $\mathbb{Z}$, which admits only two different 
orderings. It turns out that the associated ordering on $B_3$ is Dehornoy's one. Since 
the former is isolated in the space of orderings of $B_3$, this yields to the following 
result.\footnote{Added in Proof: Subsequent simpler and/or shorter proofs appear in 
\cite{DDRW} and \cite{NW} (see also \cite{hecke}).}

\vspace{0.4cm}

\noindent{\bf Theorem D.} {\em Dehornoy's ordering is an accumulation point of its orbit 
under the right action of $B_n$. (In other words, this ordering may be approximated by its 
conjugates.) Consequently, its positive cone is not finitely generated as a semigroup. 
Moreover, there exists a sequence of conjugates of Dubrovina-Dubrovin's 
ordering that converges to Dehornoy's ordering as well.}

\vspace{0.4cm}

The rough idea of the proofs of Theorems A, C, and D is that, starting from a left-invariant 
total order on a countable group, one can induce an action on the line, and from this 
action one may produce very many new order relations, except for some specific 
and well understood cases where the group structure is quite particular, 
and only finitely many orderings exist. 
Orderable groups appear in this way as a very flexible 
category despite the fact that, at first glance, it could seem very rigid because 
the underlying phase space is ordered and 1-dimensional. According to a general principle 
by Gromov \cite{gromov}, this mixture between flexibility and rigidity should contain some  
of the essence of the richness of the theory.\footnote{It is important to point out that this 
remark applies only to left-orderable groups, and not to the very interesting bi-orderable 
case: this theory remains completely out of reach of our methods. We point out, however, 
that Theorem C has no analogue in this context, since there exist bi-orderable 
groups admitting infinite but countably many bi-orderings \cite{butts}. 
Whether there is an analogue of Theorem A for bi-orderings remains as 
an open question.}

\vspace{0.1cm}

We have made an effort to make this article mostly self-contained, with the mild cost 
of having to reproduce some classical material. Several natural questions are left open. 
We hope that some of them are of genuine mathematical value and will serve as a guide 
for future research on the topic.

\vspace{0.25cm}

\noindent{\bf Acknowledgments.} It is a pleasure to thank \'E. Ghys for several 
discussions on the subject, D. Morris-Witte for drawing my attention to the 
relevant Dubrovina-Duvrovin's example, and A. Glass and A. Sikora for 
valuable remarks and corrections. It is also a pleasure to thank both 
R. Baeza for his invitation to the University of Talca (Chile), where some of the 
ideas of this article were born, and I. Liousse for her invitation to the University of 
Lille (France), where a substantial part of this paper was written out. Most of the 
contents here were presented (and some of them clarified) during a mini-course 
at the Cuernavaca Mathematical Institute (M\'exico), and I would like to thank 
A. Guillot, A. Arroyo, and A. Verjovsky, for their invitation and interest 
on the dynamical aspects of the theory of orderable groups. 

This work was funded by the PBCT/Conicyt via the 
Research Network on Low Dimensional Dynamical Systems.

%%%%%%%%%%%%%%%%%%%%%%%%%%%%%%%%%%%%%%%%%%%%%%%%%%%%%%%%%%%%%%%%%%%%%%%%%%%%%%%%%%%%%%%%%%%

\section{The space of orderings of an orderable group}
\label{basico-1}

\hspace{0.35cm} An order relation $\precede$ on a group $\Gamma$ is 
\textit{left-invariant} (resp. \textit{right-invariant}) if for all 
$g,h$ in $\Gamma$ such that $g \precede h$ one has $fg \precede fh$ 
(resp. $gf \precede hf$) for all $f \!\in\! \Gamma.$ The relation is 
\textit{bi-invariant} if it is simultaneously invariant by the left 
and by the right. To simplify, we will use the term {\em ordering} 
for referring to a left-invariant total order on a group, and we will say 
that a group $\Gamma$ is {\em orderable} (resp. {\em bi-orderable}) 
if it admits a total order which is invariant by the left (resp. 
by the right and by the left simultaneously).\footnote{Some 
authors use the term {\em orderable} for groups admitting a 
total bi-invariant order, and call {\em left orderable} 
the groups that we just call {\em orderable}.}

If $\precede$ is an order relation on a group $\Gamma$, we will say that 
$f \!\in\! \Gamma$ is \textit{positive} (resp. {\em negative}) if $f \succ id$ 
(resp. if $f \prec id$). Note that if $\precede$ is a total order relation then 
every non-trivial element is either positive or negative. Moreover, if $\precede$ 
is left-invariant and 
$P^{+} = P^{+}_{\precede}$ (resp. $P^{-}_{\precede} = P^{-}$) 
denotes the set of positive (resp. negative) elements in $\Gamma$ (sometimes 
called the positive (resp. negative) cone), then $P^{+}$ and $P^{-}$ are 
semigroups and $\Gamma$ is the disjoint union of $P^{+}, P^{-},$ and $\{id\}$. 
In fact, one can characterize the orderability in this way: a group $\Gamma$ 
is orderable if and only if it contains semigroups $P^{+}$ and $P^{-}$ such 
that $\Gamma$ is the disjoint union of them and $\{id\}$. (It suffices to 
define $\prec$ by declaring $f \prec g$ when $f^{-1}g$ belongs to $P^{+}$.) 
Moreover, $\Gamma$ is bi-orderable exactly when these semigroups may be 
taken invariant by conjugacy (that is, when they are normal subsemigroups).

\vspace{0.15cm}

\begin{ex} The category of orderable groups include torsion-free nilpotent groups, 
free groups, surface groups, etc. Another relevant example is given by braid 
groups $B_n$. Recall that the group $B_n$ has a presentation of the form
$$B_n = \langle \sigma_1,\ldots,\sigma_{n-1}\!: \quad 
\sigma_i \sigma_{i+1} \sigma_i= \sigma_{i+1} \sigma_i \sigma_{i+1}
\esp \esp \esp \esp \mbox{ for } \esp \esp \esp \esp 1 \leq i \leq n-2, 
\quad \sigma_i \sigma_j = \sigma_j \sigma_i 
\esp \esp \esp \esp \mbox{ for } \esp \esp \esp \esp |i-j| \geq 2 \rangle.$$
Following Dehornoy \cite{dehornoy-libro}, for $i \!\in\! \{1,\ldots,n-1\}$ an element 
of $B_n$ is said to be $\sigma_i$-positive if it may be written as a word of the form 
$$w_1 \sigma_i^{n_1} w_2 \sigma_i^{n_2} \cdots w_k \sigma_i^{n_k} w_{k+1},$$
where the $w_i$ are words on \esp \esp $\sigma_{i+1}^{\pm 1}, \ldots, \sigma_{n-1}^{\pm 1}$, \esp 
\esp and all the exponents $n_i$ are positive. An element in $B_n$ is said to be $\sigma$-positive 
if it is $\sigma_i$-positive for some $i \!\in\! \{1,\ldots,n-1\}$. The remarkable result 
by Dehornoy establishes that the set of $\sigma$-positive elements form the positive cone 
of a left-invariant total order $\preceq_D$ on $B_n$. We will refer to this order as the 
Dehornoy's ordering.

We remark that, for each $j \!\in\! \{2,\ldots,n\}$, the subgroup of $B_n$ generated by \esp 
$\sigma_{j}, \sigma_{j+1}, \ldots, \sigma_{n-1}$ \esp is naturally isomorphic to \esp 
$B_{n-j+1}$ \esp by an isomorphism which respects the corresponding Dehornoy's orderings. 
\label{dehor1}
\end{ex}

\begin{rem}
The characterization of orderings in terms of positive and negative cones shows immediately the 
following: if $\precede$ is an ordering on a group $\Gamma$, then the order $\bar{\precede}$ 
defined by \esp $g \esp\esp \bar{\succ} \esp\esp id$ \esp if and only if \esp $g \prec id$ 
\esp is also left-invariant and total. 
\label{la-primera}
\end{rem}

\vspace{0.1cm}

Given an orderable group $\Gamma$ we denote by $\mathcal{O}(\Gamma)$ the 
set of all the orderings on $\Gamma$. As it was pointed out to the author 
by Ghys, the group $\Gamma$ acts on $\mathcal{O}(\Gamma)$ by conjugacy 
(or equivalently, by right multiplication): given an order $\precede$ with 
positive cone $P^+$ and an element $f \!\in\! \Gamma$, the image of $\precede$ 
under $f$ is the order $\precede_{f}$ whose positive cone is \esp $f \esp P^+ f^{-1}$. 
\esp In other words, one has \esp $g \precede_f h$ \esp if and only if \esp $fgf^{-1} 
\preceq fhf^{-1}$, \esp which is equivalent to \esp $gf^{-1} \precede hf^{-1}$. 

\begin{rem} If $\Gamma$ is an orderable group, then the whole group of automorphisms of $\Gamma$ 
(and not only the conjugacies) acts on $\mathcal{O}(\Gamma)$. This may be useful for studying 
bi-orderable groups. Indeed, since the fixed points for the right action of $\Gamma$ on 
$\mathcal{O}(\Gamma)$ correspond to the bi-invariant orderings, the group of outer  
automorphisms of $\Gamma$ acts on the corresponding {\em space of bi-orderings}.
\end{rem}

\vspace{0.1cm}

The {\em space of orderings} \esp $\mathcal{O}(\Gamma)$ \esp has a natural (Hausdorff) 
topology first introduced (and exploited) by Sikora in \cite{sikora}. A sub-basis 
of this topology is the family of the sets of the form 
\esp $U_{f,g} \!=\! \{\precede : \esp f \!\prec\! g\}$.  
\esp Note that the right action of $\Gamma$ on $\mathcal{O}(\Gamma)$ 
becomes in this way an action by homeomorphisms. Similarly, the map 
sending $\precede$ to $\bar{\precede}$ from Example \ref{la-primera} is a 
continuous involution of $\mathcal{O}(\Gamma)$. To understand the topology on 
$\mathcal{O}(\Gamma)$ better, associated to the symbols $-$ and $+$ let us consider 
the space \esp $\{-,+\}^{\Gamma \setminus \{id\}}$. \esp  We claim that there exists 
a one-to-one correspondence between the set $\mathcal{O}(\Gamma)$ and the subset 
$\mathcal{X}(\Gamma)$ of $\{-,+\}^{\Gamma \setminus \{id\}}$ formed by the functions 
\esp \esp $\sg \!: \Gamma \setminus \{id\} \rightarrow \{-,+\}$ \esp \esp satisfying: 

\vspace{0.1cm}

\noindent -- for every $g \in \Gamma \setminus \{id\}$ one has \esp $\sg (g) \neq \sg (g^{-1})$,

\vspace{0.1cm} 

\noindent -- if $f,g$ in $\Gamma \setminus \{id\}$ are such that \esp $\sg (f) = \sg (g)$, 
\esp then \esp $\sg (fg) = \sg (f) = \sg (g)$.

\vspace{0.1cm}

Indeed, to each $\preceq$ in $\mathcal{O}(\Gamma)$ we may associate the function \esp 
$\sg_{\preceq}\!\!: \Gamma \setminus \{id\} \rightarrow \{-,+\}$ \esp defined by \esp $\sg_{\preceq}(g)=+$ 
\esp if and only if $g \succ id$. Conversely, given a function $\sg$ with the properties above, we may 
associate to it the unique order $\preceq_{\sg}$ in $\mathcal{O}(\Gamma)$ which satisfies 
$f \succ_{\sg} g$ if and only if \esp $\sg (g^{-1} f)$ \esp equals $+$. Now 
if we endow $\{-,+\}^{\Gamma \setminus \{id\}}$ with the product topology and $\mathcal{X} (\Gamma)$ 
with the subspace one, then the induced topology on $\mathcal{O}(\Gamma)$ via the preceding 
identification coincides with the topology previously defined by prescribing the sub-basis elements. 
As a consequence, since $\{-,+\}^{\Gamma \setminus \{id\}}$ is compact and $\mathcal{X} (\Gamma)$ 
is closed therein, this shows that the topological space $\mathcal{O}(\Gamma)$ is always 
compact.

The compactness of $\mathcal{O}(\Gamma)$ is by no means a new result. It was first established 
for countable groups by Sikora \cite{sikora}. Subsequent proofs covering the case of uncountable 
groups appear in \cite{pr} and \cite{witte}. Although our approach is not the simplest possible 
one, it allows us revisiting some classical orderability criteria essentially due to Conrad, 
Fuchs, Lo\'s, and Ohnishi (see for instance \cite{botto,glass,koko}). 
This is summarized in Proposition \ref{bur} below. For the statement, 
let us consider the following two conditions: 

\vspace{0.1cm}

\noindent (i) \esp \esp \esp For every finite family of elements $g_1,\ldots,g_k$ which are 
different from the identity, there exists a family of exponents $\eta_i \in \! \{-1,1\}$ such that 
\esp $id$ \esp does not belong to the semigroup generated by the elements of the form $g_i^{\eta_i}$,

\vspace{0.1cm}

\noindent (ii) \esp \esp \esp For every finite family of elements $g_1,\ldots,g_k$ which are different 
from the identity, there exists a family of exponents $\eta_i \in \! \{-1,1\}$ such that \esp $id$ \esp 
does not belong to the smallest semigroup which simultaneously satisfies the following two properties:\\ 

\noindent -- it contains all the elements \esp $g_i^{\eta_i}$; 

\noindent -- for all $f,g$ in the semigroup, the elements \esp 
$fgf^{-1}$ \esp and \esp $f^{-1}gf$ \esp also belong to it.

\vspace{0.1cm}

\noindent In each case such a choice of the exponents $\eta_i$ will be said to be {\em compatible}.

\vspace{0.1cm}

\begin{prop} {\em A group $\Gamma$ is orderable (resp. bi-orderable) if and only if it 
satisfies condition $\mathrm{(i)}$ (resp. condition $\mathrm{(ii)}$) above.}
\label{bur}
\end{prop}

\noindent{\bf Proof.} The necessity of the conditions (i) or (ii) is clear: it suffices 
to chose each exponent $\eta_i$ so that $g_i^{\eta_i}$ becomes a positive element. 

To prove the converse claim in case (i), for each finite family $g_1,\ldots,g_k$ of elements 
in $\Gamma$ which are different from the identity, and for each compatible choice of exponents 
$\eta_i \!\in\! \{-1,1\}$, let us consider the (closed) subset 
$\mathcal{X}(g_1,\ldots,g_k;\eta_1,\ldots,\eta_k)$ 
of $\{-,+\}^{\Gamma \setminus \{id\}}$ formed by all of the \esp $\sg$ functions \esp which 
satisfy the following property: one has \esp $\sg (g) = +$ \esp and $\sg (g^{-1}) = -$ 
\esp for every $g$ belonging to the semigroup generated by the elements $g_i^{\eta_i}$. 
(It easily follows from the hypothesis that this subset is non-empty.) Now for fixed 
$g_1,\ldots,g_k$ let $\mathcal{X}(g_1,\ldots,g_k)$ be the union of all the sets of the form 
$\mathcal{X}(g_1,\ldots,g_k;\eta_1,\ldots,\eta_k)$, where the choice of the exponents 
$\eta_i$ is compatible. Note that, if \esp $\{\mathcal{X}_i = \mathcal{X}(g_{i,1},\ldots,g_{i,k_i}), 
\esp \esp i \!\in\! \{1,\ldots,n\} \}$ \esp is a finite family of subsets of this 
form, then the intersection \esp $\mathcal{X}_1 \cap \ldots \cap \mathcal{X}_n$ 
\esp contains the (non-empty) set \esp 
$\mathcal{X}(g_{1,1},\ldots,g_{1,k_1},\ldots,g_{n,1},\ldots,g_{n,k_n})$, 
\esp and it is therefore non-empty. Since $\{-,+\}^{\Gamma \setminus \{id\}}$ is compact, 
a direct application of the Finite Intersection Property shows that the intersection 
$\mathcal{X}$ of all the sets of the form $\mathcal{X}(g_1,\ldots,g_k)$ is (closed and) 
non-empty. It is quite clear that $\mathcal{X}$ is actually contained in 
$\mathcal{X}(\Gamma)$, and this shows that $\Gamma$ is orderable.

The case of condition (ii) is similar. We just need to replace 
the sets $\mathcal{X}(g_1,\ldots,g_k;\eta_1,\ldots,\eta_k)$ by the 
sets $B \mathcal{X}(g_1\ldots,g_k;\eta_1,\ldots,\eta_k)$ formed by all of 
the \esp $\sg$ \esp functions satisfying \esp $\sg (g) \!=\! +$ \esp 
and \esp $\sg (g^{-1}) \!=\! -$ \esp for every $g$ belonging 
to the smallest semigroup satisfying simultaneously the following properties: 

\noindent -- it contains all of the elements $g_i^{\eta_i}$; 

\noindent -- for every $f,g$ in the semigroup, the elements \esp $fgf^{-1}$ 
\esp and \esp $f^{-1}gf$ \esp also belong to it. $\hfill\square$

\vspace{0.4cm}

What is relevant with the previous conditions (i) and (ii) is that they involve only finitely 
many elements. This shows in particular that the properties of being orderable or bi-orderable 
are ``local", that is, if they are satisfied by every finitely generated subgroup of 
a group $\Gamma$, then they are satisfied by $\Gamma$ itself. As we have already mentioned, 
all these facts are well-known. The classical proofs use the Axiom of Choice, and our 
approach just uses its topological equivalent, 
namely Tychonov's theorem. This point of view is more appropriate in relation to spaces 
of orderings. It will be used once again when dealing with Conradian orders, and it 
will serve to justify the pertinence of Question \ref{recurrente-local}.

\vspace{0.25cm}

If $\Gamma$ is a countable orderable group, then the topology 
on $\mathcal{O}(\Gamma)$ is metrizable. Indeed, if $\mathcal{G}_0 \subset \mathcal{G}_1 
\subset \ldots$ is a complete exhaustion of $\Gamma$ by finite sets, then we can define 
the distance between two different orderings \esp $\leq$ \esp and \esp $\precede$ \esp by 
letting \esp $dist(\leq,\precede) = e^{-n}$, \esp where $n$ is the maximum non negative 
integer number such that $\leq$ and $\precede$ coincide on $\mathcal{G}_n$. An equivalent 
metric \esp $dist'$ \esp is obtained by letting \esp $dist' (\leq,\preceq) = e^{-n'}$, 
\esp where $n'$ is the maximum non negative integer such that the positive cones 
of $\leq$ and $\preceq$ coincide on $\mathcal{G}_{n'}$, that is, \esp 
$P_{\leq} \cap \mathcal{G}_{n'} = P_{\preceq} \cap \mathcal{G}_{n'}$. \esp One 
easily checks that these metrics are ultrametric. Moreover, the fact that 
$\mathcal{O}(\Gamma)$ is compact becomes more transparent in this case.

\vspace{0.1cm}
 
When $\Gamma$ is finitely generated, one may choose $\mathcal{G}_n$ as being the ball 
of radius $n$ with respect to some finite and symmetric system of generators 
$\mathcal{G}$ of $\Gamma$, that is, the set of elements $g$ which can be written in the 
form \esp $g = g_{i_1} g_{i_2} \cdots g_{i_m}$, \esp where $g_{i_j} \in \mathcal{G}$ and 
\esp $0 \leq m \leq n$. (In this case the action of $\Gamma$ on $\mathcal{O}(\Gamma)$ is by 
bi-Lipschitz homeomorphisms.) One easily checks that the metrics on $\mathcal{O}(\Gamma)$ 
resulting from two different finite systems of generators are not only topologically 
equivalent but also H\"older equivalent. Therefore, according to Theorem A, the 
following question (suggested to the author by L. Flaminio) makes sense.

\vspace{0.1cm}

\begin{question} What can be said about the  metric structure (up to Lipschitz equivalence) 
of the Cantor set viewed as the 
space of orderings of the free groups $F_n$~? For instance, are the corresponding Hausdorff 
dimensions positive and finite~? If so, what can be said about the supremum or the infimum 
value of the Hausdorff dimensions when ranging over all finite systems of generators~? (Note 
that using the arguments of \cite{sikora}, one can easily show that the Hausdorff dimension 
of $\mathcal{O}(\mathbb{Z}^n)$ is equal to zero.)
\label{lip}
\end{question}

\vspace{0.1cm} 

In general, the study of the dynamics of the action of $\Gamma$ on $\mathcal{O}(\Gamma)$ should 
reveal useful information. This is indeed the main idea behind the proof of Morris-Witte's 
theorem \cite{witte}: see \S \ref{a-derecha}. Let us formulate two simple questions on 
this (see also Question \ref{probability}).

\vspace{0.1cm}

\begin{question} For which countable orderable groups the action of $\Gamma$ on 
$\mathcal{O}(\Gamma)$ is uniformly equicontinuous~? The same question makes 
sense for topological transitivity, or for having a dense orbit.
\end{question}

\vspace{0.01cm}

\begin{question} What can be said in general about the space $\mathcal{O}(\Gamma) / \Gamma$~? 
For instance, is the set of isolated orderings modulo the right action of $\Gamma$ always 
finite~? (Compare \cite[Theorem 3.5]{paper-in-l'enseignement}.)
\label{tuti}
\end{question}

\vspace{0.1cm}

To close this Section, we recall a short argument due to Linnell \cite{linnell} showing 
that if an ordering $\preceq$ on a group $\Gamma$ is non isolated in $\mathcal{O}(\Gamma)$, 
then its positive cone is not finitely generated as a semigroup. This shows why the 
Corollary in the Introduction of this work follows directly from Theorem A.

\vspace{0.15cm}

\begin{prop} {\em If $\preceq$ is a left-invariant total order on a group $\Gamma$ 
and $\preceq$ is non isolated in $\mathcal{O}(\Gamma)$, then the corresponding 
positive cone is not finitely generated as a semigroup.}
\label{paja}
\end{prop}

\noindent{\bf Proof.} If \esp $g_1,\ldots,g_k$ \esp generate $P_{\preceq}^{+}$, then 
the only ordering on $\Gamma$ which coincides with $\preceq$ on any set containing 
these generators and the identity element is $\preceq$ itself... $\hfill\square$

%%%%%%%%%%%%%%%%%%%%%%%%%%%%%%%%%%%%%%%%%%%%%%%%%%%%%%%%%%%%%%%%%%%%%%%%%%%%%%%%%%%%%%%%%%%%%%%%%%%%%%
%%%%%%%%%%%%%%%%%%%%%%%%%%%%%%%%%%%%%%%%%%%%%%%%%%%%%%%%%%%%%%%%%%%%%%%%%%%%%%%%%%%%%%%%%%%%%%%%%%%%%%

\section{The dynamical realization of countable orderable groups}

\subsection{A dynamical criterion for orderability}
\label{basico}

\hspace{0.35cm} The following dynamical criterion for group orderability is 
classical. We refer to \cite{ghys} for more details (see also \cite{stein} 
for an extension to the case of partially ordered groups). 

\vspace{0.12cm}

\begin{prop} \textit{For every countable group $\Gamma$, the following properties are equivalent:}

\vspace{0.05cm}

\noindent{(i) {\em $\Gamma$ acts faithfully on the real line by 
orientation-preserving homeomorphisms,}}

\vspace{0.05cm}

\noindent{(ii) {\em $\Gamma$ is an orderable group.}}
\label{thorden}
\end{prop}

\noindent{\bf Proof.} Assume that $\Gamma$ acts faithfully by 
orientation-preserving homeomorphisms of the line. Let us consider a dense sequence 
$(x_n)$ in $\mathbb{R}$, and let us define $g \prec h$ if for the smallest 
index $n$  such that $g(x_n) \neq h(x_n)$ one has $g(x_n) < h(x_n).$ One 
easily checks that $\precede$ is a total left-invariant order relation. 
(Note that this direction does not use the countability hypothesis.)

Suppose now that $\Gamma$ admits a left-invariant total order $\precede$. Choose 
a numbering $(g_i)_{i \geq 0}$ for the elements of $\Gamma$, put $t(g_0) \! = \! 0$, 
and define $t(g_k)$ by induction in the following way: assuming that 
$t(g_0),\ldots,t(g_i)$ have been already defined, if $g_{i+1}$ is bigger 
(resp. smaller) than $g_0,\ldots,g_i$ then put $t(g_{i+1}) = \mathrm{max} \{
t(g_0),\ldots,t(g_i) \} + 1$ (resp. $\mathrm{min} \{
t(g_0),\ldots,t(g_i) \} - 1$), and if $g_m \prec g_{i+1}
\prec g_n$ for some $m,n$ in $\{0,\ldots,i\}$ and $g_j$ 
is not between $g_m$ and $g_n$ for any $0 \leq j \leq i$ 
then let $t(g_{i+1})$ be equal to $(t(g_m) + t(g_n))/2$.

Note that $\Gamma$ acts naturally on $t(\Gamma)$ by $g (t(g_i)) = t(gg_i)$. 
It is not difficult to see that this action extends continuously to the closure 
of the set $t(\Gamma)$. (Compare Lemma \ref{semi-igual}.) Finally, one can extend 
the action to the whole line by extending the maps $g$ affinely to each interval 
of the complementary set of $t(\Gamma)$. $\hfill\square$

\vspace{0.45cm}

It is worth analyzing the preceding proof carefully. If $\precede$ is an 
ordering on a countable group $\Gamma$ and $(g_i)_{i \geq 0}$ is a numbering 
of the elements of $\Gamma$, then we will call the (associated) {\em dynamical 
realization} the action of $\Gamma$ on $\mathbb{R}$ constructed in this proof. 
It is easy to see that this realization has no global fixed point unless $\Gamma$  
is trivial. Moreover, if $f$ is an element of $\Gamma$ whose dynamical realization 
has two fixed points $a\!<\!b$ (which may be equal to $\pm \infty$) and has no fixed 
point in $]a,b[$, then there must exist some point of the form $t(g)$ inside $]a,b[$. 
Finally, it is not difficult to show that the dynamical realizations associated to 
different numberings of the elements of $\Gamma$ are all topologically conjugate. 
(Compare Lemma \ref{semi-igual}.) Therefore, we can speak of any dynamical property 
for the dynamical realization without referring to a particular numbering. 

More interesting is to analyze the order obtained from an action on the line. First, note that 
if the dense sequence \esp $(x_n)$ \esp is such that the orbit of the first point $x_0$ is 
{\em free} (that is, one has \esp $g (x_0) \neq x_0$ \esp for all \esp $g \neq id$), \esp 
then the tail $(x_n)_{n \geq 1}$ of the sequence is irrelevant for the definition of the 
associated order. This remark is non innocuous since many group actions on the line have 
free orbits, as the following examples show.

\vspace{0.04cm}

\begin{ex} Let $\Gamma$ be the affine group over the rationals (that is, the group of maps 
of the form \esp $x \mapsto bx + a$, \esp where $a,b$ belong to $\mathbb{Q}$). Clearly, the 
orbit of every irrational number \esp $\varepsilon$ \esp by the natural action of $\Gamma$ 
on the line is free. Therefore, we may define an ordering $\preceq_{\varepsilon}$ on 
$\Gamma$ by declaring that \esp \esp $g \succ_{\varepsilon} id$ \esp \esp if and only if 
\esp $g(1/\varepsilon) > 1/\varepsilon$. \esp Note that for \esp $g(x) = bx + a$, \esp this 
is equivalent to \esp $b + \varepsilon a > 1$. \esp The orderings $\preceq_{\varepsilon}$ 
were introduced by Smirnov in \cite{smirnov}.
\label{ex-smirnov}
\end{ex}

\begin{ex} As it is well explained in \cite{paper-in-l'enseignement}, the actions of 
braid groups on the line constructed using Nielsen's geometrical arguments have 
(plenty of) free orbits.
\end{ex}

\vspace{0.04cm}

Perhaps the most important (and somehow ``universal") case of actions with free orbits corresponds 
to dynamical realizations of left-invariant total orders $\preceq$ on countable groups: the orbit 
of the point \esp $t(id)$ \esp --and therefore the orbit of each point of the form $t(h)$-- \esp 
is free, since \esp $g(t(id)) = t(g) \neq t(id)$ \esp for every \esp $g \neq id$. \esp 

The existence of free orbits allows showing that 
not all actions without global fixed points of (countable) orderable groups appear as 
dynamical realizations. For instance, this is the case of non-Abelian groups of piecewise-linear 
homeomorphisms of the line which coincide with translations outside a compact subset, 
as for example Thompson's group F (see \cite{brin}). 
Indeed, non-trivial commutators in such a group have 
intervals of fixed points; by suitable conjugacies, the intervals so obtained cover 
the line, hence no point has free orbit.

\vspace{0.08cm}

\begin{question} What are the (countable) orderable groups all of whose actions by 
orientation-preserving homeomorphisms of the line without global fixed points 
are semiconjugate to dynamical realizations~? (For example, this is the case of the 
group \esp $(\mathbb{Z},+)$.)
\label{remark-smirnov}
\end{question}

\begin{question} For countable orderable groups, what can be said on the structure of 
the space of faithful actions on the line up to topological semiconjugacy\hspace{0.05cm}? 
(Compare Question \ref{tuti}.)
\end{question}

\vspace{0.08cm}

Remark that, for each $g \!\in\! \Gamma$, the order relation for which an element \esp 
$h \!\in\! \Gamma$ \esp is positive if and only if \esp $g (t(h)) \!>\! t(h)$ \esp is 
no other thing than the conjugate of $\preceq$ by $h^{-1}$. Indeed, by construction, 
the condition $g(t(h)) > t(h)$ is equivalent to $t(gh) > t(h)$, \esp and therefore to \esp 
$g h \succ h$, \esp that is, to \esp \esp $h^{-1} g h \succ id$. \esp \esp Letting $h = id$, 
this allows to recover the original ordering $\preceq$ from its dynamical realization. 

\vspace{0.1cm}

\begin{rem} The involution \esp \esp $\preceq \hspace{0.15cm} \mapsto \bar{\preceq}$ \esp \esp 
of \esp $\mathcal{O}(\Gamma)$ \esp introduced in Remark \ref{la-primera} has also a dynamical 
interpretation. Indeed, let $\Gamma$ be a group of orientation-preserving homeomorphisms 
of the line, and let $(x_n)$ be a dense sequence of points in $\mathbb{R}$. If $\preceq$ 
is the order on $\Gamma$ induced from this sequence and 
$\varphi\!\!: \mathbb{R} \rightarrow \mathbb{R}$ is an orientation-reversing 
homeomorphism, then the order on $\Gamma$ induced by the dense sequence 
$(\varphi(x_n))$ and the action \esp $g \mapsto \varphi \circ g \circ \varphi^{-1}$ 
\esp corresponds to $\bar{\precede}$.
\end{rem}

\vspace{0.1cm}

In general, the homeomorphisms appearing in dynamical realizations are 
not smooth. However, according to \cite[Th\'eor\`eme D]{DKN}, the dynamical 
realization of every countable orderable group is topologically conjugate 
to a group of locally Lipschitz homeomorphisms of the line.

Although faithful actions on the line contain all the algebraic information of 
the corresponding orderable group, these actions are not always easy to deal with. 
For instance, according to \cite[Proposition 5.7]{DKN}, for a countable orderable 
group $\Gamma$, none of its actions on the line provides relevant probabilistic 
information when the initial distribution is symmetric (see however \cite{kaimanovich} 
for some interesting examples in the non symmetric case; see also \cite{pickel}). 
Nevertheless, a probabilistic approach may be useful for the study of the action 
of $\Gamma$ on $\mathcal{O} (\Gamma)$. A basic question on this is the following.
\vspace{0.1cm}

\begin{question} If $\Gamma$ is a countable group having infinitely 
many left-invariant total orders, under what conditions is the space 
$\mathcal{O} (\Gamma)$ a $\Gamma$-boundary (in the sense of \cite{furman})~? 
\label{probability}
\end{question} 

%%%%%%%%%%%%%%%%%%%%%%%%%%%%%%%%%%%%%%%%%%%%%%%%%%%%%%%%%%%%%%%%%%%%%%%%%%%%%%%%%%%%%%%%%%%%%%%%%%%%

\subsection{On the space of orderings of free groups}

\hspace{0.35cm} A natural strategy for proving Theorem A is the following. Starting 
with an ordering on the free group $F_n$, one considers the corresponding dynamical 
realization. By slightly perturbing the homeomorphisms corresponding to a system of 
free generators of $F_n$, one obtains an action on the line of a group which ``in 
most cases'' will still be free \cite[Proposition 4.5]{ghys}. From the perturbed  
action one may induce a new ordering on $F_n$, which will be near the original 
one if the perturbation is very small (with respect to the compact-open 
topology). Finally, in general this new ordering should be different, 
because if not then the original action would be ``structurally stable'', 
and this cannot be the case for free group actions on the line.

To put all these ideas in practice there are some technical difficulties. Although the strategy 
that we will actually follow uses a similar idea, it does not rely on any genericity type 
argument. This will allow us to provide an elementary and self-contained  proof for Theorem A.

\vspace{0.05cm}

Recall that given two faithful actions \esp $\phi_i \!: \Gamma \rightarrow \mathrm{Homeo}_+(\mathbb{R})$, 
\esp $i \!\!\in\!\! \{1,2\}$, \esp the action $\phi_2$ is said to be {\em topologically semiconjugate} 
to $\phi_1$ if there exists a continuous non-decreasing surjective map \esp $\varphi\!: \mathbb{R} 
\rightarrow \mathbb{R}$ \esp such that \esp $\phi_1 (g) \circ \varphi = \varphi \circ \phi_2 (g)$ 
\esp for every \esp $g \!\in\! \Gamma$. \esp The following criterion will allow us to distinguish 
two orderings obtained from actions on the line.

\vspace{0.1cm}

\begin{lem} {\em Let $\preceq$ be an ordering on a non-trivial countable group \esp $\Gamma$, and 
let $\phi_1$ be the action corresponding to a dynamical realization of $\preceq$. Let $\phi_2$ 
be an action of \esp $\Gamma$ by orientation-preserving homeomorphisms of the line for which 
there is no global fixed point and such that the orbit of the origin is free. If $\preceq'$ 
denotes the ordering on $\Gamma$ induced from the $\phi_2$-orbit of the origin, then $\preceq$ 
and $\preceq'$ coincide if and only if $\phi_2$ is topologically semiconjugate to $\phi_1$.}
\label{semi-igual}
\end{lem}

\noindent{\bf Proof.} If $\phi_2$ is topologically semiconjugate to $\phi_1$, then the relative 
positions of the points in \esp $\{\phi_i (g), \esp \esp g \!\in\! \Gamma\}$ \esp are the same 
for \esp $i\!=\!1$ \esp and \esp $i\!=\!2$. \esp From this one easily concludes that the induced 
orderings $\preceq$ and $\preceq'$ coincide.

Conversely, if $\preceq$ and $\preceq'$ coincide, then we may define a map $\varphi$ from the 
$\phi_2$-orbit of the origin to the set $t(\Gamma)$ by sending $\phi_2 (g)(0)$ to $t(g)\!=\!\phi_1(g)(0)$. 
This map $\varphi$ is strictly increasing because both conditions \esp $\phi_2(g)(0) > \phi_2(h)(0)$ 
\esp and \esp $t(g) > t(h)$ \esp are equivalent to \esp $g \succ h$. \esp Moreover, $\varphi$ 
satisfies \esp $\phi_1(g) \circ \varphi = \varphi \circ \phi_2(g)$ \esp for every $g \! \in 
\!\Gamma$. 

\vspace{0.25cm}

\noindent{\underbar{Claim.}} The map $\varphi$ extends continuously to a 
non-decreasing map defined on the closure of the $\phi_2$-orbit of the origin. 

\vspace{0.1cm}

Indeed, to show that $\varphi$ has a continuous extension to the closure, 
it suffices to show that, if two sequences $(g_n)$, $(h_n)$ of elements of 
$\Gamma$, the first of which being strictly increasing and the second strictly 
decreasing, are such that $\lim_n \phi_2 (g_n) (0) = p = \lim_n \phi_2 (h_n) (0)$, then 
the points $a = \lim_n t (g_n)$ and $b = \lim_n t(h_n)$ coincide. Suppose not, and let 
$\varepsilon = b-a$. Let $n \in \mathbb{N}$ be such that $t(h_n) - b < \varepsilon/ 3$ and 
$a - t(g_n) < \varepsilon/3$. Since for each $n$ there exist elements between $g_n$ and $h_n$, 
the method of construction of the dynamical realization implies that the midpoint between 
$t(g_n)$ and $t(h_n)$ must belong to $t(\Gamma)$. By the definition of $\varepsilon$, 
this midpoint $t(f_1)$ belongs to $]a,b[$. Similarly, the midpoint of between $t(f_1)$ 
and $t(h_n)$ belongs to $]a,b[ \esp\cap\esp t(\Gamma)$, thus it is of the form 
$t(f_2)$ for some $f_2 \in \Gamma$.

Now, let $f \in \Gamma$ be any element such that $t (f) \in \esp ]a,b[$. We have, 
$t(g_n) < t(f) < t(h_n)$, hence $g_n \prec f \prec h_n$, for all $n \in \mathbb{N}$. 
As a consequence, $\phi_2 (g_n)(0) < \phi(f)(0) < \phi_2(h_n)(0)$. Passing to the limit 
this yields $\phi(f) = p$. Applying this to the elements $f_1 \neq f_2$, we obtain 
$\phi_2 (f_1) (0) = \phi_2 (f_2) (0) = p$. However, this contradicts the fact that 
the $\phi_2$-orbit of the origin is free. Thus, $\varphi$ extends continuously, 
and since it is strictly increasing when defined on $\phi_2(\Gamma)(0)$, 
its extension to the closure of this set is non-decreasing. 

\vspace{0.2cm}

Now notice that, if $t(\Gamma)$ is dense in the line, then there is only one way to 
extending $\varphi$ into a non-decreasing continuous and surjective map realizing the semiconjugacy. 
If not, let $]a,b[$ be a connected component of the complementary set of the closure of $t(\Gamma)$. 
Choosing an arbitrary orientation-preserving homeomorphism between the intervals \esp 
$[\varphi^{-1}(c),\varphi^{-1}(d)]$ \esp and \esp $[c,d]$, \esp and extending it to the 
orbits by $\Gamma$ of these intervals in an equivariant way, 
we may enlarge the domain of definition of $\varphi$ 
still preserving the semiconjugacy relation \esp $\phi_1(g) \circ \varphi = \varphi \circ \phi_2(g)$. 
\esp Doing this with all the connected components of the complementary set of the closure of \esp 
$t(\Gamma)$, \esp we can extend $\varphi$ to a semiconjugacy from $\phi_2$ to $\phi_1$ defined 
on the whole real line. $\hfill\square$
 
\vspace{0.45cm}

During the proof of Theorem A, we will need to approximate a given homeomorphism of the interval 
by a real-analytic one. Although there exist many results of this type for general compact 
manifolds with boundary, the one-dimensional version of this fact is elementary. 

\vspace{0.1cm}

\begin{lem} \label{obvio} 
{\em Every orientation-preserving homeomorphism of the interval $[0,1]$ can be approximated 
(in the $\sup$-norm) by a sequence of real-analytic orientation-preserving diffeomorphisms.}
\end{lem}

\noindent{\bf Proof.} Let $f$ be an orientation-preserving 
homeomorphism of $[0,1]$. For each $n \in \mathbb{N}$ let $f_n$ be 
a $C^1$ diffeomorphism sending the point $i/n$ into $f(i/n)$, for all
$i \in \{0,1,\ldots,n\}$. Such an $f_n$ can be easily constructed by using an interpolation 
method. Alternatively, one may use piecewise-linear homeomorphisms, and then smoothing the 
derivative at the break-points by conjugating with (a translate of) the map 
$x \mapsto \exp (-1 / x)$ (see \cite{Ts}).  

Now, for each $n \!\in\! \mathbb{N}$, let us consider the derivative 
$f_n' \!:  [0,1] \to \mathbb{R}$ of $f_n$. This is a continuous function 
satisfying $f_n'(x) \geq \lambda_n$ for some $\lambda_n > 0$ and all $x \in [0,1]$. 
By the Stone-Weierstrass Theorem, each $f_n'$ can be approximated by a 
sequence of real-analytic functions (even polynomials) $h_{n,k}$. For 
$k$ large enough we have $\big| g_{n,k} (x) - f_n'(x) \big| 
\leq \min \{ 1/n, \lambda_n / 2\}$ for all $x \in [0,1]$. 
We choose such a $k = k_n$, and we let $g_n = g_{n,k_n}$. 

By integrating $g_n$, we obtain a diffeomorphism $F_n$ from $[0,1]$ to a 
certain interval $[0,y_n]$. Since $g_n$ and $f_n'$ are close and $y_n$ is 
the total integral of $g_n$, the sequence $(y_n)$ converges to $1$. Thus, 
by rescaling the image of each $F_n$, we get the desired sequence of 
real-analytic diffeomorphisms approximating $f$. $\hfill\square$

\vspace{0.35cm}

We can now proceed to the proof of Theorem A. Let $\preceq$ be an ordering on the free group 
$F_n$. Given an arbitrary finite family of positive elements \esp $h_j \!\in\! F_n$, \esp 
where \esp $j \!\in\! \{1,\ldots,m\}$, \esp we need to show the existence of a distinct 
ordering $\preceq'$ on $F_n$ for which all of these elements are still positive. To do this,  
let us fix a free system of generators \esp $\{g_1,\ldots,g_n\}$ \esp of $F_n$. Let us also 
consider the corresponding generators \esp $g_{1,0},\ldots,g_{n,0}$ \esp of a dynamical 
realization of $\preceq$ associated to a numbering of the elements of $F_n$ starting with 
$id$. We first claim that, given \esp $i \!\in\! \{1,\ldots,n\}$, \esp there exists a sequence 
of real-analytic diffeomorphisms \esp $g_{i,k} \!\in\! \mathrm{Homeo}_{+}(\mathbb{R})$ \esp 
that converges to $g_{i,0}$ in the compact-open topology and such that, for each fixed $k$, 
the group $\Gamma_k$ generated by $g_{1,k},\ldots,g_{n,k}$ has no global fixed point. 
Indeed, let us fix a real-analytic diffeomorphism 
\esp $\varphi\!: \mathbb{R} \rightarrow \hspace{0.03cm} ]0,1[$. \esp 
By Lemma \ref{obvio}, the conjugate homeomorphisms 
\esp $\bar{g}_{i,0} = \varphi \circ g_{i,0} \circ \varphi^{-1}$, \esp 
$i \!\in\! \{1,\ldots,n\}$, \esp may be approximated in the strong topology on $[0,1]$ 
by sequences of real-analytical diffeomorphisms \esp $\bar{g}_{i,k}$ \esp of 
$[0,1]$. This easily implies that each $g_{i,0}$ may be approximated in the 
compact-open topology by the sequence of real-analytic diffeomorphisms \esp 
$g_{i,k} = \varphi^{-1} \circ \bar{g}_{i,k} \circ \varphi$. \esp Finally, by conjugating 
each of these maps by a very small translation \esp $T_{i,k}$, \esp we may assume 
that for each fixed \esp $k \! \in \! \mathbb{N}$ \esp the maps \esp $g_{i,k}$  
\esp have no common fixed point, and therefore the group $\Gamma_k$ 
generated by them has no global fixed point in the line. 

\vspace{0.25cm}

\noindent{\underbar{Case 1.}} Passing to a subsequence if necessary, for every $k$ the elements 
\hspace{0.03cm} $g_{1,k},\ldots,g_{n,k}$ \hspace{0.03cm} satisfy some non-trivial relation.

\vspace{0.15cm}

In this case \esp $\Gamma_k \!\sim\! F_n / N_k$ \esp for some non-trivial normal subgroup 
$N_k$ in $F_n$. Let us write one of the elements $h_j$ above as a product of the generators of 
$F_n$, say $h_j = g_{i_1}^{\eta_{1}} \cdots g_{i_{\ell}}^{\eta_{\ell}}$. \esp If we 
identify $F_n$ to its dynamical realization (and therefore $h_j$ to 
$g_{i_1,0}^{\eta_{1}} \cdots g_{i_{\ell},0}^{\eta_{\ell}}$), then from 
the fact that \esp $h_j (0) > 0$ \esp and that $(g_{i,k})_k$ converges to $g_i$ in the 
compact-open topology, one easily deduces that, if $k$ is large enough, then \esp 
$g_{i_1,k}^{\eta_{1}} \cdots g_{i_{\ell},k}^{\eta_{\ell}}$ \esp sends the origin 
into a positive real number. \esp This means that the element in $\Gamma_k$ 
corresponding to $h_j$ is positive with respect to any ordering obtained from 
the action of $\Gamma_k$ on the line using any dense sequence of points $(x_n)$ 
starting at the origin. Since this is true for each index \esp 
$j \!\in\! \{1,\ldots,m\}$, \esp for $k$ large enough all of the elements in 
$\Gamma_k$ corresponding to the $h_j$'s \esp are simultaneously positive for 
all of such orderings. Let us fix one of these orderings $\preceq'_k$ on $\Gamma_k$, 
as well as an ordering $\preceq_{N_k}$ on $N_k$. Denoting by $[h]$ the class modulo 
$N_k$ of an element \esp $h \in F_n$, \esp let us consider the ordering $\preceq^1_k$ 
(resp. $\preceq^2_k$) on $F_n$ defined by \esp $h \succ id$ \esp if and only if \esp 
$[h] \succ'_k id$, \esp or if \esp $h \!\in\! N_k$ \esp and \esp $h \succ_{N_k} id$ 
\esp (resp. \esp $h \prec_{N_k} id$). \esp The elements $h_j$ are still positive 
with respect to \esp $\preceq^1_k$ \esp and \esp $\preceq^2_k$ \esp for $k$ large 
enough. On the other hand, \esp $\preceq^1_k$ \esp and \esp $\preceq^2_k$ \esp are 
different, because they do not coincide on $N_k$. Therefore, at least one of them is 
distinct from $\preceq$, which concludes the proof in this case.

\vspace{0.25cm}

\noindent{\underbar{Case 2.}} Passing to a subsequence if necessary, for every $k$ the 
elements \esp $g_{1,k},\ldots,g_{n,k}$ \esp do not satisfy any non-trivial relation.
 
\vspace{0.15cm}

We first claim that it is possible to change the $g_{i,k}$'s into homeomorphisms of the real 
line so that the dynamical realization of $F_n$ is not topologically semiconjugate to the 
action of $\Gamma_k$ but the latter group still satisfies the properties above (namely, it has 
no global fixed point, and for each \esp $i \!\in\! \{1,\ldots,n\}$ \esp the maps $g_{i,k}$ 
converge to $g_{i,0}$ in the compact-open topology). To show this let us first note that, 
since the \esp $g_{i,k}$'s \esp are topologically conjugate to maps which extend to real 
analytic diffeomorphism of the closed interval $[0,1]$, they have only finitely many fixed 
points. Since topological semiconjugacies send fixed points into fixed points for 
corresponding elements, if one of the generators \esp $g_{1,0},\ldots,g_{n,0}$ \esp of the 
dynamical realization of $\preceq$ has fixed points outside every compact interval of the line, 
then this realization cannot be topologically semiconjugate to the action of $\Gamma_k$. If the 
sets of fixed points of the \esp $g_{i,0}$'s \esp are contained in some compact interval, then 
for each $k$ let us consider an increasing sequence of points \esp $y_l \geq 2^{l}$ \esp which 
are not fixed by the generators \esp $g_{1,k},\ldots,g_{n,k}$. \esp Let us change $g_{1,k}$ 
into a homeomorphisms of the real line which coincides with the original one on the interval 
$[-2^k,2^k]$ and whose set of fixed points outside $[-2^k,2^k]$ coincides with the set \esp 
$\{y_l\! : l \geq k\}$. The new maps \esp $g_{1,k}$ \esp still converge to $g_{1,0}$ 
in the compact-open topology. Moreover, by the choice of the sequence $(y_l)$, there is 
no global fixed point for the group generated by (the new homeomorphism) $g_{1,k}$ and 
$g_{2,k},\ldots,g_{n,k}$. Finally, by looking at the sets of fixed points of $g_{1,k}$ and 
$g_{1,0}$, one easily concludes the nonexistence of a topological semiconjugacy between 
the action of the (new group) $\Gamma_k$ and the dynamical realization of $\preceq$.

\vspace{0.05cm}

Now for each $k$ the new homeomorphisms \esp $g_{1,k},\ldots,g_{n,k}$ \esp may satisfy some 
non-trivial relation. If this is the case for infinitely many \esp $k \!\in\! \mathbb{N}$, 
\esp then one proceeds as in Case 1. If not, then (passing to subsequences if necessary) 
we just need to consider the following two subcases.

\vspace{0.25cm}

\noindent{\underbar{Subcase i.}} The orbit of the origin by each $\Gamma_k$ is free.
 
\vspace{0.15cm}

For each $k$ we may consider the order relation $\preceq_k$ on \esp $F_n \!\sim\! \Gamma_k$ 
\esp obtained from the corresponding action on the line using the orbit of the origin. A 
simple continuity argument as before shows that, for $k$ large enough, the elements $h_j$ are 
$\precede_k$-positive. On the other hand, since the action of $\Gamma_k$ is not topologically 
semiconjugate to the dynamical realization of $\preceq$, Proposition \ref{semi-igual} implies 
that $\preceq_k$ and $\preceq$ do not coincide, thus finishing the proof for this case.

\vspace{0.25cm}

\noindent{\underbar{Subcase ii.}} The orbit of the origin by each $\Gamma_k$ is non free.
 
\vspace{0.15cm}

For a fixed $k$ let us consider a positive element 
\esp $h \!=\! g_{i_1}^{\eta_{1}} \cdots g_{i_{\ell}}^{\eta_{\ell}} \!\in\! F_n$ 
\esp of minimal length $\ell = \ell_k$ for which the map \esp $g_{i_1,k}^{\eta_{1}} 
\cdots g_{i_{\ell},k}^{\eta_{\ell}}$ \esp fixes the origin (here the exponents 
$\eta_{i}$ belong to $\{-1,1\}$). By the choice of $h$, the points \esp\esp $0, 
g_{i_{\ell},k}^{\eta_{\ell}}(0),g_{i_{\ell-1},k}^{\eta_{\ell-1}} g_{i_{\ell},k}^{\eta_{\ell}}(0), 
\ldots, g_{i_2,k}^{\eta_{i_2}} \cdots g_{i_{\ell},k}^{\eta_{\ell}}(0)$ \esp \esp 
are two-by-two distinct. By perturbing slightly the generator $g_{i_1}$ near 
the latter point, we obtain a new group $\Gamma_k'$ such that the new map \esp 
$g_{i_1,k}^{\eta_{1}} \cdots g_{i_{\ell},k}^{\eta_{\ell}}$ \esp corresponding 
to $h$ sends the origin into a negative real number, but all of the elements in 
$\Gamma_k'$ corresponding to the \esp $h_j$'s \esp still send the origin into positive 
real numbers. If the generators of $\Gamma_k'$ satisfy no non-trivial relation, then using 
any dense sequence of points on the line starting with the origin we may induce a new 
ordering $\preceq'$ on \esp $F_n\sim\Gamma_k'$ \esp which still satisfies $h_j\succ' id$, 
but which is different from $\preceq$ since $h \succ id$ and $h \prec' id$. If there is 
some non-trivial relation between the generators of $\Gamma_k'$, then one may proceed 
as in Case 1. This finishes the proof of Theorem A.

%\begin{rem} Theorem A may be easily extended to the 
%case of infinitely generated countable free groups.
%\end{rem}

\begin{ex} In contrast to Theorem A, we will see in Examples \ref{dehor6} and \ref{dehor7} 
that braid groups admit orderings which are isolated in the corresponding space of 
orderings (although these spaces contain homeomorphic copies of the Cantor set~!). 
\label{dehor2}
\end{ex}

%%%%%%%%%%%%%%%%%%%%%%%%%%%%%%%%%%%%%%%%%%%%%%%%%%%%%%%%%%%%%%%%%%%%%%%%%%%%%%%%%%%%%%%%%%%%%%%%%%%%%%%%%

\section{A dynamical approach to some properties of left-invariant orders}

\subsection{Archimedean orders and H\"older's theorem}

\hspace{0.35cm} The main results of this Section are essentially due to H\"older. 
Roughly, they state that free actions on the line can exist only for groups 
admitting an order relation satisfying an {\em Archimedean} type property. Moreover, 
these groups are necessarily isomorphic to subgroups of $(\mathbb{R},+)$, and 
the corresponding actions are semiconjugate to actions by translations.

\vspace{0.025cm}

\begin{defn} A left-invariant total order relation $\precede$ on a group $\Gamma$ is said 
to be \textit{Archimedean} if for all $g, h$ in $\Gamma$ such that $g \!\neq\! id$ 
there exists $n \!\in\! \mathbb{Z}$ such that $g^n \!\succ\! h.$\\
\end{defn}

\vspace{0.025cm}

\begin{prop} \textit{If \esp\esp $\Gamma$ is a group acting freely by homeomorphisms 
of the real line, then $\Gamma$ admits a total bi-invariant order which is Archimedean.} 
\label{existeordenarq}
\end{prop}

\noindent{\bf Proof.} Let us consider the left-invariant order relation $\precede$ in 
$\Gamma$ such that $g \prec h$ if $g(x) < h(x)$ for some (equivalently, for all) 
$x \in \mathbb{R}$. This order relation is total, and since the action is free, 
one easily checks that it is also right-invariant and Archimedean. $\hfill\square$

\vspace{0.45cm}

The converse to the proposition above is a direct consequence to the following one. As we 
will see in the next Section, the hypothesis of bi-invariance for the order is superfluous: 
it suffices for the order to be left-invariant ({\em c.f.} Proposition \ref{conrad-hist}).

\vspace{0.1cm}

\begin{prop} \textit{Every group admitting a bi-invariant Archimedean 
order is isomorphic to a subgroup of \esp \esp $(\mathbb{R},+)$.} 
\label{completa}
\end{prop}

\noindent{\bf Proof.} Assume that a non-trivial group $\Gamma$ admits a bi-invariant Archimedean 
order $\precede$, and let us fix a positive element $f \in \Gamma$. For each $g \in \Gamma$ 
and each $p \in \mathbb{N}$ let us consider the unique integer $q=q(p)$ such that 
$\hspace{0.15cm} f^q \precede g^p \prec f^{q+1}$.

\vspace{0.25cm}

\noindent{\underbar{Claim 1.}} The sequence \esp $q(p) / p$ \esp 
converges to a real number as \esp $p$ \esp goes to infinite.

\vspace{0.15cm}

Indeed, if $\hspace{0.15cm} f^{q(p_1)} \precede g^{p_1} \prec
f^{q(p_1)+1} \hspace{0.15cm}$ and $\hspace{0.15cm} f^{q(p_2)}
\precede g^{p_2} \prec f^{q(p_2)+1} \hspace{0.15cm}$ then
$$f^{q(p_1)+q(p_2)} \precede g^{p_1+p_2} \prec f^{q(p_1)+q(p_2)+2},$$
and therefore
$\hspace{0.15cm} q(p_1) + q(p_2) \leq q(p_1+p_2) \leq
q(p_1)+q(p_2)+1.$ \esp The convergence of the sequence $(q(p) / p)$ 
to some point in $[-\infty,\infty[$ then follows from a classical 
lemma on subaditive sequences \cite[Page 277]{mane}. 
On the other hand, if we denote by 
$\phi(g)$ the limit of \esp $q(p)/p$, \esp then for the integer 
$n \in \mathbb{Z}$ satisfying \esp $f^n \precede g \prec f^{n+1}$ 
\esp one has \esp $ f^{np} \precede g^p \prec f^{(n+1)p}$, \esp and therefore
$$n = \lim\limits_{p \rightarrow \infty} \frac{np}{p} \leq \phi(g)
\leq \lim\limits_{p \rightarrow \infty} \frac{(n+1)p - 1}{p} = n+1.$$

\vspace{0.25cm}

\noindent{\underbar{Claim 2.}} The map $\phi: \Gamma
\rightarrow (\mathbb{R},+)$ is a group homomorphism.\\

\vspace{0.15cm}

Indeed, let $g_1,g_2$ be arbitrary elements in $\Gamma$. Let us suppose that 
$g_1g_2 \precede g_2g_1$ (the case where $g_2g_1 \precede g_1g_2$ is analogous). 
Since $\preceq$ is bi-invariant, if \esp $f^{q_1} \precede g_1^p \prec f^{q_1 +1}$ 
\esp and \esp $f^{q_2} \precede g_2^p \prec f^{q_2 +1}$ \esp then
$$\hspace{0.15cm} f^{q_1+q_2} \precede g_1^p g_2^p \precede
(g_1g_2)^p \precede g_2^p g_1^p \prec f^{q_1 + q_2 + 2} \hspace{0.15cm}.$$
From this one concludes that 
$$\phi(g_1) + \phi(g_2) = \lim\limits_{p \rightarrow \infty}
\frac{q_1+q_2}{p} \leq \phi(g_1g_2) \leq \lim\limits_{p
\rightarrow \infty} \frac{q_1+q_2 +1}{p} = \phi(g_1)+\phi(g_2),$$ 
and therefore $\phi(g_1g_2) = \phi(g_1) + \phi(g_2).$

\vspace{0.25cm}

\noindent{\underbar{Claim 3.}} The homomorphism $\phi$ is one to one.\\

\vspace{0.15cm}

Note that $\phi$ is order preserving, in the sense that if $g_1 \precede g_2$ 
then $\phi(g_1) \leq \phi(g_2)$. Moreover, $\phi(f)=1.$ Let $h$ be an element 
in $\Gamma$ such that $\phi(h)=0$. Assume that $h \neq id$. Then there exists 
$n \in \mathbb{Z}$ such that $h^n \succeq f$. From this one concludes that 
$0 = n \phi(h) = \phi(h^n) \geq \phi(f) = 1,$ which is absurd. Therefore, 
if $\phi(h) = 0$ then $h = id$, and this concludes the proof. $\hfill\square$

\vspace{0.5cm}

If $\Gamma$ is an infinite group acting freely on the line, then we can fix the order 
relation introduced in the proof of Proposition \ref{existeordenarq}. This order allows us 
to construct an embedding $\phi$ from $\Gamma$ into $(\mathbb{R},+)$. If $\phi(\Gamma)$ is 
isomorphic to $(\mathbb{Z},+)$ then the action of $\Gamma$ is conjugate to the action 
by integer translations. In the other case, the group $\phi(\Gamma)$ is dense in 
$(\mathrm{R},+)$. For each point $x$ in the line we define
$$\varphi(x) = \sup \{ \phi(h) \in \mathbb{R}: h(0) \leq x \}.$$
It is easy to see that $\varphi \! : \mathbb{R} \rightarrow \mathbb{R}$ 
is a non-decreasing map. Moreover, it satisfies the equality 
\esp $\varphi(h(x)) = \varphi(x) + \phi(h)$ \esp for all 
$x \in \mathbb{R}$ and all $h \in \Gamma$. Finally, $\varphi$ is continuous, 
as otherwise $\mathbb{R} \setminus \varphi(\mathbb{R})$ would be a 
non-empty open set invariant by the translations of $\phi(\Gamma)$, 
which is impossible.

\vspace{0.15cm}

To summarize, if $\Gamma$ is a group acting freely on the line, 
then its action semiconjugates to an action by translations.

%%%%%%%%%%%%%%%%%%%%%%%%%%%%%%%%%%%%%%%%%%%%%%%%%%%%%%%%%%%%%%%%%%%%%%%%%%%%%%%%%%%%%%%%%

\subsection{Almost free actions and bi-invariant orders}
\label{invariant}

\hspace{0.35cm} We will say that the action of a group $\Gamma$ of orientation-preserving 
homeomorphisms of the line is {\em almost free} if for every element $g \in \Gamma$ 
one has either $g(x) \geq x$ for all $x \in \mathbb{R}$ or $g(x) \leq x$ for all 
$x \in \mathbb{R}$. The following proposition gives the algebraic counterpart 
of this notion.

\vspace{0.15cm}

\begin{prop} {\em A countable group $\Gamma$ admits a faithful almost 
free action on the real line if and only if it is bi-orderable.}
\label{casi-obvio}
\end{prop}

\noindent{\bf Proof.} If $\Gamma$ is bi-orderable, then the action on the line of 
the dynamical realization associated to any of its numberings is almost free. 
Indeed, if $g \succ id$ then $gg_i \succ g_i$ for all $g_i \in \Gamma$, and 
therefore $g(t(g_i)) = t(gg_i) > t(g_i)$. By the construction of the dynamical 
realization, this implies that \esp $g(x) \geq x$ \esp for all $x \in \mathbb{R}$. 
In an analogous way, for $g \prec id$ one has \esp $g(x) \leq x$ \esp for all 
$x \in \mathbb{R}$, thus showing that the action is almost free.
 
Conversely, let $\Gamma$ be a group of homeomorphisms of the line 
whose action is almost free. We claim that the order $\precede$ associated 
to any dense sequence $(x_n)$ of points in $\mathbb{R}$ is bi-invariant. Indeed, 
if $f \succeq id$, then the graph of $f$ does not have any point below the diagonal. 
Obviously, if $g$ is any element in $\Gamma$, then the same is true for the graph 
of $gfg^{-1}$. This clearly implies that $gfg^{-1} \succeq id$, thus proving 
the bi-invariance of $\precede$. $\hfill\square$

\vspace{0.15cm}

\begin{ex} Groups of piecewise-linear homeomorphisms of the interval are 
bi-orderable: it suffices to define $\precede$ by \esp $f \!\succ\! id$ \esp when 
$f (x_f + \varepsilon) > x_f + \varepsilon$ for every $\varepsilon \!>\! 0$ sufficiently 
small, where \esp $x_f = \inf \{x: \esp f(x) \neq x \}$. \esp As an application 
of the previous proposition, we obtain for example a non standard action of 
Thompson's group F on the line. (Compare \cite{nos}.) A similar construction applies to 
countable groups of germs at the origin of one dimensional real-analytic diffeomorphisms. 
\end{ex}

\vspace{0.15cm}

To close this Section, we give a dynamical proof of a fact first 
remarked by Conrad in \cite{conrad}.

\vspace{0.2cm}

\begin{prop} {\em Every Archimedean left-invariant total order on a group is bi-invariant.}
\label{conrad-hist}
\end{prop}

\noindent{\bf Proof.} Let $\{f_1,\ldots,f_k\}$ be any finite family of elements in a group $\Gamma$ 
endowed with a total order relation $\precede$ which is left-invariant and Archimedean. Let us consider 
some numbering $(h_n)_{n \geq 0}$ of the group generated by them, as well as the corresponding 
dynamical realization. We claim that this action is free. Indeed, if not then there exist 
$h \! \in \! \langle f_1,\ldots,f_k \rangle$ and an interval $]a,b[$ which is not the whole 
line such that $h$ fixes $a$ and $b$ and has no fixed point in $]a,b[$. By the comments 
after Proposition \ref{thorden}, a moment reflexion shows that such an interval $]a,b[$ can 
be taken so that $b \neq +\infty$. Moreover, there exists some point of the form $t(h_i)$ 
inside $]a,b[$, and by conjugating by $h_i$ if necessary, we may assume that $t(id)$ 
belongs to $]a,b[$. Now since dynamical realizations of non-trivial orderable 
groups have no global fixed point, there must exist some 
$\bar{h} \in \langle f,g \rangle$ such that $\bar{h} (t(id)) > b$. We thus have 
$h^n(t(id)) < b < \bar{h}(t(id))$ for all $n \in \mathbb{Z}$, which implies that 
$h^n \prec \bar{h}$ for all $n \in \mathbb{Z}$. Nevertheless, 
this violates the Archimedean property for $\precede$.

Now let $f \!\prec\! g$ and $h$ be three elements in $\Gamma$. Since the dynamical 
realization associated to the group generated by them is free and $f(t(id)) \!<\! g(t(id))$, 
one has $f(t(h)) \!<\! g (t(h))$, that is, $t(fh) \!<\! t(gh)$. By construction, this implies 
that $fh \! \prec \! gh$. Since $f \!\prec\! g$ and $h$ were arbitrary elements of $\Gamma$, 
this shows that $\precede$ is right-invariant. $\hfill\square$

%%%%%%%%%%%%%%%%%%%%%%%%%%%%%%%%%%%%%%%%%%%%%%%%%%%%%%%%%%%%%%%%%%%%%%%%%%%%%%%%%%%%%%%%%

\subsection{The Conrad property and crossed elements (resilient orbits)}
\label{todo-conrad}

\subsubsection{The Conrad property} 
\label{general}

\hspace{0.35cm} A left-invariant total order relation $\precede$ on a group $\Gamma$ 
satisfies the {\em Conrad property} (or it is a Conradian order, or simply a 
$\mathcal{C}$-order) if for all positive elements $f,g$ there exists $n \!\in\! \mathbb{N}$ 
such that $f g^n \!\succ\! g$. If a group admits such an order, then it is said to be 
Conrad orderable. These notions were introduced in \cite{conrad}, where several 
characterizations are given (see also \cite{botto,glass,koko}). Nevertheless, 
the following quite simple (and unexpectedly useful) proposition does not 
seem to appear in the literature. 

\vspace{0.15cm}

\begin{prop} {\em If $\precede$ is a $\mathcal{C}$-order on a group $\Gamma$, 
then for every positive elements $f,g$ one has $f g^2 \succ g$.}
\label{lesla}
\end{prop}

\noindent{\bf Proof.} Suppose that two positive elements $f,g$ for an ordering $\preceq'$ 
on $\Gamma$ are such that $fg^2 \precede' g$. Then $(g^{-1}fg) g \precede' id$, and since 
$g$ is a positive element this implies that $g^{-1}fg$ is negative, and therefore 
$fg \prec' g$. Now for the positive element $h = fg$ and every $n \in \mathbb{N}$ one has
\begin{multline*}
f h^{n} = f (fg)^n = f (fg)^{n-2} (fg) (fg) \\ \prec' 
f (fg)^{n-2} (fg) g = f (fg)^{n-2} fg^2 \precede' f (fg)^{n-2} g 
= f (fg)^{n-3} fg^2 \precede' f (fg)^{n-3}g \precede' \ldots\\ 
\precede' f (fg) g = f fg^2 \precede' fg = h.
\end{multline*}
This shows that $\preceq'$ does not satisfy the Conrad property. $\hfill\square$

\vspace{0.35cm}

The nice argument of the proof above is due to Jim\'enez \cite{leslie}. 
Latter in \S \ref{equivalence} we will see that, in fact, \esp 
$f g^{n+1} \! \succ \! g^n$ \esp for all $n \!\in\! \mathbb{N}$. More generally, 
we will show that if \esp $W(f,g) = f^{m_1} g^{n_1} \cdots f^{m_k} g^{n_k}$ \esp 
is a word such that $\sum m_i \!>\! 0$ and $\sum n_i \!>\! 0$, then $W(f,g)$ 
is a positive element in $\Gamma$  provided that $f$ and $g$ are both positive. 
(Notice that $f g^{n+1} \succ g^n$ is equivalent to $g^{-n} f g^{n+1} \succ id$.) 
However, we were not able to extend the preceding proof for this, and we will need 
the dynamical characterization of the Conrad property (or at least its algebraic 
counterpart, which corresponds to the characterization in terms of  
convex subgroups: see Remark \ref{convex-car}). 

As a first application of Proposition \ref{lesla} we will show that, for every orderable 
group, the subset of $\mathcal{O}(\Gamma)$ formed by the Conradian orders is closed. 
Note that a similar argument to the one given below applies to the (simpler) case 
of bi-invariant orders. (Compare \cite[Proposition 2.1]{sikora}.)

\vspace{0.2cm}

\begin{prop} \label{conrad-cerrado}
{\em If \esp \esp $\Gamma$ is an orderable group, then the set 
of $\mathcal{C}$-orders on $\Gamma$ is closed in $\mathcal{O}(\Gamma)$.}
\end{prop}

\noindent{\bf Proof.} According to Proposition $\ref{lesla}$, an element $\precede$ 
of $\mathcal{O}(\Gamma)$ is not Conradian if and only if there exists two elements 
$f \succ id$ and $g \succ id$ such that $fg^2 \precede g$, which necessarily implies 
that $g^{-1}fg^2 \prec id$. Since the sets $U_{id,f}$, $U_{id,g}$, and 
$U_{id,g^{-2}f^{-1}g}$, are clopen, the set 
$$U(f,g) = U_{id,f} \cap U_{id,g} \cap U_{id,g^{-2}f^{-1}g} = \{\precede : 
\esp \esp \esp f \succ id, \esp \esp g \succ id, \esp \esp g^{-1} f g^2 \prec id\}$$
is open for every $f,g$ in $\Gamma$ different from the identity. Thus, 
the union of the \esp $U(f,g)$'s \esp is open, and therefore its complementary 
set (that is, the set of $\mathcal{C}$-orders) is closed. $\hfill\square$

\vspace{0.1cm}

\begin{question} What can be said about the topology of the set of Conradian orders~? When is 
the set of Conradian orders open or at least of non-empty interior in $\mathcal{O}(\Gamma)$~?
\footnote{Added in Proof: This has been partially answered in \cite{rivas}.}
\end{question}

\vspace{0.01cm}

As another application of Proposition \ref{lesla}, we give a criterion 
for Conrad orderability which is similar to those of Proposition \ref{bur}.

\vspace{0.1cm}

\begin{prop} {\em A group $\Gamma$ admits a Conradian order if and only if the following condition is 
satisfied:  \esp for every finite family of elements $g_1,\ldots,g_k$ which are different from the 
identity, there exists a family of exponents $\eta_i \!\in \! \{-1,1\}$ such that \esp $id$ \esp does not 
belong to the smallest semigroup \esp $\langle \langle g_1^{\eta_1},\ldots,g_k^{\eta_k} \rangle \rangle$ 
which simultaneously satisfies the following two properties:\\ 

\noindent -- it contains all the elements \esp $g_i^{\eta_i}$; 

\noindent -- for all $f,g$ in the semigroup, the element \esp $f^{-1}gf^{2}$ \esp also belongs to it.}
\label{notario}
\end{prop}

\noindent{\bf Proof.} The necessity of the condition follows as a direct application of Proposition 
\ref{lesla} after choosing $\eta_i$ in such a way that $g_i^{\eta_i}$ is a positive element of 
$\Gamma$. To prove that the condition is sufficient, one proceeds as in the case of Proposition 
\ref{bur} by introducing the sets $C \mathcal{X}(g_1\ldots,g_k;\eta_1,\ldots,\eta_k)$ formed 
by all the functions \esp $\sg$ \esp for which \esp $\sg(g) \!=\! +$ \esp and \esp 
$\sg(g^{-1}) \!=\! -$ \esp for each $g$ contained in the semigroup \esp 
$\langle \langle g_1^{\eta_1},\ldots,g_k^{\eta_k} \rangle \rangle$. \esp 
We leave the details to the reader. $\hfill\square$

\vspace{0.35cm}

It easily follows from the criterion above that residually Conrad orderable 
groups are Conrad orderable.\footnote{Recall that, if P is some group property, 
then a group $\Gamma$ is said to be {\em residually} P if for every  
$g \in \Gamma \setminus \{id\}$ there exists a surjective group 
homomorphism from $\Gamma$ to a group $\Gamma_{g}$ such that the image 
of $g$ is non-trivial.} As a more interesting application, we give a short 
proof of a theorem due to Brodskii \cite{brodski}, and independently 
obtained by Rhemtulla and Rolfsen \cite{remt}. For the statement, recall 
that a group is said to be {\em locally indicable} if for each non-trivial 
finitely generated subgroup there exists a non-trivial homomorphism into $(\mathbb{R},+)$.

\vspace{0.12cm}

\begin{prop} {\em Every locally indicable group is Conrad orderable.}
\label{notar}
\end{prop}

\noindent{\bf Proof.} We need to check that every locally indicable group $\Gamma$ satisfies 
the condition of Proposition \ref{notario}. Let $\{g_1,\ldots,g_k\}$ be any finite family 
of elements in $\Gamma$ which are different from the identity. By hypothesis, there exists a 
non-trivial homomorphism $\phi_1\!: \langle g_1,\ldots,g_k \rangle \rightarrow (\mathbb{R},+)$. 
Let \esp $i_1,\ldots,i_{k'}$ \esp be the indexes (if any) such that $\phi_1 (g_{i_{j}}) = 0$. 
Again by hypothesis, there exists a non-trivial homomorphism 
$\phi_2\!: \langle g_{i_1},\ldots,g_{i_{k'}} \rangle \rightarrow (\mathbb{R},+)$. Letting 
\esp $i_1',\ldots,i_{k''}'$ \esp be the indexes in $\{i_1,\ldots,i_{k'}\}$ for which 
$\phi_2 (g_{i_j'}) = 0$, we may choose a non-trivial homomorphism 
$\phi_3\!: \langle g_{i_1'},\ldots,g_{i_{k''}'} \rangle \rightarrow (\mathbb{R},+)$... Note that 
this process must finish in a finite number of steps (indeed, it stops in at most $k$ steps). 
Now for each $i \!\in\! \{1,\ldots,k\}$ choose the (unique) index $j(i)$ such that $\phi_{j(i)}$ 
is defined at $g_i$ and $\phi_{j(i)} (g_i) \neq 0$, and let $\eta_i \!\in\! \{-1,1\}$ be so that 
$\phi_{j(i)} (g_i^{\eta_i}) > 0$. We claim that this choice of exponents $\eta_i$ is ``compatible". 
Indeed, for every index $j$ and every $f,g$ for which $\phi_j$ are defined, one has  
\esp \hspace{0.07cm} $\phi_j(f^{-1} g f^2) = \phi_j (f) + \phi_j (g).$ \esp 
\hspace{0.07cm} Therefore, \esp $\phi_1(h) \geq 0$ \esp for every \esp 
$h \!\in\! \langle \langle g_1^{\eta_1},\ldots,g_k^{\eta_k} \rangle \rangle$. 
\esp Moreover, if \esp $\phi_1(h) = 0$, \esp then $h$ actually belongs to \esp 
$\langle\langle g_{i_1}^{\eta_{i_1}},\ldots,g_{i_{k'}}^{\eta_{i_{k'}}} \rangle\rangle$. 
\esp In this case, the preceding argument shows that \esp $\phi_2(h) \!\geq\! 0$, \esp 
with equality if and only if 
$h \!\in\! \langle\langle g_{i_1´}^{\eta_{i_1'}},\ldots,g_{i_{k''}'}^{\eta_{i_{k''}'}} \rangle\rangle$... 
Continuing in this way, one concludes that $\phi_j (h)$ must be strictly 
positive for some index $j$. Thus, the element $h$ cannot be equal to 
the identity, and this concludes the proof. $\hfill\square$

\vspace{0.4cm}

As we will see in \S \ref{equivalence}, 
the converse of Proposition \ref{notar} also 
holds ({\em c.f.} Proposition \ref{indicar-loc}).

%%%%%%%%%%%%%%%%%%%%%%%%%%%%%%%%%%%%%%%%%%%%%%%%%%%%%%%%%%%%%%%%%%%%%%%%%%%%%%%%%%%%%%%%%%%%%%

\subsubsection{Crossed elements, invariant Radon measures, and translation numbers}
\label{bp}

\hspace{0.35cm} We say that two orientation-preserving homeomorphisms of the real line 
are {\em crossed} \esp on an interval $]a,b[$ if one of them fixes $a$ and $b$ and no other 
point in $[a,b]$, while the other one sends $a$ or $b$ into $]a,b[$. Here we allow the 
case where \esp $a = - \infty$ \esp or \esp $b = +\infty$.
 
If $f$ and $g$ are homeomorphisms of the line which are contained in a group 
without crossed elements, and if $f$ has a fixed point $x_0$ which is not fixed 
by $g$, then the fixed points of $g$ immediately to the left and to the right of 
$x_0$ are also fixed by $f$. This gives a quite particular combinatorial structure 
for the dynamics of groups of homeomorphisms of the line without crossed elements. 
To understand this dynamics better, one can use an extremely useful tool for 
detecting fixed points of elements, namely the translation number associated to 
an invariant Radon measure. The Proposition below is originally due to Beklaryan 
\cite{bekl}. Here we provide a proof taken from \cite[Section 2.1]{growth}.

\vspace{0.2cm}

\begin{prop} {\em Let $\Gamma$ be a finitely generated group of orientation-preserving 
homeomorphisms of the real line. If \esp $\Gamma$ has no crossed elements, then 
$\Gamma$ preserves a (non-trivial) Radon measure on $\mathbb{R}$ (that is, a measure 
on the Borelean sets which is finite on the compact subsets of $\mathbb{R}$).}
\label{radon}
\end{prop}

\noindent{\bf Proof.} If $\Gamma$ has global fixed points in $\mathbb{R}$, then 
the Dirac delta measure on any of such points is invariant by the action. 
Assume in what follows that the $\Gamma$-action on $\mathbb{R}$ has no 
global fixed point, and take a finite system $\{f_1,\ldots,f_k\}$ 
of generators for $\Gamma$. We first claim that (at least) one of 
these generators does not have interior fixed points. Indeed, suppose by 
contradiction that all the maps $f_i$ have interior fixed points, and let 
$x_1 \!\in \mathbb{R}$ be any fixed point of $f_1$. If $f_2$ 
fixes $x_1$, then letting $x_2=x_1$ we have that $x_2$ is fixed by both $f_1$ 
and $f_2$. If not, choose a fixed point $x_2 \!\in \mathbb{R}$ for $f_2$ 
such that $f_2$ does not fix any point between $x_1$ and $x_2$. 
Since $f_1$ and $f_2$ are non crossed on 
any interval, $x_2$ must be fixed by $f_1$. Now if $x_2$ is fixed by $f_3$, 
let $x_3=x_2$; if not, take a fixed point $x_3 \!\in \mathbb{R}$ for 
$f_3$ such that $f_3$ has no fixed point between $x_2$ and $x_3$. 
The same argument as before shows that $x_3$ is fixed by $f_1,f_2$, and $f_3$. 
Continuing in this way, we find a common fixed point for all of the generators $f_i$, 
and so a global fixed point for the action of $\Gamma$, thus giving a contradiction.

Now we claim that there exists a non-empty minimal invariant 
closed set for the action of $\Gamma$ on $\mathbb{R}$. To 
prove this, consider a generator $f = f_i$ without fixed points, 
fix any point $x_0 \!\in \mathbb{R}$, and let $I$ be the interval 
$[x_0,f(x_0)]$ if $f(x_0) > x_0$, and $[f(x_0),x_0]$ if $f(x_0) < x_0$. On 
the family $\mathcal{F}$ of non-empty closed invariant subsets of 
$\mathbb{R}$, let us consider the order relation $\precede$ given by \esp 
$K_1 \succeq K_2$ if \esp $K_1 \cap I \subset K_2 \cap I$. \esp 
Since $f$ has no fixed point, every orbit by $\Gamma$ must 
intersect the interval $I$, and so 
$K \cap I$ is a non-empty compact set for all 
$K \in \mathcal{F}$. Therefore, we can apply Zorn Lemma to obtain a maximal 
element for the order $\precede$, and this element is the intersection with 
$I$ of a minimal $\Gamma$-invariant non-empty closed subset of $\mathbb{R}$.

Consider now the non-empty minimal invariant closed set $K$ obtained above. 
Note that its boundary $\partial K$ as well as the set of its accumulation 
points $K'$ are also closed sets invariant by $\Gamma$. Because of 
the minimality of $K$, there are three possibilities:

\vspace{0.25cm}

\noindent \underbar{Case 1.} $K' = \emptyset$.
 
\vspace{0.1cm}

In this case, $K$ is discrete, that is, $K$ coincides 
with the set of points of a sequence $(y_n)_{n \in \mathbb{Z}}$ satisfying 
$y_{n} \!<\! y_{n+1}$ for all $n$ and without accumulation points inside 
$\mathbb{R}$. It is then easy to see that the Radon measure \esp 
$\mu = \sum_{n \in \mathbb{Z}} \delta_{y_n}$ \esp is invariant by $\Gamma$.

\vspace{0.25cm}

\noindent \underbar{Case 2.} $\partial K = \emptyset$.

\vspace{0.1cm}
 
In this case, $K$ coincides with the whole line. We claim that the action of $\Gamma$ 
is free. Indeed, if not let $]u,v[$ be an interval strictly contained in $\mathbb{R}$ 
and for which there exists an element $g \in \Gamma$ fixing $]u,v[$ and with no fixed 
point inside it. Since the action is minimal, there must be some $h \in \Gamma$ sending 
a real endpoint of $]u,v[$ inside $]u,v[$; however, this implies that $g$ and $h$ are 
crossed on $[u,v]$, contradicting our assumption. Now the action of $\Gamma$ being 
free, H\"older's theorem implies that $\Gamma$ is topologically conjugate to a (in 
this case dense) group of translations. Pulling back the Lebesgue measure by 
this conjugacy, we obtain an invariant Radon measure for the action of $\Gamma$.

\vspace{0.25cm}

\noindent \underbar{Case 3.} $\partial K = K' = K$.
 
\vspace{0.1cm}

In this case, $K$ is ``locally'' a Cantor set. Collapsing to a point the closure of 
each connected component of the complementary set of $K$, we obtain a topological 
line on which the original action induces (by semi-conjugacy) an action of $\Gamma$. 
As in the second case, one easily checks that the induced action is free, hence 
it preserves a Radon measure. Pulling back this measure by the semi-conjugacy, 
one obtains a Radon measure on $\mathbb{R}$ which is invariant by the original 
action. $\hfill\square$

\vspace{0.4cm}

Recall that for (non necessarily finitely generated) groups of orientation-preserving 
homeomorphisms of the line preserving a (non-trivial) Radon measure $\mu$, there is an 
associated {\em translation number} function $\tau_{\mu}: \Gamma \rightarrow \mathbb{R}$ 
defined by 
$$\tau_{\mu}(g) = \left \{ \begin{array} {l} 
\mu([x_0,g(x_0)[) \hspace{0.67cm} \mbox{ if } \esp g(x_0) > x_0,\\ 
0 \hspace{2.8cm} \mbox{ if } \esp g(x_0) = x_0,\\
- \mu([g(x_0),x_0[) \hspace{0.4cm} \mbox{ if } \esp g(x_0) < x_0, \end{array} \right.$$
where $x_0$ is any point of the line \cite{plante}. (One easily checks that 
this definition is independent of $x_0$.) The following properties are satisfied 
(the verification is easy and may be left to the reader):\\

\vspace{0.07cm}

\noindent (i) $\tau_{\mu}$ is a group homomorphism;\\

\vspace{0.07cm}

\noindent (ii) $\tau_{\mu}(g) = 0$ if and only if $g$ has fixed points; in this 
case, the support of $\mu$ is contained in the set of these points;\\

\vspace{0.07cm}

\noindent (iii) $\tau_{\mu}$ is trivial if and only if there is no global fixed 
point for the action of $\Gamma$.\\

\vspace{0.07cm}

\begin{rem} For codimension-one foliations, the notion of crossed elements corresponds to 
that of {\em resilient leaves} ({\em feuilles ressort}). In this context, an analogous of 
Proposition \ref{radon} holds, but its proof is more difficult and uses completely 
different ideas (see \cite[Th\'eor\`eme E]{DKN}).
\end{rem}

%%%%%%%%%%%%%%%%%%%%%%%%%%%%%%%%%%%%%%%%%%%%%%%%%%%%%%%%%%%%%%%%%%%%%%%%%%%%%%%%%%%%%%%%%

\subsubsection{The equivalence}
\label{equivalence}

\hspace{0.35cm} Propositions \ref{conrad-yo} and \ref{conrad-yo2} below give the equivalence 
between the Conrad property and the nonexistence of crossed elements for the actions on the line.

\vspace{0.1cm}

\begin{prop} {\em Let $\Gamma$ be a countable group with a $\mathcal{C}$-order $\precede$. 
For any numbering $(g_n)_{n \geq 0}$ of $\Gamma$, the corresponding dynamical realization 
is a subgroup of $\mathrm{Homeo}_+(\mathbb{R})$ without crossed elements.}
\label{conrad-yo}
\end{prop}

\noindent{\bf Proof.} The claim is obvious if $\Gamma$ is trivial; thus, we will 
assume in the sequel that $\Gamma$ contains infinitely many elements. Let us 
suppose that there exist $f,g$ in $\Gamma$ and an interval $[a,b]$ such that 
(for their dynamical realizations one has) $Fix(f) \cap [a,b] = \{a,b\}$ and 
$g(a) \! \in ]a,b[$ (the case where $g(b)$ belongs to $]a,b[$ is analogous). 
Changing $f$ by its inverse if necessary, we can suppose that $f(x) \!<\! x$ 
for all $x \!\! \in ]a,b[$. As we already observed after the proof of Proposition 
\ref{thorden}, there must exist some element $g_i \!\in\! \Gamma$ such that $t(g_i)$ 
belongs to the interval $]a,b[$. Let $j \geq 0$ be the index such that $g_j \!=\! id$. 
By conjugating $f$ and $g$ by the element $g_i^{-1}$ if necessary, we may 
assume that $t(g_j) \!=\! t(id)$ belongs to $]a,b[$. Furthermore, changing 
$g$ by $f^{-n}g$ for $n$ large enough, we may assume that $g(a)\!>\!t(g_j)$. 
Let us define $c \!=\! g(a) \! \in ]t(g_j),b[$, and let us fix a point 
$d\! \in ]c,b[$. Since $gf^n(a) \!=\! c$ for all $n\!\in\!\mathbb{N}$, 
and since $gf^n(d)$ converges to $c \!<\! d$ as $n$ goes to infinity, 
for $n\!\in\!\mathbb{N}$ sufficiently big the map $h_n \!=\! gf^n$ satisfies 
$h_n(a)\!>\!a$, \esp $h_n(d) \!<\! d$, \esp $\mathrm{Fix}(h_n) \cap ]a,d[ \subset \! 
[c_n,c_n'] \! \subset ]c,h_n(d)[$ \esp and $\{c_n,c_n'\} \subset \mathrm{Fix}(h_n)$ 
for some sequences $(c_n)$ and $(c_n')$ converging to $c$ by the right. (See Figure 
1 below.) Note that each $h_n$ satisfying the preceding properties is positive, because 
from \esp $h_n \big( t(g_j) \big) \!>\! h_n (a) \!=\! c \!>\! t(g_j)$ \esp one 
concludes that $t(h_n) > t(id)$, and by the construction of the 
dynamical realization this implies that $h_n \succ id$. 

Let us fix $m > n$ large enough so that the preceding properties are satisfied for 
$h_m$ and  $h_n$, and such that \esp $[c_m,c_m'] \! \subset ]c,c_n[$. \esp Let us fix 
$k \in \mathbb{N}$ sufficiently big so that \esp $h_n^k (a) > h_m(c_n)$, \esp and 
let us define $h = h_n^k$. For each $i \!\in\! \mathbb{N}$ one has \esp 
$h^i \big( t(g_j) \big) \! \in \hspace{0.08cm} ]h_m(c_n),c_n[$, \esp 
and therefore 
$$h_m h^i \big( t(g_j) \big) \!<\! h_m(c_n) 
\!<\! h(a) \!<\! h \big( t(g_j) \big).$$ 
Thus, $h_m h^i \prec h$ for each $i \!\in\! \mathbb{N}$. Nevertheless, this 
in contradiction with the Conrad property for the order $\precede$. $\hfill\square$

\vspace{0.6cm}

%%%%%%%%%%%%%%%%%%%%%%%%%%%%%%%%%%%%%%%%%%%%%%%%%%%%%%%%%%%%%%%%%%%%%%%%%%%%%

\beginpicture

\setcoordinatesystem units <1cm,1cm>

\putrule from 0 0 to 0 8
\putrule from 0 0 to 8 0
\putrule from 8 0 to 8 8 
\putrule from 0 8 to 8 8

\put{$a$} at 0 -0.4
\put{$t(g_j)$} at 1.2 -0.4
\put{$c$} at 2 -0.4
\put{$c_m$} at 2.7 -0.4
\put{$c_m'$} at 3.3 -0.4
\put{$c_n$} at 6 -0.4
\put{$c_n'$} at 6.6 -0.4
\put{$d$} at 7.2 -0.4
\put{$b$} at 8 -0.4
\put{$c$} at -0.4 2 
\put{$f$} at 5 1 
\put{$g$} at 1.75 7 
\put{$h \!=\! h_n^k$} at 2.25 5.4 
\put{$h_n$} at 2.8 4.32   
\put{$h_m$} at 2.6 3.2 
\put{$h_m(c_n)$} at -0.87 3.65
\put{$h(a)$} at -0.55 4.32 

\setquadratic 

\plot 
0 2 1.6 2.55 2.7 2.7 /

\plot 
2.7 2.7 2.92 2.82 3 3 /

\plot 
3 3 3.1 3.2 3.3 3.3 /

\plot 
0 0 5.8 2.2 8 8 /

\plot 
0 2 1 6 2.3 8 /

\plot 
6 6 6.25 6.1 6.3 6.3 /

\plot 
6.3 6.3 6.36 6.5 6.6 6.6 /

\plot 
6 6 6.15 6.05 6.3 6.3 /

\plot 
6.3 6.3 6.42 6.52 6.6 6.6 /

%%%%%%%%%%%%%%%%%%%%%%%%%%%%%%%%%%%%%%%%%%%%%%%%%%%%%%%%%%%%%%%%%%%%%

\setlinear

\plot
0 2 
6 6 /

\plot 
0 4.2  
6 6 /

\plot 
6 6 
6.6 6.6 /

\plot 
3.3 3.3 
6 3.7 /

\plot 
6 3.7 
7.2 3.9 /

\plot 
6.6 6.6 
7.2 6.8 / 

\plot 
6.6 6.6 
7.2 6.7 /

%%%%%%%%%%%%%%%%%%%%%%%%%%%%%%%%%%%%%%%%%%%%%%%%%%%%%%%%%%%%%%%%%%%%%%%%%%%%%%%%%%%%%%%%

\plot 3.3 3.3  3.5 3.3 /

\plot 3.7 3.3  3.9 3.3 /

\plot 4.1 3.3 4.3 3.3 /

\plot 4.5 3.3  4.7 3.3 /

\plot 4.9 3.3  5.1 3.3 /

\plot 5.3 3.3 5.5 3.3 /

\plot 5.7 3.3 5.9 3.3 /

%%%%%%%%%%%%%%%%%%%%%%%%%%%%%%%%%%%%%%%%%%%%%%%%%%%%%%%%%%%%%%%%%%%%%%%%%%%

\plot 
3.3 3.3  3.3 3.5 /

\plot 
3.3 3.7  
3.3 3.9 /

\plot
3.3 4.1 
3.3 4.3 /

\plot 
3.3 4.5
3.3 4.7 /

\plot 
3.3 4.9 
3.3 5.1 /

\plot
3.3 5.3 
3.3 5.5 /

\plot 
3.3 5.7 
3.3 5.9 /

%%%%%%%%%%%%%%%%%%%%%%%%%%%%%%%%%%%%%%%%%%%%%%%%%%%%%%%%%%%%%%%%%%%%%%%%%%%

\plot 
3.3 6  
3.5 6 /

\plot 
3.7 6  
3.9 6 /

\plot
4.1 6 
4.3 6 /

\plot 
4.5 6  
4.7 6 /

\plot 
4.9 6  
5.1 6 /

\plot
5.3 6 
5.5 6 /

\plot 
5.7 6 
5.9 6 /

%%%%%%%%%%%%%%%%%%%%%%%%%%%%%%%%%%%%%%%%%%%%%%%%%%%%%%%%%%%%%%%%%%%%%%%%%%%

\plot 
6 3.3 
6 3.5 /

\plot 
6 3.7   
6 3.9 /

\plot
6 4.1  
6 4.3 /

\plot 
6 4.5   
6 4.7 /

\plot 
6 4.9   
6 5.1 /

\plot
6 5.3  
6 5.5 /

\plot 
6 5.7 
6 5.9 /

%%%%%%%%%%%%%%%%%%%%%%%%%%%%%%%%%%%%%%%%%%%%%%%%%%%%%%%%%%%%%%%%%%%%%%%%%%%%%

\setdots

\plot 
0 0 
8 8 /

\putrule from 1.2 0 to 1.2 3.7 
\putrule from 2 0 to 2 2 
\putrule from 0 2 to 2 2 
\putrule from 2.72 0 to 2.72 2.7 
\putrule from 3.3 0 to 3.3 3.3 
\putrule from 6 0 to 6 3.3 
\putrule from 6.6 0 to 6.6 6.6 
\putrule from 0 6 to 3.3 6 
\putrule from 0 3.7 to 6 3.7 
\putrule from 7.2 0 to 7.2 7.2 

\put{Figure 1} at 4 -1 

\put{} at -4.24 0 

\endpicture

%%%%%%%%%%%%%%%%%%%%%%%%%%%%%%%%%%%%%%%%%%%%%%%%%%%%%%%%%%%%%%%%%%%%%%%%%%%%%%%%%%%%%%%%%%

\vspace{0.76cm}

The reader should note that, for the positive elements $h$ and $\bar{h} = h_m$ 
that we found, one has $W_1 (h,\bar{h}) \prec W_2 (h,\bar{h})$ 
for all reduced words $W_1,W_2$ in positive powers such that $W_1$ (resp. 
$W_2$) begins with a power of $\bar{h}$ (resp. $h$). Therefore, the 
following general characterization for the Conrad property holds: 
a left-invariant total order relation $\precede$ on a group $\Gamma$ is a 
$\mathcal{C}$-order if and only if for every pair of positive elements $f,g$ 
in $\Gamma$ one has $W_1 (f,g) \succeq W_2 (f,g)$ for some reduced words 
$W_1,W_2$ in positive powers such that $W_1$ (resp. $W_2$) begins with a 
power of $f$ (resp. $g$). This shows in particular that all orderings on an 
orderable group without free semigroups on two generators are $\mathcal{C}$-orders. 
(This fact was first proved by Longobardi, Maj, and Rhemtulla in \cite{semi}.) However, 
a more transparent argument showing this consists in applying the positive Ping-Pong 
Lemma to the restrictions of the elements $h_m$ and $h$ to the interval $[c_m',c_n]$ 
(see \cite{harpe}, Chapter VII).

\begin{question} What are the orderable groups all of whose orderings are Conradian~?
\label{rasca}
\end{question}

\vspace{0.05cm}

Using Proposition \ref{conrad-yo}, one can provide a dynamical 
proof for the converse of Proposition \ref{notar}. The next 
proposition is originally due to Conrad \cite{conrad}.

\vspace{0.2cm}

\begin{prop} {\em Every group admitting a Conradian ordering is locally indicable.}
\label{indicar-loc}
\end{prop}

\noindent{\bf Proof.} Let \esp $\Gamma$ be a finitely generated subgroup of 
a group provided with a Conradian ordering $\preceq$. The restriction of $\preceq$ 
to $\Gamma$ is still Conradian. By Proposition \ref{conrad-yo}, the dynamical realization 
of $\Gamma$ is a group without crossed elements. By Proposition \ref{radon}, this dynamical 
realization preserves a Radon measure $\mu$. To get a non-trivial homomorphisms from $\Gamma$ 
into $(\mathbb{R},+)$, just take the translation number homomorphism associated to $\mu$. 
$\hfill\square$

\vspace{0.3cm}

For another application of Proposition \ref{conrad-yo}, recall that, by Thurston's stability 
theorem, the group $\mathrm{Diff}^1_+([0,1])$ (as well as the group of germs of $C^1$ 
diffeomorphisms at the origin) is locally indicable \cite{thurston}. As a consequence, 
these groups admit faithful actions on $[0,1]$ without crossed elements.

\begin{rem} \label{th-question} 
For interesting obstructions to $C^1$ smoothing of many actions on the 
line of some locally indicable groups (as for instance free groups), see 
\cite{cale} and references therein. However, 
we should point out that the following question remains open: does there 
exist a finitely generated locally indicable group having no faithful action by 
$C^1$ diffeomorphisms of the interval\hspace{0.05cm}?\footnote{Added in proof: 
This has been recently answered by the affirmative in \cite{yo-th}.} It is already 
interesting to know whether surface groups do admit such an action. See also 
Remark \ref{T-relative}.
\end{rem}

The following is a kind of converse to Proposition \ref{conrad-yo}.

\begin{prop} {\em Let $\Gamma$ be a subgroup of $\mathrm{Homeo}_+(\mathbb{R})$ 
without crossed elements. If $(x_n)$ is any dense sequence of points in the real 
line, then the order relation associated to this sequence is a $\mathcal{C}$-order.}
\label{conrad-yo2}
\end{prop}

\noindent{\bf Proof.} Let $f$ and $g$ be two positive elements in $\Gamma$, 
and let $\Gamma_0$ be the subgroup generated by them. Let $i \geq 0$ and 
$j \geq 0$ be the smallest indexes for which $f(x_i) \neq x_i$ and $g(x_j) \neq 
x_j$. Assume for instance that $i < j$. (The cases where $i = j$ or $i > j$ are 
similar and are left to the reader.) Let $I$ be the minimal open interval invariant 
by $\Gamma_0$ and containing $x_i$. Since $\Gamma$ does not contain crossed 
elements, there exists a (non-trivial) Radon measure $\mu$ on $I$ which is 
invariant by $\Gamma_0$. Moreover,  there is no global fixed point for 
the action of $\Gamma_0$ on it.

By the definition of $i$ and $j$, one has $f(x_n) = g(x_n) = x_n$ for all $n < i$; 
moreover, $g(x_i) = x_i$ and $f(x_i) > x_i$. Since $f$ has no fixed point on 
$I$, this easily implies that $\tau_{\mu} (f) \!>\! 0$ and $\tau_{\mu} (g) = 0$. 
Therefore, $\esp \tau_{\mu}(g^{-1}fg^2) = \tau_{\mu}(f) + \tau_{\mu}(g) = 
\tau_{\mu}(f) > 0$, \esp which implies that \esp $g^{-1}fg^2 (x) > x$ \esp 
for all \esp $x \in I$. \esp In particular, $g^{-1} f g^2$ is a positive 
element of $\Gamma$, which shows that $f g^2 \succ g$.  $\hfill\square$

\vspace{0.3cm}

As an application of the preceding equivalence, we will prove the property concerning 
positive words in $\mathcal{C}$-ordered groups announced in \S \ref{general}.

\vspace{0.1cm}

\begin{prop} {\em Let $\Gamma$ be any group with a $\mathcal{C}$-order $\precede$. 
Let $W(f,g) = f^{m_1} g^{n_1} \cdots f^{m_k} g^{n_k}$ be a word such that 
$\sum m_i > 0$ and $\sum n_i > 0$. If $f$ and $g$ are positive 
elements in $\Gamma$, then $W(f,g)$ also represents a 
positive element in $\Gamma$.}
\label{palabra-positiva}
\end{prop}

\noindent{\bf Proof.} Let us enumerate the elements of the subgroup $\Gamma_0$ 
generated by $f$ and $g$, and let us consider the dynamical realization corresponding 
to this numbering. If $\tau_{\mu}$ denotes the translation number function associated 
to an invariant Radon measure $\mu$, then one has $\tau_{\mu}(f) \geq 0$ and 
$\tau_{\mu}(g) \geq 0$. Moreover, at least one of these values is strictly 
greater than zero, as otherwise there would be global fixed points for the 
dynamical realization. Therefore, denoting $m = \sum m_i > 0$ and $n = \sum n_i > 0$, 
we have $\tau_{\mu} (W(f,g)) = m \tau_{\mu}(f) + n \tau_{\mu}(g) > 0$, and this 
implies that $W(f,g)$ is a positive element of $\Gamma$. $\hfill\square$

\vspace{0.07cm}

\begin{ex} Dehornoy's ordering is not Conradian ({\em c.f.} Example \ref{dehor1}). Indeed, 
for every \esp $i \!\in\! \{1,\ldots,n-2\}$ \esp the elements \esp 
$u \!=\! \sigma_{i} \sigma_{i+1}$ \esp and \esp $v \!=\! \sigma_{i+1}$ 
\esp are positive, but the product
\begin{eqnarray*}
u^{-1} v^{-2} u^2 v^3 \!\!
&=& \!\! \sigma_{i+1}^{-1} \sigma_{i}^{-1} \sigma_{i+1}^{-2} 
(\sigma_{i} \sigma_{i+1} \sigma_{i}) \esp \sigma_{i+1} \sigma_{i+1}^3 
= \sigma_{i+1}^{-1} \sigma_{i}^{-1} \sigma_{i+1}^{-2} (\sigma_{i+1} \sigma_{i} 
\sigma_{i+1}) \sigma_{i+1} \sigma_{i+1}^3 \\
\!\! &=& \!\! (\sigma_{i+1}^{-1} \sigma_{i}^{-1} \sigma_{i+1}^{-1}) \sigma_{i} \sigma_{i+1}^5 
= (\sigma_{i}^{-1} \sigma_{i+1}^{-1} \sigma_{i}^{-1}) \sigma_{i} \sigma_{i+1}^5 
= \sigma_{i}^{-1} \sigma_{i+1}^4
\end{eqnarray*}
is negative.
\label{dehor3}
\end{ex}

\begin{question} Let $W(f,g)$ be a word as in Proposition \ref{palabra-positiva}. Assume that for 
an ordering $\precede$ on a group $\Gamma$ one has $W(f,g) \succ id$ for all positive elements $f,g$. 
Under what conditions on $W$ one can ensure that \esp $\precede$ \esp is a 
$\mathcal{C}$-order\hspace{0.05cm}? 
(The reader may easily check that this is for instance the case of 
$W(f,g) = f^{-1} g^{-1} fgfg$.)
\end{question}

\vspace{0.05cm}

For future reference, we give a slight modification of Proposition \ref{conrad-yo2} which 
involves subgroups of countable groups endowed with a non necessarily Conradian order.

\vspace{0.15cm}

\begin{prop} {\em Let $\precede$ be an ordering on a countable group \esp $\Gamma$, 
and let \esp $\Gamma_*$ be a subgroup of \esp $\Gamma$. Let $(g_n)_{n \geq 0}$ be any 
numbering of the elements of \esp $\Gamma$ starting with \esp $g_0 \! = \! id$. \esp 
Assume that, for the corresponding dynamical realization of $\precede$, there exists an 
interval $]\alpha,\beta[$ containing the origin and which is globally fixed by \esp 
$\Gamma_*$. If the restriction of \esp $\Gamma_*$ to $]\alpha,\beta[$ has no crossed 
elements, then the order $\precede$ restricted to $\Gamma_*$ is Conradian.}
\label{conradito}
\end{prop}

\noindent{\bf Proof.} Since for each $g\!\in\!\Gamma$ one has \esp $t(g) \!=\!g(0)$, 
\esp for every $g \! \in \! \Gamma_*$ the point $t(g)$ must belong to $]\alpha,\beta[$. 
Moreover, an element $g \! \in \! \Gamma$ is positive if and only if $g(0) \!>\! 0$. 
With these facts in mind one may proceed to the proof as in the case of Proposition 
\ref{conrad-yo2}. We leave the details to the reader. $\hfill\square$

\vspace{0.45cm}

We do not know whether there exists an analogous extension (or modification) 
of Proposition \ref{conrad-yo}. However, in the next Section we will show 
such an statement under a convexity hypothesis (see Lemma \ref{salida-digna}), 
and this will be enough for our purposes. 

We close this Section with a useful definition.

\vspace{0.05cm}

\begin{defn} Two orientation-preserving homeomorphisms $f,g$ of the real line are said 
to be {\em in transversal position} on an interval $[a,b] \subset \mathbb{R}$ if \esp 
$f(x) \!<\! x$ \esp for all \esp $x \! \in ]a,b]$ \esp and \esp $f(a) \!=\! a$, \esp 
and \esp $g(x) \!>\! x$ \esp for all $x \! \in [a,b[$ and \esp $g(b) \!=\! b$. 
\end{defn}

\vspace{0.05cm}

The reader can easily check that some of the arguments used in the 
proof of Proposition \ref{conrad-yo} actually show the following.

\vspace{0.05cm}

\begin{prop} {\em A subgroup of $\mathrm{Homeo}_+(\mathbb{R})$ has no crossed 
elements if and only if it does not contain elements in transversal position.}
\label{transversal}
\end{prop}

%%%%%%%%%%%%%%%%%%%%%%%%%%%%%%%%%%%%%%%%%%%%%%%%%%%%%%%%%%%%%%%%%%%%%%%%%%%%%%%%%%%%%%%%%

\subsubsection{The Conradian soul of an order}
\label{conrad-soul}

\hspace{0.35cm} Let $\precede$ be a left-invariant total order on a (non necessarily 
countable) group $\Gamma$. A subgroup $\Gamma_*$ of $\Gamma$ is said to be {\em convex} 
with respect to $\precede$ (or just $\precede$-convex) if, for all $f \prec g$ in $\Gamma_*$, 
every element $h \!\in\! \Gamma$ satisfying $f \prec h \prec g$ belongs to $\Gamma_*$. 
Equivalently, $\Gamma_*$ is convex if, for each $f \succ id$ in $\Gamma_{*}$, every 
$g \in \Gamma$ such that $id \prec g \prec f$ belongs to $\Gamma_*$.

\vspace{0.02cm}

\begin{ex} From the definition one easily checks that, 
for each $n \!\geq\! 2$ and each $j \!\in\! \{1,\ldots,n-1\}$, the subgroup 
\esp $\langle \sigma_j,\ldots,\sigma_{n-1}\rangle \sim B_{n-j+1}$ \esp of $B_n$ 
is convex with respect to Dehornoy's ordering ({\em c.f.} Example \ref{dehor1}).
\label{dehor4}
\end{ex}

\vspace{0.02cm}

Note that for every ordering $\preceq$ on a group $\Gamma$, the family of $\preceq$-convex 
subgroups coincides with that of $\bar{\preceq}$-convex ones 
({\em c.f.} Remark \ref{la-primera}). 
A more important (and  
also easy to check) fact is that this family is linearly ordered (by inclusion). 
More precisely, if $\Gamma_0$ and $\Gamma_1$ are $\precede$-convex, then either \esp 
$\Gamma_0 \!\subset\! \Gamma_1$ \esp or \esp $\Gamma_1 \!\subset\! \Gamma_0$. In particular, 
the union and the intersection of any family of convex subgroups is a convex subgroup.

\begin{rem} Let $\preceq$ be an ordering on a group $\Gamma$. For each non-trivial 
element $g \in \Gamma$ one may define $\Gamma_g$ (resp. $\Gamma^g$) as the largest 
(resp. smallest) convex subgroup which does not contain $g$ (resp. which contains $g$). 
It turns out that $\preceq$ is Conradian if and only if for each $g \neq id$ the group 
$\Gamma_g$ is normal in $\Gamma^g$ and the order on $\Gamma^g / \Gamma_g$ induced by 
$\preceq$ is Archimedean (and in particular the quotient $\Gamma^g / \Gamma_g$ is 
torsion-free Abelian), see \cite{botto,glass,koko}. The reader should note a close 
relationship between this characterization and the dynamical one given in the previous 
Section. For instance, a good exercise is to prove Proposition \ref{palabra-positiva} 
using the characterization of $\mathcal{C}$-orders in terms of convex subgroups. 
(See \cite{leslie} for more on this.)
\label{convex-car}
\end{rem}

We will say that a subgroup $\Gamma_*$ of $\Gamma$ is Conradian 
with respect to an ordering $\precede$ 
on $\Gamma$ (or just $\precede$-Conradian) if the restriction of $\precede$ to $\Gamma_*$ 
is a $\mathcal{C}$-order. Note that if $\{\Gamma_i\}_{i \in \mathcal{I}}$ is a linearly 
ordered family of $\precede$-Conradian subgroups of $\Gamma$, then the union 
$\Gamma_{*} \!=\! \cup_{i \in \mathcal{I}} \Gamma_i$ is still $\precede$-Conradian. 
Therefore, the following definition makes sense.

\vspace{0.04cm} 

\begin{defn} The {\em Conradian soul} of $\Gamma$ with respect to $\precede$ 
(or just the $\preceq$-Conradian soul of \esp $\Gamma$) is the maximal 
subgroup $\Gamma_{\preceq}^{\esp \esp c}$ of $\Gamma$ which is 
simultaneously $\precede$-convex and $\precede$-Conradian.
\end{defn}

%\vspace{0.02cm}

\begin{ex} We will see in Example \ref{dehor8} that the Conradian soul of $B_n$ with 
respect to Dehornoy's ordering is the cyclic subgroup generated by $\sigma_{n-1}$ 
({\em c.f.} Examples \ref{dehor1} and \ref{dehor3}).
\label{dehor5}
\end{ex}

For the case where $\Gamma$ is countable, the Conradian soul has a very simple 
dynamical description. Indeed, fix a numbering $(g_n)_{n \geq 0}$ of $\Gamma$ 
such that $g_0 \!=\! id$, and for the corresponding dynamical realization define 
$$\alpha = \sup \{ b < 0: \esp \mbox{ there exist } f,g \mbox{ in } 
\Gamma \mbox{ such that } f,g \mbox{ are crossed on } ]a,b[ \},$$
$$\beta = \inf \{ a > 0: \esp \mbox{ there exist } f,g \mbox{ in } 
\Gamma \mbox{ such that } f,g \mbox{ are crossed on } ]a,b[ \},$$
where we let \esp $\alpha \!=\! -\infty$ \esp (resp. \esp $\beta \!=\! +\infty$) 
\esp if the corresponding set of $b$'s (resp. $a$'s) in $\mathbb{R}$ is empty. 
Note that the arguments of the proof of Proposition \ref{conrad-yo} show that, 
in the previous definitions, we can replace ``are crossed on $]a,b[$'' by 
``are in transversal position on $[a,b]$'' without changing the values of 
$\alpha$ and $\beta$. The following lemma will be implicitly used in what 
follows, and helps to understand the situation better.

\vspace{0.15cm}

\begin{lem} {\em The equality \esp $\alpha \!=\! - \infty$ \esp 
holds if and \esp and only if \esp $\beta \!=\! +\infty$. Similarly, 
one has \esp $\alpha \!<\! 0$ \esp if and only if $\beta \!>\! 0$.}
\end{lem}

\noindent{\bf Proof.} Assume that $\beta \!<\! +\infty$. Then there exists $f,g$ which 
are in transversal position on some interval $[a,b]$ satisfying $a \!\geq\! \beta$. Let 
$h \!\in\! \Gamma$ be such that $h(b) \!<\! 0$. Then the elements $hfh^{-1}$ and $hgh^{-1}$ 
are in transversal position on $[h(a),h(b)]$, and since $h(b) \!<\! 0$ this shows that 
$\alpha \!>\! -\infty$. A similar argument shows that the condition 
$\alpha \!>\! -\infty$ implies $\beta \!<\! +\infty$.

Now suppose that $\beta \!=\! 0$. Then given any $h \!\succ\! id$ there are elements 
$f,g$ which are in transversal position on an interval $[a,b]$ satisfying $a \!\in ]0,t(h)[$. 
After conjugacy by $f^k$ for $k \!\in\! \mathbb{N}$ large enough, we may suppose that the point 
$b$ also belongs to $]a,t(h)[$. If this is the case, the elements $h^{-1}fh$ and $h^{-1}gh$ are 
in transversal position on $[h^{-1}(a),h^{-1}(b)] \!\subset \esp \esp \esp ]t(h^{-1}),0[$. Since 
this construction can be performed for any positive element $h \!\in\! \Gamma$, this implies 
that $\alpha\!=\!0$. A similar argument shows that, if $\alpha \!=\! 0$, then $\beta \!=\! 0$. 
$\hfill\square$

\vspace{0.35cm}

Note that the equalities \esp $\alpha = -\infty$ \esp and \esp $\beta = +\infty$ \esp hold 
if and only if $\con = \Gamma$, \esp that is, if $\precede$ is a $\mathcal{C}$-order.  

\vspace{0.15cm}

\begin{prop} {\em With the previous notations, the $\precede$-Conradian soul of 
\esp \esp $\Gamma$ coincides with the stabilizer of the interval $]\alpha,\beta[$.} 
\label{reca}
\end{prop}

\vspace{0.15cm}

To prove this proposition, we will need the following general lemma.

\vspace{0.15cm}

\begin{lem} {\em Let $\Gamma$ be a countable group, and let $(g_n)_{n \geq 0}$ be a numbering 
of its elements starting with $g_0 \!=\! id$. Let us consider the dynamical realization 
associated to an ordering $\precede$ on $\Gamma$ and corresponding to this 
numbering. Suppose that \esp $\Gamma_*$ is a convex subgroup, and that $]\alpha,\beta[$ is an 
interval which is fixed by $\Gamma_*$ and which does not contain any global fixed point 
of $\Gamma_*$. If the restriction of $\Gamma_*$ to $]\alpha,\beta[$ has crossed elements 
and $]\alpha,\beta[$ contains the origin, then $\Gamma_*$ is not $\precede$-Conradian.}
\label{salida-digna}
\end{lem}

\noindent{\bf Proof.} We would like to use similar arguments as those of the proof 
of Proposition \ref{conrad-yo}. Note that those arguments still apply and involve 
only elements of $\Gamma_*$, except perhaps the one concerning the element $g_i$. 
More precisely, we need to ensure that an element $g_i \!\in\! \Gamma$ such that 
$t(g_i)$ is in \esp $]a,b[ \subset ]\alpha,\beta[$ \esp actually belongs to 
$\Gamma_*$. For this we use the convexity hypothesis. Indeed, since the 
supermom of the orbit by $\Gamma_*$ of the origin is a point which is globally 
fixed by $\Gamma_*$, it must coincide with $\beta$. In particular, there exists 
$h_1 \! \in \! \Gamma_*$ such that $h_1(0) \!>\! t(g_i)$. In an analogous way, 
one obtains $h_2(0) \!<\! t(g_i)$ for some $h_2 \!\in\! \Gamma_*$. Now since 
$h_i(0) \!=\! t(h_i)$, this gives \esp $h_2 \!\prec\! g_i \!\prec\! h_1$. 
\esp By the convexity of $\Gamma_*$, this implies that $g_i$ is 
contained in $\Gamma_*$, thus finishing the proof. $\hfill\square$

\vspace{0.5cm}

Now we can pass to the proof of Proposition \ref{reca}. Denote by $\Gamma_*$ 
the stabilizer of $]\alpha,\beta[$. We need to verify several facts.

\vspace{0.25cm}

\noindent{\underbar{Claim 1.}} The group $\Gamma_*$ 
is a $\precede$-convex subgroup of $\Gamma$.

\vspace{0.25cm}

We first claim that there is no element $h \! \in \! \Gamma$ 
sending $\alpha$ or $\beta$ into $]\alpha,\beta[$. Indeed, assume that $h(\beta)$ 
belongs to $]\alpha,\beta[$. (The case $h(\alpha) \! \in \esp ]\alpha,\beta[$ 
is analogous.) If $h(\beta)$ is in $[0,\beta[$, then let $\varepsilon > 0$ be such 
that $h([\beta,\beta+\varepsilon]) \! \subset [0,\beta[$. By the definition of 
$\beta$, there exist $a < b$ and elements $f,g$ in $\Gamma$ such that 
$\beta \leq a < \beta + \varepsilon$ and such that $f,g$ are in transversal 
position on $[a,b]$. Changing (if necessary) $g$ by $f^ngf^{-n}$ for $n$ large 
enough, we may assume that $[a,b]$ is contained in $[\beta,\beta+\varepsilon[$;  
then changing $f$ by $g^k f g^{-k}$ for $k$ large enough, we may suppose 
that $[a,b]$ is actually contained in $]\beta,\beta + \varepsilon[$. Now the 
elements $hfh^{-1}$ and $hgh^{-1}$ are in transversal position on $[h(a),h(b)]$, 
and since $0 \!<\! h(a) \!<\! \beta$, this contradicts the definition of 
$\beta$. 

When $h(\beta)$ is in $]\alpha,0[$, the situation is slightly more 
complicated. Fix $\varepsilon \! > \! 0$ such that $h([\beta,\beta+\varepsilon]) 
\! \subset ]\alpha,0[$. Again by the definition of $\beta$, there exist $a < b$ 
and elements $f,g$ in $\Gamma$ such that $\beta \leq a < \beta + \varepsilon$ 
and such that $f,g$ are crossed on $]a,b[$, where for concreteness we assume 
that $Fix (f) \cap [a,b] = \{a,b\}$ and $f(x) < x$ for all $x \! \in ]a,b[$. 
Now refer to Figure 1, where for \esp $m >\!\!> n$ \esp big enough the 
elements $h_n$ and $h_m$ are in transversal position on the interval $[c_m',c_n]$. 
Fix  $k \in \mathbb{N}$ large enough in such a way $f^{k}(c_n)$ is near to $a$ so that 
$h( f^k(c_n) ) \!\in\! [h(\beta),0[$. Then the elements $h f^k h_n f^{-k} h^{-1}$ 
and $h f^k h_m f^{-k} h^{-1}$ are in transversal position on the interval \esp 
$[hf^k (c_m'),hf^k (c_n)]$, \esp and since $\alpha \! < \! h(\beta) \! 
< \! hf^k(c_n) \! <  \! 0$, this contradicts the definition of $\alpha$.

Now to conclude the proof of the $\precede$-convexity of $\Gamma_*$, let 
$h \! \in \! \Gamma$ be such that $f \! \prec \! h \! \prec \! g$ for 
some elements $f,g$ in $\Gamma_*$. We then have 
$\alpha < t(f) < t(h) < t(g) < \beta$, and therefore 
$\alpha < h(0) < \beta$. Since both $h$ and $h^{-1}$ do not send 
neither $\alpha$ nor $\beta$ into $]\alpha,\beta[$, this easily implies that 
$h (\alpha) \!=\! \alpha$ and $h (\beta) \!=\! \beta$. 
Therefore, $h$ belongs to $\Gamma_*$. 

\vspace{0.3cm}

\noindent{\underbar{Claim 2.}} The restriction of $\precede$ to $\Gamma_*$ is Conradian.

\vspace{0.25cm}

This follows as a direct application of Proposition \ref{conradito}.

\vspace{0.3cm}

\noindent{\underbar{Claim 3.}} The group $\Gamma_*$ is a maximal subgroup for 
the property of being simultaneously $\precede$-convex and $\precede$-Conradian.

\vspace{0.25cm} 

Let $\hat{\Gamma}$ be a convex subgroup of $\Gamma$ strictly containing $\Gamma_*$. 
Fix a positive element $h \! \in \! \hat{\Gamma} \setminus \Gamma_*$. One has 
$h(\alpha) \geq \beta$, and therefore $h (0) \!>\! \beta$. Let $\varepsilon 
\!=\! h(0) \!-\! \beta$. As in the proof of Claim 1, there exist $f,g$ 
in $\Gamma$ which are in transversal position on an interval $[a,b]$ such 
that $[a,b] \! \subset  ]\beta,\beta + \varepsilon[$. We then have 
$$t(h) = h(0) = \beta \!+\! \varepsilon > t(f) \qquad \mbox{ and } 
\qquad t(h) > t(f^{-1}),$$
and similarly \esp $t(h) > t(g)$ \esp and \esp $t(h) > t(g^{-1})$. \esp 
From the $\precede$-convexity of $\hat{\Gamma}$ one easily deduces from this that both 
elements \esp $f$ \esp and \esp $g$ \esp belong to $\hat{\Gamma}$. Now the first global 
fixed point of $\hat{\Gamma}$ immediately to the right of the origin is to the right 
of \esp $h(0) \!\geq\! b$. \esp Therefore, by Lemma \ref{salida-digna}, the subgroup 
$\hat{\Gamma}$ is not $\precede$-Conradian. This proves Claim 3 and finishes 
the proof of Proposition \ref{reca}.

\vspace{0.15cm}

\begin{rem} The reader should have no problem in adapting some of the arguments 
above to prove that, if $\Gamma$ is infinite, then $\con$ is non-trivial 
if and only if \esp $\alpha \!<\! 0$, \esp which is equivalent to 
\esp $\beta \!>\! 0$.
\end{rem}

%%%%%%%%%%%%%%%%%%%%%%%%%%%%%%%%%%%%%%%%%%%%%%%%%%%%%%%%%%%%%%%%%%%%%%%%%%%%%%%%%%%%%%%%

\subsubsection{Extensions of orders and stability of Conradian souls}
\label{extension}

\hspace{0.35cm} Let $\preceq$ be an ordering on a group $\Gamma$, 
and let $\Gamma_*$ be a $\precede$-convex subgroup of $\Gamma$. Let $\precede_*$ be 
any (total and left-invariant) order on $\Gamma_*$. The {\em extension of $\precede_*$ 
by $\precede$} is the order relation $\precede'$ on $\Gamma$ whose positive cone is 
$(P_{\precede}^+ \setminus \Gamma_*) \cup P^+_{\precede_*}$. It is easy to check 
that $\precede'$ is also a left-invariant total order relation, and that $\Gamma_*$ 
remains convex in $\Gamma$ (that is, it is a $\precede'$-convex subgroup of $\Gamma$).

\vspace{0.05cm}

\begin{rem} With the notations above, one easily checks that the family of 
$\preceq'$-convex subgroups of $\Gamma$ is formed by the $\preceq_*$-convex 
subgroups of $\Gamma_*$ and the $\preceq$-convex of $\Gamma$ which contain $\Gamma_*$. 
\label{convex-by-extension}
\end{rem}

\vspace{0.05cm}

The extension procedure is a classical and useful technique which allows 
for instance to give an alternative approach to the orderings on 
braid groups introduced by Dubrovina and Dubrovin in \cite{dub}.

\begin{ex} Since the cyclic subgroup \esp $\langle \sigma_2 \rangle$ \esp is convex in $B_3$ with 
respect to Dehornoy's ordering $\preceq_D$ \esp\esp ({\em c.f.} Example \ref{dehor4}), one can define 
the order $\preceq_3$ on $B_3$ as being the extension by $\preceq_D$ of the restriction to \esp 
$\langle \sigma_2 \rangle$ \esp of $\bar{\preceq}_D$ \esp ({\em c.f.} Remark \ref{la-primera}). We claim 
that the positive cone of $\preceq_3$ is generated by the elements \esp $u_1 \!=\! \sigma_1 \sigma_2$ 
\esp and \esp $u_2 \!=\! \sigma_2^{-1}$. \esp Indeed, by definition these elements are positive 
with respect to $\preceq_3$, and therefore it suffices to show that for every \esp $u \neq id$ 
\esp in $B_3$ either $u$ or $u^{-1}$ belongs to the semigroup \esp $\langle u_1,u_2 \rangle^+$ 
\esp generated by $u_1$ and $u_2$. Now if $u$ or $u^{-1}$ is $\sigma_2$-positive for Dehornoy's 
ordering, then there exists an integer $m \neq 0$ such that \esp $u = \sigma_2^m = u_2^{-m}$, \esp 
and therefore \esp $u \in \langle u_2 \rangle^+ \subset \langle u_1,u_2 \rangle^+$ \esp if $m < 0$ 
and \esp $u^{-1} \in \langle u_2 \rangle^+ \subset \langle u_1,u_2 \rangle^+$ \esp if $m > 0$. If 
$u$ is $\sigma_1$-positive, then for a certain choice of integers $m_1'', \ldots,m_{k''+1}''$ one has 
$$u = \sigma_2^{m_1''} \sigma_1 \sigma_2^{m_2''} \sigma_1 
\cdots \sigma_2^{m_{k''}''} \sigma_1 \sigma_2^{m_{k''+1}''}.$$
Using the identity \esp \esp $\sigma_1 \!=\! u_1u_2$, \esp \esp 
this allows us to writte $u$ in the form 
$$u = u_2^{m_1'} u_1 u_2^{m_2'} u_1 \ldots u_2^{m_{k'}'} u_1 u_2^{m_{k'+1}'}$$
for some integers $m_1',\ldots,m_{k'+1}'$. Now using several times the (easy to 
check) identity \esp $u_2 u_1^2 u_2 =u_1$, \esp one may express $u$ as a product 
$$u = u_2^{m_1} u_1 u_2^{m_2} u_1 \ldots u_2^{m_{k}} u_1 u_2^{m_{k+1}}$$
in which all the exponents $m_i$ are non negative, and this shows that $u$ belongs to 
\esp $\langle u_1,u_2 \rangle^+$. \esp Finally, if $u^{-1}$ is $\sigma_1$-positive 
then $u^{-1}$ belongs to \esp $\langle u_1,u_2 \rangle^+$. 
\label{dehor6}
\end{ex}

\begin{ex}
The generalization of the previous example to all braid groups proceeds inductively as follows. Let 
us see $B_{n-1} = \langle \tilde{\sigma}_1, \ldots, \tilde{\sigma}_{n-2} \rangle$ as a subgroup 
of \esp $B_n = \langle \sigma_1, \ldots, \sigma_{n-1} \rangle$ \esp via the monomorphism \esp 
$\tilde{\sigma}_i \mapsto \sigma_{i+1}$. \esp Via this identification, we obtain from $\preceq_{n-1}$ 
an order on $\langle \sigma_2,\ldots,\sigma_{n-1} \rangle \subset B_n$, which we still denote by 
$\preceq_{n-1}$. We then let $\preceq_{n}$ be the extension of $\bar{\preceq}_{n-1}$ by the 
Dehornoy's ordering $\preceq_D$. Once again, an important property of $\preceq_n$ is that its 
positive cone is finitely generated as a semigroup (and therefore, by Proposition \ref{paja}, 
the ordering $\preceq_n$ is an isolated point of the space of orderings of $B_n$.) More 
precisely, letting 
$$v_1 = \sigma_1 \sigma_2 \cdots \sigma_{n-1}, \quad 
v_2 = \sigma_2 \sigma_3 \cdots \sigma_{n-1}, \quad \ldots \ldots \hspace{0.2cm} , 
\quad v_{n-2} = \sigma_{n-2}\sigma_{n-1}, \quad v_{n-1} = \sigma_{n-1},$$
and \esp $u_i \!=\! v_i^{(-1)^{i-1}}$ \esp (where $i \!\in\! \{1,\ldots,n-1\}$), 
the semigroup  $P_{\preceq_n}^{+}$ is generated by the elements \esp 
$u_1,\ldots,u_{n-1}$. \esp To check this, one proceeds by induction 
using (as in the case $n \!=\! 3$) the remarkable identities
$$(u_2 u_3^{-1} \cdots u_{n-1}^{(-1)^{n-1}}) u_1^{n-1} 
(u_2 u_3^{-1} \cdots u_{n-1}^{(-1)^{n-1}}) = u_1 \quad \mbox{ and } 
\quad (u_2 u_3^{-1} \cdots u_{n-1}^{(-1)^{n-1}})^2 = u_2^{n-1}.$$
For the sake of clarity, we will denote by $\preceq_{DD}$ the orderings 
constructed above (called Dubrovina-Dubrovin's orderings in the Introduction).
\label{dehor7}
\end{ex}

For countable groups, the extension procedure can be described in pure dynamical 
terms. Roughly, it corresponds to consider the dynamical realization of $\precede$, then to 
change the action of $\Gamma_*$ on the smallest interval $]\alpha,\beta[$ containing the origin 
and which is fixed by $\Gamma_*$ by (a conjugate of) the action associated to a dynamical 
realization of $\precede_*$, and then to extend the new action to the whole group $\Gamma$ 
in an equivariant way. This approach 
naturally leads to the following stability type property for Conradian souls: if \esp 
$\Gamma_*$ coincides with the $\precede$-Conradian soul of $\Gamma$ and $\precede_*$ 
is a $\mathcal{C}$-order on $\Gamma_*$, then $\Gamma_*$ also corresponds to the 
$\precede'$-Conradian soul  of $\Gamma$. However, the algebraic presentation of 
the extension operation being more concise, it allows to give a short proof 
of this fact which also covers the case of uncountable orderable groups.

\vspace{0.15cm}

\begin{lem} {\em Let $\preceq$ be an ordering on a group \esp \esp $\Gamma$, \esp and 
let $\preceq_*$ be any left-invariant total order on the $\preceq$-Conradian soul $\con$ of 
\esp \esp $\Gamma$ which is still a $\mathcal{C}$-order. If $\preceq'$ denotes the extension 
of $\preceq_*$ by $\preceq$, then the $\preceq'$-Conradian soul of \esp \esp $\Gamma$ 
coincides with $\con$.}
\label{esta-con}
\end{lem}

\noindent{\bf Proof.} Since $\con$ is a convex and Conradian subgroup of $\Gamma$ with 
respect to $\preceq'$, we just need to check the maximality property. So let $\Gamma_*$ be 
any $\preceq'$-convex subgroup of $\Gamma$ strictly containing $\con$. We first claim that 
$\Gamma_*$ is also $\preceq$-convex. Indeed, assume that \esp $f \prec h \prec g$ \esp for 
some $f,g$ in $\Gamma_*$ and $h \in \Gamma$. If either \esp $f^{-1}h$ \esp 
or \esp $g^{-1}h$ \esp belongs to 
$\con$ then, since $\con$ is contained in $\Gamma_*$ and \esp $h = f(f^{-1}h) = g(g^{-1}h)$, 
\esp the element $h$ belongs to $\Gamma_*$. If neither \esp $f^{-1}h$ \esp 
nor \esp $g^{-1}h$ \esp does belong to 
$\con$ then, since \esp $id \prec f^{-1}h$ \esp and \esp 
$g^{-1}h \prec id$, \esp one has \esp 
$id \prec' f^{-1}h$ \esp and \esp $g^{-1}h \prec' id$, \esp 
that is, \esp $f \prec' h \prec' g$. \esp By the $\preceq'$ convexity 
of $\Gamma_*$, this still implies that $h$ is contained in $\Gamma_*$, thus showing the 
$\preceq$-convexity of $\Gamma_*$.

Since $\Gamma_*$ is $\preceq$-convex and strictly contains $\con$, there exist positive 
elements $f,g$ in $\Gamma_*$ such that \esp $fg^n \preceq g$ \esp for all $n \!\in\! \mathbb{N}$. 
We claim that $g$ does not belong to $\con$. Indeed, if not then one has $f \notin \con$, and 
therefore $f^{-1} \prec g$, that is, $fg \succ id$. Again, since $fg \notin \con$, 
this implies that $fg \succ g$, which contradicts our choice.

We now claim that, for every $n \geq 0$, the element $g^{-1}fg^{n}$ does not belong to $\con$. 
Indeed, since $g$ is a positive element not contained in $\con$, if \esp $g^{-1}fg^n$ \esp 
is in $\con$ then \esp $g \succ (g^{-1} f g^n)^{-1}$, \esp and therefore \esp 
$g^{-1} f g^{n+1} \succ id$, \esp contradicting again our choice.

Now we remark that, independently if $f$ does belong or not to $\con$, the element $h \!=\! fg$ 
(is positive and) is not contained in $\con$. Therefore, both $g$ and $h$ are still positive 
with respect to the ordering $\preceq'$. Moreover, since \esp $g^{-1} f g^n \preceq id$ \esp and \esp 
$g^{-1} f g^{n} \notin \con$ \esp for all \esp $n \geq 0$, \esp one necessarily has \esp 
$g^{-1} h g^n \prec' id$ \esp for all $n \geq 0$. In particular, $\Gamma_*$ is not a 
$\preceq'$-Conradian subgroup of $\Gamma$. Since this is true for any $\preceq'$-convex 
subgroup of $\Gamma$ strictly containing $\con$, this shows that the $\preceq'$-Conradian 
soul of $\Gamma$ coincides with $\con$. $\hfill\square$

\vspace{0.1cm}

\begin{ex} The only $\preceq_{n}$-convex subgroups of 
\esp $B_{n}$ \esp are \esp $B^1\!=\!\{id\}$, \esp 
$B^2\!=\!\langle u_{n-1} \rangle \!=\! \langle \sigma_{n-1} \rangle$, \esp 
$B^{3}\!=\!\langle u_{n-2},u_{n-1} \rangle \!=\! \langle \sigma_{n-2},\sigma_{n-1} 
\rangle$, \esp \ldots, \esp $B^{n-1}\!=\!\langle u_2,\ldots,u_{n-1} \rangle\!=\!
\langle \sigma_2,\ldots,\sigma_{n-1} \rangle$ \esp and \esp $B^{n}\!=\!B_{n}$. 
Indeed, suppose that there exists a $\preceq_{n}$-convex subgroup $B$ of 
$B_{n}$ such that \esp $B^{i} \subsetneq B \subsetneq B^{i+1}$ for some \esp 
$i \!\in\! \{1,\ldots,n-1\}$. \esp Let $\preceq^1$, $\preceq^2$, and $\preceq^3$, 
be the orderings respectively defined on $B^{i}$, $B$, and $B_n$, by:

\vspace{0.1cm}

\noindent -- $\preceq^1$ is the restriction of $\preceq_{n}$ to $B^{i}$,

\vspace{0.1cm}

\noindent -- $\preceq^2$ is the extension of $\preceq^1$ by the restriction of 
$\bar{\preceq}_{n}$ to $B$,

\vspace{0.1cm}

\noindent -- $\preceq^3$ is the extension of $\preceq^2$ by $\preceq_{n}$.

\vspace{0.1cm}

\noindent The order $\preceq^3$ is different from $\preceq_{n}$ (the $\preceq_{n}$-negative 
elements in $B \setminus B^i$ are $\preceq^3$-positive), but its positive cone still 
contains the elements \esp $u_1,\ldots,u_i,u_{i+1},\ldots,u_{n-1}$. \esp Nevertheless, this is 
impossible, since these elements generate the positive cone of $\preceq_{n}$.

Note that, by Remark \ref{convex-by-extension}, the $\preceq_D$-convex subgroups 
of $B_n$ coincide with the $\preceq_{n}$-convex subgroups listed above.
\label{dehor9}
\end{ex}

%\vspace{0.05cm}

\begin{ex} Since the smallest $\preceq$-convex subgroup strictly containing 
$\langle \sigma_{n-1} \rangle$ is $\langle \sigma_{n-2}, \sigma_{n-1} \rangle$, and 
since the restriction of $\preceq_D$ to $\langle \sigma_{n-2}, \sigma_{n-1} \rangle$ 
is not Conradian ({\em c.f.} Example \ref{dehor3}), the Conradian soul of $B_n$ with 
respect to Dehornoy's ordering is the infinite cyclic subgroup generated by $\sigma_{n-1}$.
\label{dehor8}
\end{ex}

%\vspace{0.05cm}

\begin{rem} 
In \cite{paper-in-l'enseignement}, Short and Wiest study the orderings on braid groups (and 
more generally on some mapping class groups) which arise from Nielsen's geometrical methods. 
They define two different families of such orderings, namely those of {\em finite} and {\em 
infinite type}. They distinguish these families by showing that the former ones are {\em 
discrete} (that is, there exists a minimal positive element for them), and the latter 
ones are non discrete. (Dehornoy's ordering belongs to the first family.) It would be nice 
to pursue a little bit on this point for explicitly determining the Conradian soul in 
each case.\footnote{Added in proof: This has been recently done in \cite{NW}.} 
\label{primero-l'enseignement}
\end{rem}

%%%%%%%%%%%%%%%%%%%%%%%%%%%%%%%%%%%%%%%%%%%%%%%%%%%%%%%%%%%%%%%%%%%%%%%%%%%%%%%%%%%%%%%%%%%%%%

\subsection{Right-recurrent orders}
\label{a-derecha}

\hspace{0.35cm} A left-invariant total order relation $\precede$ on a group 
$\Gamma$ is {\em right-recurrent} if for all positive elements $f,g$ there 
exists $n \in \mathbb{N}$ such that $g f^n \succ f^n$. Clearly, every such 
order satisfies the Conrad property, but the converse is not true. Remark 
that both the sets of $\mathcal{C}$-orders and right-recurrent orders 
are invariant under the action of \esp $\Gamma$ by conjugacy.

The property of right-recurrence for left-invariant orders is not so clear as the Conradian property 
or the bi-invariance. For instance, as the following example shows, there is no analogue of 
neither Proposition \ref{lesla} nor Proposition \ref{conrad-cerrado} for right-recurrent orders.

\begin{ex} Let $f$ be the translation $x \mapsto x+1$, and let $g$ be any orientation-preserving 
homeomorphism of the unit interval such that $g(x) > x$ for all $x \! \in ]0,1[$. Fix an increasing 
sequence $(n_i)$ of non negative integers such that $n_0 = 0$ and such that $n_{2k+1} \!-\! n_{2k}$ 
goes to infinite with $k$. Extend $g$ into a homeomorphism of the whole line by defining, 
for $n \in \mathbb{Z}$ and $x \in [n,n+1]$,
$$g(x) = \left \{ \begin{array} {l} 
f^{n}gf^{-n}(x) \hspace{1cm} \mbox{if } n = n_{2k},\\ 
f^n g^{-1} f^{-n} (x) \hspace{0.62cm} \mbox{if } n = n_{2k + 1},\\ 
x \hspace{2.6cm} \mbox{otherwise.} \end{array} \right.$$
It is not difficult to check that the group $\Gamma$ generated by $f$ and $g$ is isomorphic 
to the wreath product $\mathbb{Z} \wr \mathbb{Z}$. For each $k$ let $\precede_k$ be the 
order relation on $\Gamma$ defined by \esp $h_1 \! \prec_k \! h_2$ \esp if and only if 
the minimum integer $i \geq n_{2k}$ for which \esp \esp $h_1 (i +1/2) \neq h_2 (i + 1/2)$ 
\esp \esp is such that \esp $h_1 (i +1/2) < h_2 (i + 1/2)$. \esp One can easily show that 
each $\precede_k$ is total, left-invariant, and right-recurrent. (Note that $\precede_k$ 
coincides with the image of $\precede_0$ by $f^{-n_{2k}}$.)  Nevertheless, 
no accumulation point $\precede$ of the sequence of orders $\precede_k$ is 
right-recurrent. Indeed, the elements $f$ and $g$ are positive for all the 
orders $\precede_k$. On the other hand, one has \esp $gf^{n} \prec_k f^{n}$ 
\esp for all \esp $n \! \in \! \{1,\ldots,n_{2k+1} - n_{2k}\}$, \esp and passing 
to the limit this gives \esp $gf^n \prec f^n$ \esp for all $n \in \mathbb{N}$.
\end{ex}

\vspace{0.1cm}

Although the set of right-recurrent orders is contained in the set of $\mathcal{C}$-orders, 
it is not necessarily dense therein. (See however Question \ref{densidad}.) Indeed, according 
to \cite[Example 4.6]{witte}, if $F$ is a finite index free subgroup of $\mathrm{SL}(2,\mathbb{Z})$, 
then the group $\Gamma \!=\! F \ltimes \mathbb{Z}^2$ admits no right-recurrent order. However, $\Gamma$ 
is locally indicable, and therefore by Proposition \ref{notar} it admits a $\mathcal{C}$-order. 
(By Proposition \ref{conrad-yo}, it also admits a faithful action on the interval without 
crossed elements.) 

\begin{rem} The group $\Gamma$ above satisfies the relative Kazhdan's property (T) with 
respect to the {\em normal} subgroup $\mathbb{Z}^2$. By \cite[Th\'eor\`eme A]{comment}, 
for no $\varepsilon > 1/2$ this group can act faithfully 
by $C^{3/2 + \varepsilon}$ diffeomorphisms 
of the interval.\footnote{Added in Proof: This has been recently extended in 
\cite{thurston} to actions by $C^1$ diffeomorphisms.}
\label{T-relative}
\end{rem}

\vspace{0.01cm}

\begin{question} Is the property of admitting a rigth-recurrent order a ``local" 
property~? (See the comments after the proof of Proposition \ref{bur}.)
\label{recurrente-local}
\end{question}

\begin{question} What are the orderable groups all of whose orderings
are right-recurrent~? (This should be compared with Question \ref{rasca} 
as well as Tararin's theorem in \S \ref{caso-conrad}; see also 
\cite[Theorem 6.L]{glass})
\end{question}

\vspace{0.05cm}

Somehow related to the preceding question is the following well-known lemma, 
for which we provide a short proof based on the notion of right-recurrence.

\vspace{0.1cm}

\begin{lem} {\em If an orderable group $\Gamma$ admits only finitely many 
left-invariant total orders, then every element of $\mathcal{O}(\Gamma)$ 
is Conradian.}
\label{finito-conrad}
\end{lem}

\noindent{\bf Proof.} Since $\mathcal{O}(\Gamma)$ is finite, its points are periodic 
for the action of every element of $\Gamma$. This obviously implies that every order 
in $\mathcal{O}(\Gamma)$ is right-recurrent, hence Conradian. $\hfill\square$

\vspace{0.25cm}

\begin{rem} Using Tararin's theorem which describes all orderable groups admitting only finitely 
many orderings (see \S \ref{caso-conrad}), one can show that every ordering $\preceq$ on 
such a group satisfies the following: if $f$ is positive and $g$ is any group element, 
then \esp $fg^2 \succ g^2$. (This should be compared with Proposition \ref{lesla}.)
\end{rem}

\vspace{0.1cm}

The notion of right-recurrence for left-invariant orders was introduced by Morris-Witte 
in \cite{witte}, where he proves that every countable amenable orderable group is locally 
indicable. Actually, Morris-Witte proves that such a group always admits a right-recurrent 
ordering. His strategy shows how the dynamical properties of the action of an orderable 
group on its space of orderings can reveal some of its algebraic properties. 
His brilliant argument may be summarized as follows: 

\vspace{0.1cm}

\noindent -- since $\Gamma$ is amenable and $\mathcal{O}(\Gamma)$ is a compact metric space, 
the right action of $\Gamma$ on $\mathcal{O}(\Gamma)$ must preserve a probability measure 
(see for instance \cite{wagon});

\vspace{0.1cm}

\noindent -- if the right action of a countable orderable group $\Gamma$ on 
$\mathcal{O}(\Gamma)$ preserves a probability measure $\mu$, then the set 
of right-recurrent orderings has full $\mu$-measure, and in particular is 
non-empty (this follows by applying the Poincar\'e Recurrence Theorem).

\vspace{0.05cm}

\begin{question} If $\Gamma$ is countable amenable and orderable, is the set of 
right-recurrent orderings on $\Gamma$ dense inside the set of $\mathcal{C}$-orders~?
\label{densidad}
\end{question}

\vspace{0.05cm}

Since (countable) amenable groups do not contain free subgroups on two generators, 
it is natural to ask whether Morris-Witte's theorem is still true under the last 
(weaker) hypothesis. Partial evidence for an affirmative answer to this question 
is the result obtained by Linnell in \cite{linel}. The (apparently easier) question 
of the local indicability for orderable groups satisfying a non-trivial law (or identity) 
is still interesting. For instance, an affirmative answer for this case would allow to 
conclude that orderable groups satisfying an Engel type identity are locally nilpotent 
(see \cite[Theorem 6.G]{glass}). 

%%%%%%%%%%%%%%%%%%%%%%%%%%%%%%%%%%%%%%%%%%%%%%%%%%%%%%%%%%%%%%%%%%%%%%%%%%%%%%%%%%%%%%%%%%%%%%
%%%%%%%%%%%%%%%%%%%%%%%%%%%%%%%%%%%%%%%%%%%%%%%%%%%%%%%%%%%%%%%%%%%%%%%%%%%%%%%%%%%%%%%%%%%%%%

\section{Finitely many or a Cantor set of orders}

\subsection{The case of Conradian orders}
\label{caso-conrad}

\hspace{0.35cm} The approximation of Conradian orders is a problem of algebraic nature. In 
order to deal with it, we will use an elegant result by Tararin \cite{tararin} (see \cite{koko} 
for a detailed proof). For its statement, recall that a 
{\em rational series} for a group $\Gamma$ is a finite sequence of subgroups 
$$\{id\} = \Gamma^k \subset \Gamma^{k-1} \subset \ldots \subset \Gamma^0 = \Gamma$$
which is {\em subnormal} (that is, each $\Gamma^{i}$ is normal in $\Gamma^{i-1}$), 
and such that each quotient $\Gamma^{i-1} / \Gamma^i$ is torsion-free rank-one 
Abelian. Note that every group admitting a rational series is orderable.

\vspace{0.4cm}

\noindent{\bf Theorem [Tararin].} {\em If \esp \esp 
\esp $\Gamma$ is a group admitting a rational series}
$$\{id\} = \Gamma^k \subset \Gamma^{k-1} \subset \ldots \subset \Gamma^0 = \Gamma,$$ 
{\em then its space of orderings $\mathcal{O}(\Gamma)$ is finite if and only the subgroups $\Gamma^i$ 
are normal in $\Gamma$ and no quotient $\Gamma^{i-2} / \Gamma^i$ is bi-orderable. If this is the 
case, then $\Gamma$ admits a unique rational series, and for every left-invariant total order on 
$\Gamma$, the convex subgroups are precisely $\Gamma^0, \Gamma^1, \ldots, \Gamma^k$.}

\vspace{0.4cm}

Indeed, the number of orderings on a group satisfying the properties above equals $2^k$. 
Moreover, by choosing \esp $g_i \!\in\! \Gamma^i \setminus \Gamma^{i-1}$, \esp each of 
such orderings is uniquely determined by the sequence of {\em signs} of the elements 
\esp $g_i$. \esp Tararin's theorem will be fundamental for establishing the following 
proposition. (Note that there is no countability hypothesis for the group in the 
result below.)

\vspace{0.2cm}

\begin{prop} {\em If \esp $\Gamma$ is a Conrad orderable group having infinitely many 
left-invariant total orders, then all neighborhoods in \esp $\mathcal{O}(\Gamma)$ 
\esp of Conradian orders on $\Gamma$ do contain homeomorphic copies of the Cantor set.}
\label{dorilita}
\end{prop}

\vspace{0.15cm}

To prove this proposition we need to show that, if $\Gamma$ is an orderable group which admits 
a Conradian order having a neighborhood in $\mathcal{O}(\Gamma)$ which does not contain any
homeomorphic copy of the Cantor set, then $\Gamma$ admits a rational series as in the 
statement of Tararin's theorem. 
%To do this we will proceed in several steps. 

\vspace{0.15cm}

\begin{lem} {\em If a $\mathcal{C}$-order $\preceq$ on a group $\Gamma$ has a neighborhood 
in $\mathcal{O}(\Gamma)$ which does not contain any homeomorphic copy of the Cantor set, then 
$\Gamma$ admits a (finite) subnormal sequence formed by $\preceq$-convex subgroups so that the 
corresponding successive quotients are torsion-free Abelian.}
\label{hay}
\end{lem}

\noindent{\bf Proof.} Since the family of $\preceq$-convex subgroups is completely ordered 
by inclusion, referring to Remark \ref{convex-car} we just need to show that there exist 
only finitely many distinct subgroups of the form $\Gamma^g$. Let $\{f_1,\ldots,f_k\}$ 
be any finite family of elements of $\Gamma$. If there exist infinitely many distinct 
groups of the form $\Gamma^g$, then one may obtain an infinite ascending or 
descending sequence of these groups $\Gamma^{g_i}$ in such a way that \esp 
$f_m^{-1}f_n \notin \Gamma^{g_i} \setminus \Gamma_{g_i}$ for every \esp 
$m \neq n$ \esp in \esp $\{1,\ldots,k\}$ \esp and every \esp $i \!\in\! \mathbb{N}$. 
\esp Both cases being similar, we will consider only the former one. Following 
Zenkov \cite{zenkov}, for each \esp $i \!\in\! \mathbb{N}$ \esp and each \esp 
$\omega \!=\! (\ell_1,\ldots,\ell_i) \!\in\! \{0,1\}^i$ \esp let us inductively 
define the order 
\esp $\preceq_{\omega} = \preceq_{(\ell_1,\ldots,\ell_i)}$ \esp on $\Gamma^{g_i}$ by 
letting $\preceq_{\omega}$ be the extension of $\preceq_{(\ell_1,\ldots,\ell_{i-1})}$ 
by $\preceq$ (resp. by $\bar{\preceq}$) if $\ell_i = 0$ (resp. if $\ell_i = 1$). Passing 
to the limit, this allows to define a continuous embedding of the Cantor set \esp 
$\{0,1\}^{\mathbb{N}}$ \esp into the space of orderings of the subgroup \esp 
$\Gamma_* = \cup_{i \in \mathbb{N}} \esp \Gamma^{g_i}$, \esp which in its turn 
induces (just extending each resulting order on $\Gamma_*$ by $\preceq$) a continuous 
embedding of $\{0,1\}^{\mathbb{N}}$ into \esp $\mathcal{O} (\Gamma)$. \esp Moreover, 
since \esp $f_m^{-1}f_n \notin \Gamma^{g_i} \setminus \Gamma_{g_i}$ \esp for every \esp 
$m \neq n$ \esp in \esp $\{1,\ldots,k\}$ \esp and every \esp $i \!\in\! \mathbb{N}$, 
\esp the image of the latter embedding is contained in the neighborhood of $\preceq$ 
consisting of all orderings which do coincide with $\preceq$ on $\{f_1,\ldots,f_k\}$. 
Since this finite family of elements was arbitrary, this proves the lemma. $\hfill\square$

\vspace{0.45cm}

The lemma below concerns the rank of the quotients $\Gamma^{i-1} / \Gamma^i$.

\vspace{0.1cm}

\begin{lem} {\em Let $\preceq$ be a $\mathcal{C}$-order on a group $\Gamma$ having a 
neighborhood in $\mathcal{O}(\Gamma)$ which does not contain any homeomorphic copy 
of the Cantor set. If \esp \esp \esp $\{id\} = 
\Gamma^k \subset \Gamma^{k-1} \subset \ldots \subset \Gamma^0 = \Gamma$ \esp \esp 
\esp is a subnormal sequence of \esp $\Gamma$ formed by $\preceq$-convex subgroups   
so that each quotient $\Gamma^{i-1} / \Gamma^i$ is torsion-free Abelian, then the 
rank of each of these quotients equals one.}
\label{rango1}
\end{lem}

\noindent{\bf Proof.} For the proof we will use an elegant result by Sikora \cite{sikora} 
which establishes that $\mathcal{O}(\mathbb{Z}^n)$ has no isolated point (and it 
is therefore homeomorphic to the Cantor set) for every integer $n \geq 2$. 

Assume that some of the quotients \esp $\Gamma^{i-1} / \Gamma^i$ \esp 
has rank greater than or equal to 2. We will show that in this case every neighborhood 
of $\preceq$ contains a homeomorphic copy of the Cantor set. To do this, let \esp 
$\{f_1,\ldots,f_k\}$ \esp be any finite family of elements of $\Gamma$. Denoting by \esp 
$\pi: \Gamma^{i-1} \rightarrow \Gamma^{i-1} / \Gamma^i$ \esp the projection map, let $\Gamma_*$ 
be a subgroup of $\Gamma^{i-1}$ containing $\Gamma^i$, such that the rank of the quotient 
$\Gamma_* / \Gamma^{i}$ is finite and greater than or equal to 2, and such that each  
$f_{i}^{-1} f_j$ is contained in $\Gamma_* \cup (\Gamma \setminus \Gamma^{i-1})$. 
Let $\Gamma_{**}$ be the subgroup of $\Gamma^{i-1}$ containing $\Gamma^{i}$ and such that 
$\Gamma^{i-1}/\Gamma^i$ is the direct sum of $\Gamma_*/\Gamma^i$ and $\Gamma_{**}/\Gamma^i$.
By Sikora's result, the space of orderings of the quotient $\Gamma_* / \Gamma^{i}$ is 
homeomorphic to the Cantor set. For each $\preceq'$ in this space we may define 
an ordering $\preceq^*$ on $\Gamma$ by letting:

\vspace{0.1cm}

\noindent -- $\preceq^1$ be the order on $\Gamma^{i-1} / \Gamma^i$ defined by \esp 
$[g_1] + [h_1] \prec^1 [g_2] + [h_2]$ \esp if and only if either $[g_1] \prec' [g_2]$,  
or $[g_1] = [g_2]$, $[h_1] \neq [h_2]$, and $h_1 \prec h_2$. Here, for $i \! \in \! \{1,2\}$ 
the elements $g_i$ (resp. $h_i$) belong to $\Gamma_*$ (resp. $\Gamma_{**}$), and $[\cdot]$ 
stands for their class modulo $\Gamma^i$;

\vspace{0.1cm}

\noindent -- $\preceq^2$ be the order on $\Gamma^{i-1}$ for which an element $g$ is positive if and 
only if either $g \in \Gamma^{i}$ and $g \succ id$, or $g \notin \Gamma^{i}$ and $id \prec^1 [g]$;

\vspace{0.1cm}

\noindent -- $\preceq^*$ be the extension of $\preceq^2$ by $\preceq$.

\vspace{0.1cm}

The map \esp $\preceq' \esp \mapsto \esp \preceq^*$ \esp is continuous and injective. 
Therefore, the intersection of its image with the subset of \esp $\mathcal{O}(\Gamma)$ 
\esp consisting of all orderings which do coincide with $\preceq$ on 
$\{f_1,\ldots,f_k\}$ corresponds to a homeomorphic copy of the Cantor set inside the 
corresponding neighborhood of $\preceq$ in $\mathcal{O}(\Gamma)$. Once again, since 
this finite family of elements was arbitrary, this proves the lemma. $\hfill\square$

\vspace{0.5cm}

The next lemma is essentially due to Linnell \cite{linnell} (see also \cite{zenkov}). 

\vspace{0.1cm}

\begin{lem} {\em Let $\Gamma$ be a group and $\Gamma^1$ a normal subgroup such that $\Gamma^1$ and 
$\Gamma / \Gamma^1$ are torsion-free Abelian of rank one. Let $\preceq$ be a Conradian order on 
$\Gamma$ respect to which $\Gamma^1$ is a convex subgroup. If $\Gamma$ is bi-orderable, then 
every neighborhood of $\preceq$ in $\mathcal{O}(\Gamma)$ contains a homeomorphic copy of 
the Cantor set.}
\label{dos}
\end{lem}

\noindent{\bf Proof.} Let us consider the action by conjugacy \esp $\alpha: \Gamma / \Gamma^1 \rightarrow 
\mathrm{Aut}(\Gamma^1)$, \esp namely $\alpha (g\Gamma^1) (h) = ghg^{-1}$, where $g \!\in\! \Gamma$ and 
$h \!\in\! \Gamma^1$. If $\alpha$ is trivial then $\Gamma$ is Abelian and its rank is necessarily 
greater than or equal to 2. However, this together with the hypothesis is in contradiction 
with Sikora's theorem. If \esp $\{id\} \neq Ker(\alpha) \neq \Gamma / \Gamma^1$ \esp then 
$(\Gamma/ \Gamma^1)/Ker(\alpha)$ is a non-trivial torsion group, and since the only non-trivial 
finite order automorphism of $\Gamma^1$ is the inversion, there must exist $g \!\in\! \Gamma$ such 
that $ghg^{-1} = h^{-1}$ for every $h \!\in\! \Gamma$. This obviously implies that $\Gamma$ is not 
bi-orderable. Therefore, $Ker(\alpha) = \{id\}$ and $\Gamma / \Gamma^1 \sim (\mathbb{Z},+)$. Viewing 
$\Gamma^1$ as a subgroup of $\mathbb{Q}$, the action of $(\mathbb{Z},+)$ is generated by the 
multiplication by a non zero rational number $q$. If $q$ is negative then $\Gamma$ is still non 
bi-orderable. It just remains the case where $q$ is positive. Note that in this case $\Gamma$ 
embeds in the affine group; more precisely, $\Gamma$ can be identified with the group whose 
elements are of the form 
$$(k,a) \sim 
\left(
\begin{array}
{cc}
q^k & a  \\
0   & 1  \\
\end{array}
\right),$$
where $a \in \Gamma^1$ and $k \!\in\! (\mathbb{Z},+)$. Let $(k_1,a_1),\ldots,(k_n,a_n)$ 
be an arbitrary family of positive elements of $\Gamma$ indexed in such a way that 
$k_1 = k_2 = \ldots = k_r = 0$ and $k_{r+1} \neq 0, \ldots, k_n \neq 0$ for some 
$r \!\in\! \{1,\ldots,n\}$. Four cases are possible:

\vspace{0.1cm}

\noindent (i) \esp \esp \esp $a_1 > 0, \ldots, a_r > 0$ \esp \esp and 
\esp \esp $k_{r+1} > 0, \ldots, k_n > 0$, 

\vspace{0.1cm}

\noindent (ii) \esp \esp \esp $a_1 < 0, \ldots, a_r < 0$ \esp \esp and 
\esp \esp $k_{r+1} > 0, \ldots, k_n > 0$, 

\vspace{0.1cm}

\noindent (iii) \esp \esp \esp $a_1 > 0, \ldots, a_r > 0$ \esp \esp and 
\esp \esp $k_{r+1} < 0, \ldots, k_n < 0$, 

\vspace{0.1cm}

\noindent (iv) \esp \esp \esp $a_1 < 0, \ldots, a_r < 0$ \esp \esp and 
\esp \esp $k_{r+1} < 0, \ldots, k_n < 0$. 

\vspace{0.1cm}

\noindent As in Example \ref{ex-smirnov}, for each irrational number \esp $\varepsilon$ \esp 
let us consider the ordering $\preceq_{\varepsilon}$ on $\Gamma$ whose positive cone is 
$$P_{\preceq_{\varepsilon}} = \{(k,a) \! : \esp \esp \esp q^k + \varepsilon a > 1\}.$$
Note that if $\varepsilon_1 \neq \varepsilon_2$ then $\preceq_{\varepsilon_1}$ is 
different from $\preceq_{\varepsilon_2}$. (Remark also that no order $\preceq_{\varepsilon}$ 
is Conradian.) Now in case (i), for $\varepsilon$ positive and very small the order 
$\preceq_{\varepsilon}$ is different from $\preceq$ but still makes all the elements 
$(k_i,a_i)$ positive. The same is true in case (ii) for $\varepsilon$ negative and near 
zero. In case (iii) this still holds for the order $\bar{\preceq_{\varepsilon}}$ 
when $\varepsilon$ is negative and near zero. Finally, in case (iv) one needs to consider 
again the order $\bar{\preceq_{\varepsilon}}$ but for $\varepsilon$ positive and small. 
Now letting $\varepsilon$ vary over a Cantor set formed by irrational 
numbers\footnote{Take for instance the set of numbers of the form 
$\sum_{i \geq 1} \frac{i_k}{4^k}$, where $i_k \!\in\! \{0,1\}$, and 
translate it by $\sum_{j \geq 1} \frac{2}{4^{j^2}}$.} very near to $0$ 
(and which are positive or negative according to the case), this shows that the 
neighborhood of $\preceq$ consisting of the orderings on $\Gamma$ which make all of the 
elements $(k_i,a_i)$ positive contains a homeomorphic copy of the Cantor set. Since the 
finite family of elements $(k_i,a_i)$ which are positive for $\preceq$ was arbitrary, 
this proves the lemma. $\hfill\square$

\vspace{0.55cm}

We may now pass to the proof of Proposition \ref{dorilita}. By 
Lemmas \ref{hay} and \ref{rango1}, every countable group $\bar{\Gamma}$ admitting 
a $\mathcal{C}$-order $\preceq'$ having a neighborhood in $\mathcal{O}(\bar{\Gamma})$ 
which does not contain any homeomorphic copy of the Cantor set admits a rational series 
$$\{id\} = \bar{\Gamma}^k \subset \bar{\Gamma}^{k-1} \subset \ldots 
                    \subset \bar{\Gamma}^1 \subset \bar{\Gamma}^0 = \bar{\Gamma}$$
formed by $\preceq'$-convex subgroups. Assume by contradiction that the family $\mathcal{F}$ 
of these groups $\bar{\Gamma}$ having an infinite space of orderings is non-empty. For each 
$\bar{\Gamma}$ in $\mathcal{F}$ let $k(\bar{\Gamma}) \!\in\! \mathbb{N}$ be the minimum 
possible length for a rational series formed by $\preceq'$-convex subgroups with respect 
to some $\mathcal{C}$-order $\preceq'$ having a neighborhood in $\mathcal{O}(\bar{\Gamma})$ 
which does not contain any homeomorphic copy of the Cantor set. Let $k$ the minimum 
of $k(\bar{\Gamma})$ for $\bar{\Gamma}$ ranging over all groups in $\mathcal{F}$, 
and let $\Gamma$ and $\preceq$ be respectively a countable group in $\mathcal{F}$ 
and a $\mathcal{C}$-order on it realizing this value $k$. Clearly, one has $k \neq 0$ and 
$k \neq 1$. Moreover, Lemma \ref{dos} together with Tararin's theorem implies that $k \neq 2$. 

To get a contradiction in the other cases, we fist claim that all the corresponding subgroups 
$\Gamma^i$ are normal in $\Gamma$. Indeed, the restriction of $\preceq$ to $\Gamma^1$ is 
Conradian, and it clearly has a neighborhood in $\mathcal{O}(\Gamma^1)$ which does not 
contain any homeomorphic image of the Cantor set. Since 
$$\{id\} = \Gamma^k \subset \Gamma^{k-1} \subset \ldots \subset \Gamma^1$$
is a rational series of length $k-1$ formed by $\preceq$-convex subgroups of $\Gamma^1$, the 
minimality of the index $k$ implies that $\mathcal{O}(\Gamma^1)$ is finite. By Tararin's theorem, 
the rational series for $\Gamma^1$ is unique. Therefore, since $\Gamma^1$ is already normal in 
$\Gamma$, for every $g \!\in\! \Gamma$ the rational series for $\Gamma^1$ given by   
$$\{id\} = g \Gamma^k g^{-1} \subset g \Gamma^{k-1} g^{-1} \subset \ldots 
\subset g \Gamma^1 g^{-1} = \Gamma^1$$
must coincide with the original one. Since the 
element $g \!\in\! \Gamma$ was arbitrary, this shows 
that all the subgroups $\Gamma^i$ are normal in $\Gamma$. 

We now claim that no quotient $\Gamma^{i-2} / \Gamma^i$ is bi-orderable. Indeed, for the normal sequence 
$$\{id\} = \Gamma^i / \Gamma^i \subset \Gamma^{i-1} / \Gamma^i \subset \Gamma^{i-2} / \Gamma^i$$
the groups \esp \esp $\Gamma^{i-1} / \Gamma^i$ \esp \esp and 
$$(\Gamma^{i-2} / \Gamma^i) / (\Gamma^{i-1} / \Gamma^i) \esp \esp 
\sim \esp \esp \Gamma^{i-2} / \Gamma^{i-1}$$ 
are torsion-free rank-one Abelian. Moreover, $\preceq$ induces a Conradian order $\preceq'$ 
on the quotient $\Gamma^{i-2} / \Gamma^i$ respect to which $\Gamma^{i-1} / \Gamma^i$ is convex. 
Since $\preceq$ has a neighborhood in $\mathcal{O}(\Gamma)$ which does not contain any 
homeomorphic copy of the Cantor set, an extension type argument shows that a similar property 
holds for $\preceq'$ inside $\mathcal{O} (\Gamma^{i-2} / \Gamma^i)$. The fact that 
$\Gamma^{i-2} / \Gamma^i$ is not bi-orderable then follows from Lemma \ref{dos}.

We already know that each $\Gamma^i$ is normal in $\Gamma$ and no quotient $\Gamma^{i-2}/\Gamma^i$ 
is bi-orderable. As another application of Tararin's theorem we obtain that the space of orders 
$\mathcal{O}(\Gamma)$ is finite, thus finishing the proof of Proposition \ref{dorilita}. 

\vspace{0.5cm}

\noindent{\bf Proof of Theorem B.} An easy consequence of Tararin's theorem is that a 
non-trivial torsion-free nilpotent group which admit only finitely many orderings is 
rank-one Abelian. By the comments just before Figure 1, every ordering on an orderable group 
without free semigroups on two generators (and therefore, every ordering on a torsion-free 
nilpotent group) is Conradian. It follows from Proposition \ref{dorilita} that if 
$\Gamma$ is a non-trivial torsion-free nilpotent group which is not rank-one Abelian, 
then $\mathcal{O}(\Gamma)$ has no isolated point. As a consequence, if $\Gamma$ 
is countable, then $\mathcal{O}(\Gamma)$ is a totally disconnected compact metric 
space without isolated points, and therefore homeomorphic to the Cantor set 
(see \cite[Theorem 2-80]{HY})). This proves the first claim of Theorem B. 
The second claim of the theorem follows directly from the first 
one and Proposition \ref{paja}. $\hfill\square$

\vspace{0.2cm}

\begin{rem} The main property used in the proof above is that every ordering 
on a torsion-free nilpotent group is Conradian. This holds more generally for  
orderable groups without free semigroups on two generators. Actually, the 
conclusion of Theorem B applies to all these groups, provided they are 
countable and orderable. A relevant example, namely Grigorchuk-Maki's 
group of intermediate growth, was extensively studied in \cite{growth}. 
\end{rem}

%%%%%%%%%%%%%%%%%%%%%%%%%%%%%%%%%%%%%%%%%%%%%%%%%%%%%%%%%%%%%%%%%%%%%%%%%%%%%%%%%%%%%%%%%%%%%%%%%%%%%%%%%
%%%%%%%%%%%%%%%%%%%%%%%%%%%%%%%%%%%%%%%%%%%%%%%%%%%%%%%%%%%%%%%%%%%%%%%%%%%%%%%%%%%%%%%%%%%%%%%%%%%%%%%%%

\subsection{The case of orders with trivial Conradian soul}
\label{caso-trivial}

\hspace{0.35cm} In the ``pure non Conradian case'' (that is, when the Conradian soul is 
trivial), our method for approximating a given ordering on a (countable infinite) group will 
consist in taking conjugates of it. More precisely, 
given a countable orderable group $\Gamma$ and an element $\preceq$ of 
$\mathcal{O}(\Gamma)$, we will denote by $\mathrm{orb} (\preceq)$ the orbit of $\preceq$ 
by the right action of $\Gamma$. We begin by noting that, if $\preceq$ is non isolated in 
$\mathrm{orb} (\preceq)$, then the closure $\overline{\mathrm{orb}(\preceq)}$ is a $\Gamma$-invariant 
closed subset of $\mathcal{O}(\Gamma)$ without isolated points, and 
therefore homeomorphic to the Cantor set (because $\mathcal{O}(\Gamma)$ is metrizable and totally 
disconnected). To show that a particular order is non isolated inside its orbit (that is, it 
may be approximated by its conjugates), the following elementary lemma will be very useful.

\vspace{0.35cm}

\begin{lem} {\em Let $\precede$ be an ordering on a countable group $\Gamma$. 
Assume that the following property holds for the dynamical realization of 
\esp $\precede$ \esp associated to a numbering $(g_i)_{i \geq 0}$ of \esp 
$\Gamma$ such that $g_0 \!=\! id$: for every $\varepsilon > 0$ 
there exists \esp $g \succ id$ \esp and \esp $x \!\in \![-\varepsilon,\varepsilon]$ \esp such 
that \esp $g(x) < x$. \esp Then $\precede$ is a non isolated point of \esp $\mathrm{orb}(\preceq)$.}
\label{idea-tonta}
\end{lem}

\noindent{\bf Proof.} Fix a complete exhaustion $\mathcal{G}_0 \subset \mathcal{G}_1 \subset \ldots$ 
of $\Gamma$ by symmetric finite sets. We need to show that for all fixed $n \!\in\! \mathbb{N}$ 
there exists $\precede_n$ in $\mathrm{orb}(\preceq)$ different from $\precede$ such that an element 
$g \!\in\! \mathcal{G}_n$ satisfies \esp $g \succ_n id$ \esp if and only if \esp $g \succ id$. 
\esp Now recall that, for all $h \in \Gamma$, the value of \esp $h (0) = h(t(id)) = t(h)$ \esp is 
positive (resp. negative) if and only if \esp $h \succ id$ \esp (resp. \esp $h \prec id$). \esp For 
each $h \succ id$ denote \esp $\varepsilon (h) = \inf \{ |x|: h(x) \leq x \}$. \esp 
(We remark that $\varepsilon (h)$ is strictly positive, perhaps equal to infinite.)  Now let  
$$\varepsilon_n = \min \{ \varepsilon(g)\!: \esp g \succ id, \esp g \in \mathcal{G}_n \}.$$ 
By the ``transversality'' hypothesis, there exists an element $g_n \!\succ\! id$ in $\Gamma$ 
such that \esp $g_n (x_n) < x_n$ \esp for some \esp $x_n \!\! \in ]-\varepsilon_n,\varepsilon_n[$. 
\esp Moreover, according to the comments after Proposition \ref{thorden}, such a point $x_n$ may 
be taken equal to $t(h_n^{-1})$ for some element $h_n \in \Gamma$. Now consider the order relation 
\esp $\precede_n = h_n (\precede)$, \esp that is,  \esp $g \succ_n id$ \esp if and only if \esp 
$g(x_n) > x_n$. \esp The equivalence between the conditions \esp $g \succ id$ \esp and 
\esp $g \succ_n id$ \esp holds for every \esp $g \in \mathcal{G}_n$ \esp by the definition of 
$\varepsilon_n$. On the other hand, one has \esp $g_n \succ id$ \esp and \esp $g_n \prec_n id$, 
\esp thus showing that \esp $\precede$ \esp and \esp $\precede_n$ \esp are different. $\hfill\square$

\vspace{0.5cm}

The transversality hypothesis does not hold for all dynamical realizations. Indeed, according 
to \S \ref{invariant}, if the order $\precede$ is bi-invariant then (for the associated 
dynamical realization) the graph of no element crosses the diagonal. It seems also difficult 
to apply directly the previous argument for general $\mathcal{C}$-orders. However, according 
to \S \ref{conrad-soul}, the transversality condition clearly holds when the Conradian 
soul of $\preceq$ is trivial. As a consequence, we obtain the following proposition.

\vspace{0.15cm}

\begin{prop} {\em If an ordering $\preceq$ on a non-trivial countable group 
$\Gamma$ has trivial Conradian soul, then $\preceq$ is an accumulation point of 
its set of conjugates. In particular, the closure of the orbit of $\preceq$ 
under the right action of \esp\esp $\Gamma$ is homeomorphic to the Cantor set.}
\label{unilita}
\end{prop}

\vspace{0.05cm}

\begin{question} Does there exist a pure algebraic characterization of the elements of 
$\mathcal{O}(\Gamma)$ which are not accumulation points of their orbits by the action of 
$\Gamma$ (equivalently, of the orderings which are non approximable by their conjugates)~?
\end{question}

%%%%%%%%%%%%%%%%%%%%%%%%%%%%%%%%%%%%%%%%%%%%%%%%%%%%%%%%%%%%%%%%%%%%%%%%%%%%%%%%%%%%%%%%%%%%%%%%%%%%
%%%%%%%%%%%%%%%%%%%%%%%%%%%%%%%%%%%%%%%%%%%%%%%%%%%%%%%%%%%%%%%%%%%%%%%%%%%%%%%%%%%%%%%%%%%%%%%%%%%%

\subsection{The general case}
\label{caso-general}

\hspace{0.35cm} For Conrad orderable groups, Theorem C follows immediately from Proposition 
\ref{dorilita}. If $\Gamma$ has an ordering $\preceq$ having a Conradian soul $\con$ admitting 
infinitely many orders, then $\mathcal{O}(\con)$ contains a homeomorphic copy of the Cantor set. 
Therefore, extending by $\preceq$ all the orderings on $\con$ to the whole group $\Gamma$, 
we obtain a homeomorphic copy of the Cantor set inside $\mathcal{O}(\Gamma)$. 

Since for the 
case of trivial Conradian soul Proposition \ref{unilita} applies, it just remains the case 
of a non Conradian ordering $\preceq$ whose Conradian soul is non-trivial but admits 
only finitely many orderings. Let \esp $\preceq_1,\ldots,\preceq_{2^k}$ \esp be all 
of the elements of $\mathcal{O} (\Gamma_{\preceq}^{\esp \esp c})$. \esp 
For $j \!\in\! \{1,\ldots,2^k\}$ denote by $\preceq^j$ the extension of $\preceq_j$  
by $\preceq$. Note that, by Lemmas \ref{esta-con} and \ref{finito-conrad}, the subgroup $\con$ 
coincides with the Conradian soul of $\Gamma$ with respect to all of the orderings $\preceq^j$. 
To finish the proof of Theorem C, it suffices to show the following.

\vspace{0.35cm}

\begin{prop} {\em With the notations above, at least one of the 
orderings $\preceq^j$ is an accumulation point of its orbit.}
\label{sandra}
\end{prop}

\vspace{0.2cm}

For the proof of this proposition, fix a numbering $(g_i)_{i \geq 0}$ of the elements 
of \esp $\Gamma$ such that \esp $g_0 \!=\! id$, \esp and denote by $\alpha \!<\! 0$ 
and $\beta \!>\! 0$ the constants appearing in the corresponding dynamical realization 
of $\preceq$ associated to the Conradian soul $\con$ ({\em c.f.} Proposition \ref{reca}). 

\vspace{0.37cm}

\noindent{\underbar{Claim 1.}} For every $\varepsilon \!>\! 0$ there 
exist \esp $f_{\varepsilon},g_{\varepsilon}$ \esp in $\Gamma$ and  
\esp $a_{\varepsilon} , b_{\varepsilon}$ \esp in $]\beta, \beta + \varepsilon[$ 
\esp such that $f_{\varepsilon},g_{\varepsilon}$ are in transversal 
position on $[a_{\varepsilon},b_{\varepsilon}]$.

\vspace{0.35cm}

Indeed, by the definition of $\beta$, there exist elements $f,g$ in $\Gamma$ which are in transversal 
position on some interval $[a,b]$ such that \esp $\beta \leq a < \beta + \varepsilon$. Changing $g$ by 
$f^n g f^{-n}$ for $n \in \mathbb{N}$ large enough, we may suppose that \esp $b < \beta + \varepsilon$. 
\esp Similarly, changing $f$ by $gfg^{-1}$ if necessary, we may also assume that \esp $a > \beta$. 

%\esp Denote by $a_m$ (resp. $b_m$) the first (resp. the last) fixed point of $gf^m$ in $]a,b[$. If \esp 
%$m' >\!\!> m$ \esp are large enough then \esp $b_{m'} < a_m$, \esp and therefore the claim holds for the 
%constants \esp $a_{\varepsilon} =  b_{m'}$ \esp and \esp $b_{\varepsilon} = a_m$, \esp and for the elements 
%\esp $f_{\varepsilon} = g f^{m'}$ \esp and \esp $g_{\varepsilon} = g f^m$ (see Figure 2). 

\vspace{0.6cm}

%%%%%%%%%%%%%%%%%%%%%%%%%%%%%%%%%%%%%%%%%%%%%%%%%%%%%%%%%%%%%%%%%%%%%%%%%%%%%%%%%%%%%%%%%%
%%%%%%%%%%%%%%%%%%%%%%%%%%%%%%%%%%%%%%%%%%%%%%%%%%%%%%%%%%%%%%%%%%%%%%%%%%%%%%%%%%%%%%%%%%

\beginpicture

\setcoordinatesystem units <1cm,1cm>

\putrule from 0 0 to 0 8
\putrule from 0 0 to 8 0
\putrule from 8 0 to 8 8 
\putrule from 0 8 to 8 8

\put{$a_{\varepsilon_n}$} at 0 -0.4
\put{$b_{\varepsilon_n}$} at 8.1 -0.4 
\put{$a_{m_{i+1}}$} at 2.6 -0.4
\put{$b_{m_{i+1}}$} at 3.6 -0.4
\put{$a_{m_i}$} at 6.0 -0.4
\put{$b_{m_i}$} at 6.7 -0.4 
\put{$\bar{h}_{m_i}$} at 7.3 6.45 
\put{$\bar{h}_{m_{i+1}}$} at 7.4 3.5 

\plot 
0 2.4 
2.7 2.7 /

\plot
0 2.4 
6 6 /

\plot 
3.3 3.3 
8 4 /

\plot 
6.6 6.6 
8 7.2  /

\put{$\bullet$} at 4.13 4.13 
\put{$t(h_i^{-1})$} at -0.7 4.13 

%%%%%%%%%%%%%%%%%%%%%%%%%%%%%%%%%%%%%%%%%%%%%%%%%%%%%%%%%%%%%%%%%%%%%%%%%%%
%%%%%%%%%%%%%%%%%%%%%%%%%%%%%%%%%%%%%%%%%%%%%%%%%%%%%%%%%%%%%%%%%%%%%%%%%%%
%%%%%%%%%%%%%%%%%%%%%%%%%%%%%%%%%%%%%%%%%%%%%%%%%%%%%%%%%%%%%%%%%%%%%%%%%%%
%%%%%%%%%%%%%%%%%%%%%%%%%%%%%%%%%%%%%%%%%%%%%%%%%%%%%%%%%%%%%%%%%%%%%%%%%%%%

\setquadratic 

\plot 
2.7 2.7 2.92 2.8 3 3 /

\plot 
3 3 3.1 3.2 3.3 3.3 /

\plot 
6 6 6.2 6.12 6.3 6.3 /

\plot 
6.3 6.3 6.4 6.47 6.6 6.6 /

\setlinear

\setdots

\plot 
0 0 
8 8 /

\putrule from 3.3 3.75 to 6 3.75   
\putrule from 3.3 4.38 to 6 4.38  

\putrule from 3.3 3.3 to 3.3 6 
\putrule from 3.3 3.3 to 6 3.3 
\putrule from 3.3 6 to 6 6 
\putrule from 2.7 0 to 2.7 2.7 
\putrule from 3.3 0 to 3.3 3.3 
\putrule from 6 0 to 6 6 
\putrule from 6.6 0 to 6.6 6.6 

\putrule from 4.13 0 to 4.13 4.13
\putrule from 0 4.13 to 4.13 4.13 

\put{Figure 2} at 4 -1

\put{} at -4.2 0

\endpicture

%%%%%%%%%%%%%%%%%%%%%%%%%%%%%%%%%%%%%%%%%%%%%%%%%%%%%%%%%%%%%%%%%%%%%%%%%%%%%%%%
%%%%%%%%%%%%%%%%%%%%%%%%%%%%%%%%%%%%%%%%%%%%%%%%%%%%%%%%%%%%%%%%%%%%%%%%%%%%%%%%

\vspace{0.4cm}

For $g \in \Gamma \setminus \Gamma_{\preceq}^{\esp \esp c}$ \esp such 
that \esp $g \succ id$, \esp let \esp $\varepsilon (g) \!> 0$ \esp be the 
positive number defined by \esp $\varepsilon (g) = g(0) - \beta$. Let 
$\mathcal{G}_0 \subset \mathcal{G}_1 \subset \ldots$ be a complete 
exhaustion of $\Gamma$ by finite sets. Given $n \in \mathbb{N}$ 
let $\varepsilon_n$ be the (positive) number defined by
\begin{equation}
\varepsilon_n = \min \big\{ \varepsilon (g): \esp \esp g \succ id, 
\esp \esp g \in \mathcal{G}_n \setminus \con \big\}.
\label{def-em}
\end{equation}
Put $\bar{f} = f_{\varepsilon_n}$ and $\bar{g} = g_{\varepsilon_n}$. For \esp 
$m \!\geq\! 1$ \esp let $a_m$ (resp. $b_m$) be the first (resp. the last) 
fixed point of the element \esp $ \bar{h}_m = \bar{g} \bar{f}^m$ \esp in 
\esp $]a_{\varepsilon_n}, b_{\varepsilon_n}[$. \esp It is not difficult to 
check that, choosing an appropriate subsequence $(m_i)$, we may ensure 
that for each $i \in \mathbb{N}$ the following hold (see Figure 2):

\vspace{0.1cm}

\noindent -- \esp $a_{m_i} \! > \! b_{m_{i+1}}$, 

\vspace{0.1cm}

\noindent -- \esp $\bar{h}_{m_{i+1}} (a_{m_{i}}) < \bar{h}_{m_{i}} (b_{m_{i+1}})$, 

\vspace{0.1cm}

\noindent -- there exists $h_{i} \in \Gamma$ such that $t(h_i^{-1})$ belongs to 
the interval \esp $]\bar{h}_{m_{i+1}} (a_{m_{i}}),\bar{h}_{m_{i}} (b_{m_{i+1}})[$.

\vspace{0.35cm}

\noindent{\underbar{Claim 2.}} For each $i \! \in \! \mathbb{N}$ and each 
$j \!\in\! \{1,\ldots,2^k\}$, an element in $\mathcal{G}_n \setminus \con$ 
belongs to the positive cone of $(\preceq^j)_{h_i}$ if and only if it 
belongs to the positive cone of $\preceq$.

\vspace{0.3cm}

Indeed, for any element \esp $h \in \mathcal{G}_n \setminus \con$ 
\esp which is positive with respect to $\preceq$ one has 
$$t(h h_i^{-1}) = h (t(h_i^{-1})) > h(0) \geq \beta + \varepsilon_n > a_{m_{i-1}} > t(h_i^{-1}).$$
This implies that \esp $h h_i^{-1} \succ h_i^{-1}$, \esp and therefore $h_i h h_i^{-1} \succ id$. 
If we show that the element $h_i h h_i^{-1}$ is not contained in $\con$, then this would give 
\esp \esp $h_i h h_i^{-1} \succ^j id$, \esp \esp that is, $h$ is positive with respect to 
$(\succ^j)_{h_i}$. Now, if $h_i h h_i^{-1}$ was equal to 
some element $\bar{h} \in \con$, then the interval 
$$h_i ([t(h_i^{-1}),t(hh_i^{-1})]) = [0,t(\bar{h})] \subset \esp ]\alpha,\beta[$$ 
would contain in its interior the interval $[h_i(b_{m_{i}}), h_i(a_{m_{i-1}})]$ over 
which the elements $h_i \bar{h}_{m_{i}} h_i^{-1}$ and $h_i \bar{h}_{m_{i-1}} h_i^{-1}$ 
are crossed. However, this contradicts the definition of the interval $]\alpha,\beta[$. 

If \esp $h \in \mathcal{G}_n \setminus \con$ \esp is negative with respect to $\preceq$, the 
above argument shows that $h^{-1}$ is positive with respect to $(\succ^{j})_{h_i}$, and therefore 
$h$ is negative with respect to this ordering as well. This finishes the proof of Claim 2.

\vspace{0.5cm}

%%%%%%%%%%%%%%%%%%%%%%%%%%%%%%%%%%%%%%%%%%%%%%%%%%%%%%%%%%%%%%%%%%%%%%%%%%%%%

\beginpicture

\setcoordinatesystem units <1cm,1cm>

\setquadratic 

\plot 
0 2 1.6 2.55 2.7 2.7 /

\plot 
2.7 2.7 2.92 2.82 3 3 /

\plot 
3 3 3.1 3.2 3.3 3.3 /

\plot 
5 5 5.2 5.12 5.3 5.3 /

\plot 
5.3 5.3 5.4 5.47 5.6 5.6 /

\plot 
0 2 2 6.16 8 8 /

\plot 0 0 4 0.98 8 1.4 /

%%%%%%%%%%%%%%%%%%%%%%%%%%%%

\setlinear

\putrule from 0 0 to 0 8
\putrule from 0 0 to 8 0
\putrule from 8 0 to 8 8 
\putrule from 0 8 to 8 8

\put{$b_{m_{i+1}}$} at 0 -0.4
\put{$a_{m_i}$} at 8 -0.4 
\put{$\bullet$} at 1.7 1.7
\put{$a$} at 3.3 -0.4
\put{$b$} at 5.0 -0.4
\put{$t(h_i^{-1})$} at 1.7 -0.4 
\put{$\bar{h}_{m_{i+1}}$} at 7 1 
\put{$\bar{h}_{m_i} \bar{h}_{m_{i+1}}^{n}$} at 7 5.4 
\put{$\bar{h}_{m_i}$} at 2.6 7  
\put{$\bar{h}_{m_i} \bar{h}_{m_{i+1}}^{n'}$} at 7 3.859 
\put{$\bar{h}_{m_i}(t(h_i^{-1}))$} at -1 5.9 

\plot
0 2 
5 5 /

\plot 
3.3 3.3 
8 4.7 /

\plot 
5.6 5.6 
8 6.2  /

%%%%%%%%%%%%%%%%%%%%%%%%%%%%%%%%%%%%%%%%%%%%%%%%%%%%%%%%%%%%%%%%%%%%%%%%%%%
%%%%%%%%%%%%%%%%%%%%%%%%%%%%%%%%%%%%%%%%%%%%%%%%%%%%%%%%%%%%%%%%%%%%%%%%%%%
%%%%%%%%%%%%%%%%%%%%%%%%%%%%%%%%%%%%%%%%%%%%%%%%%%%%%%%%%%%%%%%%%%%%%%%%%%%
%%%%%%%%%%%%%%%%%%%%%%%%%%%%%%%%%%%%%%%%%%%%%%%%%%%%%%%%%%%%%%%%%%%%%%%%%%%%

\setdots

\plot 
0 0 
8 8 /

\putrule from 5.9 0 to 5.9 5.9    
\putrule from 3.3 3.3 to 3.3 5 
\putrule from 3.3 3.3 to 5 3.3 
\putrule from 3.3 5 to 5 5 
\putrule from 3.3 0 to 3.3 3.3 
\putrule from 5 0 to 5 5 
\putrule from 1.7 0 to 1.7 5.9 
\putrule from 0 5.9 to 5.9 5.9 

\putrule from 0 1.4 to 8 1.4
\putrule from 0 2 to 8 2

\put{Figure 3} at 4.1 -1

\put{} at -4 0

\endpicture

%%%%%%%%%%%%%%%%%%%%%%%%%%%%%%%%%%%%%%%%%%%%%%%%%%%%%%%%%%%%%%%%%%%%%%%%%%%
%%%%%%%%%%%%%%%%%%%%%%%%%%%%%%%%%%%%%%%%%%%%%%%%%%%%%%%%%%%%%%%%%%%%%%%%%%%

\vspace{0.5cm}

\noindent{\underbar{Claim 3.}} For each fixed \esp $j \! \in \! \{1,\ldots,2^k\}$ \esp 
the orders $(\preceq^j)_{h_i}$ are two-by-two distinct (for \esp $i \!\in\! \mathbb{N}$).

\vspace{0.3cm}

It easily follows from the construction that the inequality \esp 
$\bar{h}_{m_{\ell}} (t(h_i^{-1})) > t(h_i^{-1})$ \esp holds if and only if \esp 
$\ell \!\leq\! i$. \esp If this is the case, then \esp $\bar{h}_{m_{\ell}} (t (h_i^{-1})) 
> \bar{h}_{m_{i}} (b_{m_{i+1}})$. \esp Therefore, for \esp $n' >\!\!> n$ \esp 
large enough, the elements $f_{n'} = \bar{h}_{m_{i}} \bar{h}_{m_{i + 1}}^{n'}$ 
and $f_n = \bar{h}_{m_{i}} \bar{h}_{m_{i + 1}}^{n}$ are 
in transversal position on some closed interval $[a,b]$ 
contained in $]t(h_i^{-1}),\bar{h}_{m_{i}} (t(h_i^{-1}))[$ 
(see Figure 3). We claim that this 
implies that the element \esp $h_i \bar{h}_{m_{\ell}} h_i^{-1}$ \esp does not belong to $\con$ 
for all \esp $\ell \leq i$. \esp Indeed, if $h_i \bar{h}_{m_{\ell}} h_i^{-1}$ was equal to 
some element \esp $\bar{h} \in \con$ \esp then, since \esp $a > t(h_i^{-1})$ \esp and \esp 
$b < t(\bar{h}_{m_i} h_i^{-1}) \leq t(\bar{h}_{m_{\ell}} h_i^{-1})$, \esp the interval 
$$[0,t(\bar{h})] = [0,t(h_i \bar{h}_{m_{\ell}} h_i^{-1})] = 
h_i ([t(h_i^{-1}),t(\bar{h}_{m_{\ell}} h_i^{-1})])$$ 
would be contained in $[0,\beta]$ and would contain in its interior the interval 
\esp $[h_i(a),h_i(b)]$. \esp However, on the last interval the elements 
$h_i f_{n'} h_i^{-1}$ and $h_i f_n h_i^{-1}$ are in transversal position, 
and this contradicts the definition of the interval $]\alpha,\beta[$.

Now since \esp $h_i \bar{h}_{m_{\ell}} h_i^{-1} \succ id$ \esp for all \esp $\ell \leq i$, 
one also has \esp $h_i \bar{h}_{m_{\ell}} h_i^{-1} \succ^j id$ \esp for all \esp 
$j \! \in \! \{1,\ldots,2^k\}$. \esp In other words, the element $\bar{h}_{m_{\ell}}$ is 
positive with respect to \esp $(\succ^j)_{h_i}$ \esp for every $\ell \leq i$. In an analogous 
way, one proves that $\bar{h}_{m_{\ell}}$ is negative with respect to $(\succ^j)_{h_i}$ \esp 
for all \esp $\ell \! > \! i$. \esp These two facts together obviously imply that the orders 
$(\preceq^j)_{h_i}$ are two-by-two different.

\vspace{0.43cm}

\noindent{\bf Proof of Proposition \ref{sandra}.} Let \esp $(\varepsilon_m)$ 
\esp be the decreasing sequence of positive numbers converging to 0 defined 
by (\ref{def-em}). With respect to this sequence we may perform the construction 
given in Claim 1. By Claim 2, for each \esp $m \!\in\! \mathbb{N}$ \esp 
we may then fix an element \esp $g_m \!\in\! \Gamma$ \esp such that, for each 
$j \!\in\! \{1,\ldots,2^k\}$, an element in $\mathcal{G}_m \setminus \con$ belongs 
to the positive cone of $(\preceq^j)_{g_m}$ if and only if it belongs to the positive cone 
of \esp $\preceq$. \esp Moreover, by Claim 3, the sequence $(g_m)$ may be taken in such a way 
that, for each fixed \esp $j \!\in\! \{1,\ldots,2^k\}$, \esp the orderings \esp $(\preceq^j)_{g_m}$ 
\esp are two-by-two different. Passing to a subsequence if necessary, Claim 2 allows 
to ensure that each sequence of orderings \esp $(\preceq^j)_{g_m}$ \esp converges to some 
ordering of the form $\preceq^{j'}$. Thus, \esp $\preceq^{j'}$ belongs to the set of accumulation 
points \esp \esp $\mathrm{acc}(\mathrm{orb}(\preceq^j))$ \esp \esp of the orbit of $\preceq^j$. 
Let us fix \esp $j_0 \!\in\! \{1,\ldots,2^k\}$. \esp By the above one has \esp 
$\preceq^{j_1} \in \mathrm{acc}(\mathrm{orb}(\preceq^{j_0}))$ \esp for some \esp 
$j_1 \!\in\! \{1,\ldots,2^k\}$. \esp If $j_0 \!=\! j_1$ then we are done. If not, 
then for a certain \esp $j_2 \!\in\! \{1,\ldots,2^k\}$ \esp one has \esp 
$\preceq^{j_2} \in \mathrm{acc}(\mathrm{orb}(\preceq^{j_1}))$, \esp and therefore 
\esp $\preceq^{j_2} \in \mathrm{acc}(\mathrm{orb}(\preceq^{j_0}))$. \esp If $j_2$ 
equals $j_0$ or $j_1$ then we are done. If not, we continue the process... Clearly,
in no more than $2^k$ steps we will find an index $j$ such that \esp 
$\preceq^{j} \esp \in \mathrm{acc}(\mathrm{orb}(\preceq^{j}))$, \esp 
and this concludes the proof. $\hfill\square$

\vspace{0.5cm}

Although very natural, our proof of Theorem C in the case of an ordering having a non-trivial 
Conradian soul with finitely many orders is quite elaborate. However, an affirmative answer 
to the following question would allow to reduce the general case to those of Propositions 
\ref{dorilita} and \ref{unilita}.

\begin{question} Let $\Gamma$ be a countable orderable group. If $\Gamma$ admits a non Conradian 
ordering, is it necessarily true that $\Gamma$ admits an ordering having trivial Conradian soul~?
\end{question}

%%%%%%%%%%%%%%%%%%%%%%%%%%%%%%%%%%%%%%%%%%%%%%%%%%%%%%%%%%%%%%%%%%%%%%%%%%%%%%%%%%%%%%%%%%%%%%%%%

\subsection{An application to braid groups}
\label{caso-cuerdas}

\hspace{0.35cm} For the proof of Theorem D we first consider the case of the braid 
group $B_3$. According to Examples \ref{dehor6}, \ref{dehor7}, and \ref{dehor8}, 
the Conradian soul of Dehornoy's ordering coincides with the cyclic subgroup 
generated by $\sigma_2$. Since this subgroup admits finitely many (namely, 
two) different orderings, we are under the hypothesis of Proposition 
\ref{sandra} for the orderings \esp $\preceq^1 = \preceq_D$ \esp and \esp 
$\preceq^2 = \preceq_{DD}$. \esp Now the conjugates of $\preceq_D$ cannot 
approximate $\preceq_{DD}$, because the latter ordering is isolated in $\mathcal{O}(B_3)$. 
Therefore, according to the proof of Proposition \ref{sandra}, there exists a sequence 
of elements \esp\esp $g_m \!\in\! B_3$ \esp\esp such that both sequences of orderings \esp 
$(\preceq_D)_{g_m}$ \esp and \esp $(\preceq_{DD})_{g_m}$ \esp converge to $\preceq_D$.

Now, for the case of general braid groups $B_n$, recall that the subgroup \esp 
$\langle \sigma_{n-2},\sigma_{n-1} \rangle$ \esp is isomorphic to $B_3$ via the map \esp 
$\sigma_{n-2} \mapsto \sigma_1$, \esp $\sigma_{n-1} \mapsto \sigma_2$, \esp which 
respects Dehornoy's orderings. By the argument above, there exists a sequence of elements 
$g_m$ in \esp $\langle \sigma_{n-2},\sigma_{n-1} \rangle$ \esp such that the restrictions 
to \esp $\langle \sigma_{n-2},\sigma_{n-1} \rangle$ \esp of the orderings $(\preceq_D)_{g_m}$ 
converge to the restriction of $\preceq_D$ to the same subgroup. We claim that actually 
\esp $(\preceq_D)_{g_m}$ \esp converges to $\preceq_D$ over the whole group $B_n$. 
Indeed, if $g$ belongs to $B_3 \setminus \langle \sigma_{n-2},\sigma_{n-1} \rangle$ 
and \esp $h \!\in\! B_n$ \esp is $\sigma_i$-positive (resp. $\sigma_i$-negative) for 
some \esp $i \!\in\! \{1,\ldots,n-3\}$, \esp then each of the elements $g_m h g_m^{-1}$ 
is still $\sigma_i$-positive (resp. $\sigma_i$-negative). Since the orderings \esp 
$(\preceq_D)_{g_m}$ \esp are two-by-two distinct, this finishes the proof of 
Theorem D.

\vspace{0.05cm}

\begin{rem} It would be interesting to obtain a proof of 
Theorem A using the methods of that of Theorem D.
\end{rem}

%%%%%%%%%%%%%%%%%%%%%%%%%%%%%%%%%%%%%%%%%%%%%%%%%%%%%%%%%%%%%%%%%%%%%%%%%%%%%%%%%%%%%%%%%%%%%%%%%

\begin{footnotesize}

\end{footnotesize}

\vspace{0.1cm}

\noindent Andr\'es Navas\\

\noindent Univ. de Santiago de Chile\\ 

\noindent Alameda 3363, Est. Central, Santiago, Chile\\ 

\noindent E-mail address: andres.navas@usach.cl\\

\end{document}